\documentclass[12pt,twoside]{article}

\usepackage{amssymb}
\font\teneufm=eufm10 scaled \magstep1
\font\seveneufm=eufm7 scaled \magstep1
\font\fiveeufm=eufm5  scaled \magstep1
\newfam\eufmfam
\textfont\eufmfam=\teneufm
\scriptfont\eufmfam=\seveneufm
\scriptscriptfont\eufmfam=\fiveeufm
\def\frak#1{{\fam\eufmfam\relax#1}}

\usepackage{mathrsfs}

\usepackage{amsmath}

\usepackage{ucs}
\usepackage[utf8x]{inputenc}

\usepackage[all,cmtip]{xy}

\newfam\msbfam
\font\tenmsb=msbm10 scaled \magstep1  \textfont\msbfam=\tenmsb
\font\sevenmsb=msbm7 scaled \magstep1 \scriptfont\msbfam=\sevenmsb
\font\fivemsb=msbm5 scaled \magstep1  \scriptscriptfont\msbfam=\fivemsb
\def\Bbb{\fam\msbfam \tenmsb}

\def\RR{{\Bbb R}}
\def\CC{{\Bbb C}}
\def\BB{{\Bbb B}}
\def\QQ{{\Bbb Q}}
\def\NN{{\Bbb N}}
\def\ZZ{{\Bbb Z}}
\def\TT{{\Bbb T}}
\def\PP{{\Bbb P}}
\def\HH{{\Bbb H}}

\def\ra{\rightarrow}

 \def\HollowBoxx #1#2#3{{\dimen0=#1 \advance\dimen0 by -#2
       \dimen1=#1 \advance\dimen1 by #3
        \vrule height 0pt depth #3 width #2
       \hskip -#3
       \vrule height #1 depth #3 width #3}}
 \def\LeftContraction{\mathord{\kern1.45pt \HollowBoxx{6pt}{3.5pt}{.4pt}}\,}

 \def\HollowBox #1#2#3{{\dimen0=#1 \advance\dimen0 by -#3
       \dimen1=#1 \advance\dimen1 by #3
        \vrule height #1 depth #3 width #3
        \vrule height 0pt depth #3 width #2
        \hskip -#3}}
 \def\RightContraction{\mathord{\, \HollowBox{6pt}{3.1pt}{.4pt}} \kern1.6pt}

\def\qed{{\hfill $\Box$}}
\newtheorem{theorem}{THEOREM}[section]

\newtheorem{lemma}[theorem]{Lemma}

\newtheorem{remark}[theorem]{Remark}
\newtheorem{proposition}[theorem]{Proposition}

\begin{document}

\begin{center}
{\Large \bf Proper Actions of Lie Groups
\vspace{0.1cm}\\
of Dimension $n^2+1$
\vspace{0.4cm}\\
on $n$-Dimensional Complex Manifolds}
\medskip \\
\normalsize A. V. Isaev and N. G. Kruzhilin
\end{center}

\begin{quotation} \small \sl In this paper we continue to study actions of high-dimensional Lie groups on complex manifolds. We give a complete explicit description of all pairs $(M,G)$, where $M$ is a connected complex manifold $M$ of dimension $n\ge 2$, and $G$ is a connected Lie group of dimension $n^2+1$ acting effectively and properly on $M$ by holomorphic transformations. This result complements a classification obtained earlier by the first author for $n^2+2\le\hbox{dim}\,G<n^2+2n$ and a classical result due to W. Kaup for the maximal group dimension $n^2+2n$.
\end{quotation}

\thispagestyle{empty}

\pagestyle{myheadings}
\markboth{A. V. Isaev and N. G. Kruzhilin}{Proper Actions on Complex Manifolds}

\setcounter{section}{-1}

\section{Introduction}
\setcounter{equation}{0}

Let $M$ be a connected $C^{\infty}$-smooth manifold and $\hbox{Diff}(M)$ the group of $C^{\infty}$-smooth diffeomorphisms of $M$ endowed with the compact-open topology. A topological group $G$ is said to act continuously on $M$ by diffeomorphisms, if a continuous homomorphism $\Phi: G\ra\hbox{Diff}(M)$ is specified. The continuity of $\Phi$ is equivalent to the continuity of the action map
$$
\hat\Phi:\,G\times M\ra M,\quad (g,p)\mapsto \Phi(g)(p)=:gp.
$$ 
We only consider effective actions, that is, assume that the kernel of $\Phi$ is trivial.

The action of $G$ on $M$ is called {\it proper}, if the map
$$
\Psi:\, G\times M\ra M\times M, \quad (g,p)\mapsto (gp,p),
$$
is proper, i.e. for every compact subset $C\subset M\times M$ its inverse image $\Psi^{-1}(C)\subset G\times M$ is compact as well. For example, the action is proper if $G$ is compact. The properness of the action implies that: (i) $G$ is locally compact, hence by \cite{BM1}, \cite{BM2} (see also \cite{MZ}) it carries the structure of a Lie group and the action map $\hat\Phi$ is smooth; (ii) $\Phi$ is a topological group isomorphism between $G$ and $\Phi(G)$; (iii) $\Phi(G)$ is a closed subgroup of $\hbox{Diff}(M)$ (see \cite{Bi} for a brief survey on proper actions). Thus, one can assume that $G$ is a Lie group acting smoothly and properly on the manifold $M$, and that it is realized as a closed subgroup of $\hbox{Diff}(M)$.

Suppose now that $M$ is equipped with a Riemannian metric ${\mathscr G}$, and let $\hbox{Isom}(M,{\mathscr G})$ be the group of all isometries of $M$ with respect to ${\mathscr G}$. It was shown in \cite{MS} that $\hbox{Isom}(M,{\mathscr G})$ acts properly on $M$ (and so does its every closed subgroup). Conversely, by \cite{Pal} (see also \cite{Al}), for any Lie group acting properly on $M$ there exists a $C^{\infty}$-smooth $G$-invariant metric on $M$. It then follows that Lie groups acting properly and effectively on the manifold $M$ by diffeomorphisms are precisely closed subgroups of $\hbox{Isom}(M,{\mathscr G})$ for all possible smooth Riemannian metrics ${\mathscr G}$ on $M$.

If $G$ acts properly on $M$, then for every $p\in M$ its isotropy subgroup
$$
G_p:=\left\{g\in G: gp=p\right\}
$$
is compact in $G$. Then by \cite{Bo} the isotropy representation
$$
\alpha_p:\, G_p\ra GL(\RR,T_p(M)),\quad g\mapsto dg(p)
$$
is continuous and faithful, where $T_p(M)$ denotes the tangent space to $M$ at $p$ and $dg(p)$ is the differential of $g$ at $p$. In particular, the linear isotropy subgroup
$$
LG_p:=\alpha_p(G_p)
$$
is a compact subgroup of $GL(\RR,T_p(M))$ isomorphic to $G_p$. In some coordinates in $T_p(M)$ the group $LG_p$ becomes a subgroup of the orthogonal group $O_m(\RR)$, where $m:=\hbox{dim}\,M$. Hence $\hbox{dim}\,G_p\le\hbox{dim}\, O_m(\RR)=m(m-1)/2$. Furthermore, for every $p\in M$ its orbit
$$
Gp:=\left\{gp:g\in G\right\}
$$
is a closed submanifold of $M$, and $\hbox{dim}\, Gp\le m$. Thus, setting $d_G:=\hbox{dim}\,G$, we obtain
$$
d_G=\hbox{dim}\, G_p+\hbox{dim}\,Gp\le m(m+1)/2.
$$

It is a classical result (see \cite{F}, \cite{Ca}, \cite{Ei}) that if $G$ acts properly on a smooth manifold $M$ of dimension $m\ge 2$ and $d_G=m(m+1)/2$, then $M$ is isometric (with respect to some $G$-invariant metric) either to one of the standard complete simply-connected spaces of constant sectional curvature $\RR^m$, $S^m$, $\HH^m$ (where $\HH^m$ is the hyperbolic space), or to $\RR\PP^m$. Next, it was shown in \cite{Wa} (see also \cite{Eg}, \cite{Y1}) that a group $G$ with $m(m-1)/2+1<d_G<m(m+1)/2$ cannot act properly on a smooth manifold $M$ of dimension $m\ne 4$. The exceptional 4-dimensional case was considered in \cite{Ish}, and it turned out that a group of dimension 9 cannot act properly on a 4-dimensional manifold.
Furthermore, it was shown in \cite{Ish} that if a 4-dimensional manifold admits a proper action of an 8-dimensional group $G$, then it has a $G$-invariant complex structure;  this is the case  of complex $n$-dimensional manifolds (for $n=2$) admitting proper actions of $(n^2+2n)$-dimensional groups that will be discussed later.

There exists also an explicit classification of pairs $(M,G)$, where $m\ge 4$, $G$ is connected, and  $d_G=m(m-1)/2+1$ (see \cite{Y1}, \cite{Ku}, \cite{Ob}, \cite{Ish}). Further, in \cite{KN} a reasonably explicit classification of pairs $(M,G)$, where $m\ge 6$, $G$ is connected, and  $(m-1)(m-2)/2+2\le d_G\le m(m-1)/2$, was given. We also mention a classification of $G$-homogeneous manifolds for $m=4$, $d_G=6$ (see \cite{Ish}) and a classifications of $G$-homogeneous simply-connected manifolds in the cases $m=3$, $d_G=3,4$ and $m=4$, $d_G=5$ (see \cite{Ca}, \cite{Pat}) obtained by E. Cartan's method of adapted frames introduced in \cite{Ca}. There are many other results, especially for compact subgroups, but -- to the best of our knowledge -- no complete classifications exist beyond dimension $(m-1)(m-2)/2+2$ (see \cite{Ko2}, \cite{Y2} and references therein for further details).

We study proper group actions in the complex setting. From now on, $M$ will denote a complex manifold of complex dimension $n$ (hence $m=2n$)  and $G$ will be a subgroup of $\hbox{Aut}(M)$, the group of all holomorphic automorphisms of $M$. We will be classifying pairs $(M,G)$, but we will not be concerned with determining $G$-invariant Riemannian metrics on $M$. Proper actions by holomorphic transformations are found in abundance. A fundamental result due to Kaup (see \cite{Ka}) states that every closed subgroup of $\hbox{Aut}(M)$ that preserves a continuous distance on $M$ acts properly on $M$. Thus, Lie groups acting properly and effectively on $M$ by holomorphic transformations are precisely those closed subgroups of $\hbox{Aut}(M)$ that preserve continuous distances on $M$. In particular, if $M$ is a Kobayashi-hyperbolic manifold, then $\hbox{Aut}(M)$ is a Lie group acting properly on $M$ (see also \cite{Ko1}).

In the complex setting, in some coordinates in $T_p(M)$ the group $LG_p$ becomes a subgroup of the unitary group $U_n$. Hence $\hbox{dim}\,G_p\le\hbox{dim}\, U_n=n^2$, and therefore
$$
d_G\le n^2+2n.
$$
We note that $n^2+2n<(m-1)(m-2)/2+2$ for $m=2n$ and $n\ge 5$. Thus, the group dimension range that arises in the complex case, for $n\ge 5$ lies strictly below the dimension range considered in the classical real case and therefore is not covered by the existing results. Furthermore, overlaps with these results for $n=3,4$ and $n=2$, $d_G=6$ occur only in relatively easy situations and do not lead to any significant simplifications in the complex case. The only interesting overlap with the real case occurs for $n=2$, $d_G=5$ (see \cite{Pat}); we will briefly discuss it below. Note that in the situations when overlaps do occur, the existing classifications in the real case do not necessarily immediately lead to classifications in the complex case, since the determination of all $G$-invariant complex structures on the corresponding real manifolds may be a non-trivial task.

It was shown by Kaup in \cite{Ka} that if $d_G=n^2+2n$, then $M$ is holomorphically equivalent (in fact, holomorphically isometric with respect to some $G$-invariant metric) to one of $\BB^n:=\left\{z\in\CC^n: |z|<1\right\}$, $\CC^n$, $\CC\PP^n$, and an equivalence map $F$ can be chosen so that the group\linebreak $F\circ G\circ F^{-1}:=\left\{F\circ g\circ F^{-1}:g\in G\right\}$ is, respectively, one of the groups $\hbox{Aut}(\BB^n)$, $G(\CC^n)$, $G(\CC\PP^n)$. Here $\hbox{Aut}(\BB^n)\simeq PSU_{n,1}:=SU_{n,1}/\hbox{(center)}$ is the group of all transformations 
$$
z\mapsto\displaystyle\frac{Az+b}{cz+d},\\
$$
where
$$
\left(\begin{array}{cc}
A& b\\
c& d
\end{array}
\right)
\in SU_{n,1};
$$
$G\left(\CC^n\right)\simeq U_n\ltimes\CC^n$ is the group of all holomorphic automorphisms of $\CC^n$ of the form
\begin{equation}
z\mapsto Uz+a,\label{groupcn}
\end{equation}
where $U\in U_n$, $a\in\CC^n$ (we usually write $G\left(\CC\right)$ instead of $G\left(\CC^1\right)$); and $G\left(\CC\PP^n\right)\simeq PSU_{n+1}:=SU_{n+1}/\hbox{(center)}$ is the group of all holomorphic automorphisms of $\CC\PP^n$ of the form
\begin{equation}
\zeta\mapsto U\zeta,\label{groupcpn}
\end{equation}
where $\zeta$ is a point in $\CC\PP^n$ written in homogeneous coordinates, and $U\in SU_{n+1}$ (this group is a maximal compact subgroup of the complex Lie group $\hbox{Aut}(\CC\PP^n)\simeq PSL_{n+1}(\CC):=SL_{n+1}(\CC)/\hbox{(center)}$). We remark that the groups $\hbox{Aut}(\BB^n)$, $G(\CC^n)$, $G(\CC\PP^n)$ are the full groups of holomorphic isometries of the Bergman metric on $\BB^n$, the flat metric on $\CC^n$, and the Fubini-Study metric on $\CC\PP^n$, respectively, and that the above result due to Kaup can be obtained directly from the classification of Hermitian symmetric spaces (cf. \cite{Ak}, pp. 49--50). 

In the above situation we say for brevity that $F$ {\it transforms}\, $G$ into one of $\hbox{Aut}(\BB^n)$, $G(\CC^n)$, $G(\CC\PP^n)$, respectively, and, in general, if $F: M_1\ra M_2$ is a biholomorphic map, $G_j\subset\hbox{Aut}(M_j)$, $j=1,2$, are subgroups and $F\circ G_1\circ F^{-1}=G_2$, we say that $F$ transforms $G_1$ into $G_2$. If for two pairs $(M_1,G_1)$ and $(M_2,G_2)$ one can find such a map, we say that the pairs are {\it equivalent}.

We are interested in characterizing pairs $(M,G)$ up to this equivalence relation for $d_G<n^2+2n$, where $G\subset\hbox{Aut}(M)$ is connected and acts on $M$ properly. In \cite{IKra}, \cite{I5} a complete classification was obtained for $n^2+2\le d_G<n^2+2n$. Furthermore, in \cite{IKra}, \cite{I1}, \cite{I2}, \cite{I3} we considered the special case where $M$ is a Kobayashi-hyperbolic manifold and $G=\hbox{Aut}(M)$, and determined all manifolds with $n^2-1\le d_{\hbox{\small Aut}(M)}< n^2+2n$, $n\ge 2$ (see \cite{I4} for a comprehensive exposition of these results). 

In the present paper we assume that $d_G=n^2+1$. This is the next group dimension down from the range $n^2+2\le d_G<n^2+2n$ considered in \cite{IKra}, \cite{I5}. Note that $n^2+1$ is the lowest group dimension for which proper actions are necessarily transitive (it is shown in \cite{Ka} that $G$-homogeneity always takes place for $d_G>n^2$); indeed, for $d_G=n^2$ both homogeneous and non-homogeneous manifold occur (see \cite{I4}). For $d_G=n^2+1$ we have $\hbox{dim}\, G_p=(n-1)^2$, and we start by describing connected subgroups of the unitary group $U_n$ of dimension $(n-1)^2$ in Proposition \ref{un} (see Section \ref{classtab}), thus determining the connected identity components of all possible linear isotropy subgroups. According to this description, every action falls into one of three types. In Sections \ref{type1}, \ref{type2} we deal with actions of type I and II, respectively, and obtain complete lists of the corresponding pairs $(M,G)$ in Theorems \ref{classtype1} and \ref{classtype2}. Actions of type III are more difficult to deal with. In Section \ref{type3} we give a large number of examples of such actions and show that these examples provide a complete description of actions of type III (see Theorem \ref{classtype3}). This is the main result of the paper. Taken together, Theorems \ref{classtype1}, \ref{classtype2}, \ref{classtype3} describe all pairs $(M,G)$ with $d_G=n^2+1$.

Regarding Theorems \ref{classtype1}, \ref{classtype2}, \ref{classtype3} for $n=2,3$ some remarks are in order. Firstly, all connected 2- and 3-dimensional complex manifolds that admit transitive actions of Lie  groups by holomorphic transformations were determined in \cite{HL}, \cite{OR}, \cite{Wi}. However, it was not the aim of those articles to give a description of {\it all}\, possible transitive actions, and, indeed, most actions listed in Theorems \ref{classtype1}, \ref{classtype2}, \ref{classtype3} do not occur there. Hence our classification for $n=2,3$ does not follow from \cite{HL}, \cite{OR}, \cite{Wi}. Secondly, as we have already mentioned, a classification of all effective proper transitive actions of connected 5-dimensional Lie groups on simply-connected real 4-dimensional manifolds was given in \cite{Pat} (see also \cite{Ish}). Therefore, one can attempt to obtain Theorem \ref{classtype3} for $n=2$ by, firstly, determining all invariant complex structures on the real manifolds occurring in \cite{Pat} and, secondly, by passing to their quotients to produce a complete list of manifolds from that of simply-connected ones. Finally, we have been informed by G. Fels that he has recently obtained Theorem \ref{classtype3} for $n=2$ by an alternative method. 
\vspace{0.5cm}

\noindent {\bf Acknowledgements.} Part of this work was done while both authors were visiting the Ruhr-Universit\"at Bochum in January-February 2007, with the first author's visit partially supported by the Alexander von Humboldt-Stiftung. Also, a large amount of work on the paper was done during the first author's visit to the Max-Plank Institut f\"ur Mathematik in Bonn in April-May 2007. The research of the second author was supported by the RFBR (grant no. 05-01-0981a) and the Program of Support of Scientific Schools of the RF. We would like to thank G. Fels for making a large number of useful comments and interest in our work.

\section{Classification of Linear Isotropy Subgroups}\label{classtab}
\setcounter{equation}{0}

In this section we prove the following proposition that extends Lemma 2.1 of \cite{IKru}.

\begin{proposition}\label{un} \sl Let $H$ be a connected closed subgroup of $U_n$ of dimension $(n-1)^2$, $n\ge 2$. Then $H$ is conjugate in $U_n$ to one of the following subgroups:
\vspace{0.3cm}

\noindent I. $e^{i\RR}SO_3(\RR)$ (here $n=3$); 
\vspace{0.3cm}

\noindent II. $SU_{n-1}\times U_1$ realized as the subgroup of all matrices
\begin{equation}
\left(\begin{array}{cc}
A & 0\\
0& e^{i\theta}
\end{array}\right),\label{spec}
\end{equation}
where $A \in SU_{n-1}$ and $\theta\in\RR$, for $n\ge 3$; 
\vspace{0.3cm}

\noindent III. the subgroup $H_{k_1,k_2}^n$ of all matrices
\begin{equation}
\left(\begin{array}{cc}
A & 0\\
0 & a
\end{array}\right),\label{mat}
\end{equation}
where $k_1,k_2$ are fixed integers such that $(k_1,k_2)=1$, $k_1>0$, and
$A\in U_{n-1}$, $a\in
(\det A)^{\frac{k_2}{k_1}}:=\exp(k_2/k_1\, \hbox{Ln}\,(\det A))$.\footnote{For $k_2\ne 0$ the group $H_{k_1,k_2}^n$ is a $k_1$-sheeted cover of $U_{n-1}$.} 
\end{proposition}

\noindent {\bf Proof:} Since $H$ is compact, it is completely
reducible, i.e. $\CC^n$ splits into the sum of $H$-invariant
pairwise orthogonal complex subspaces, $\CC^n=V_1\oplus\dots\oplus V_m$,
such that the restriction $H_j$ of $H$ to each  $V_j$ is irreducible. Let
$n_j:=\hbox{dim}_{\CC}V_j$ (hence $n_1+\dots+n_m=n$) and let
$U_{n_j}$ be the group of unitary
transformations of $V_j$. Clearly, $H_j\subset U_{n_j}$, and therefore
$\hbox{dim}\,H\le n_1^2+\dots+n_m^2$.
On the other hand  $\hbox{dim}\,H=(n-1)^2$, which shows that
$m\le 2$.

Let $m=2$. Then there exists a unitary change of coordinates
in $\CC^n$ such all elements of $H$ take the form (\ref{mat}), where $A \in U_{n-1}$ and $a \in U_1$. We note that the groups $H_1$, $H_2$ consist of all possible $A$ and $a$, respectively.   

If $\hbox{dim}\,H_2=0$, then $H_2=\{1\}$, and
therefore $H_1=U_{n-1}$. In this case we obtain the group $H_{1,0}^n$.

Assume  that $\hbox{dim}\,H_2=1$, i.e., $H_2=U_1$. Then
$(n-1)^2-1\le\hbox{dim}\,H_1\le (n-1)^2$.  Let
$\hbox{dim}\,H_1=(n-1)^2-1$ first. The only connected subgroup of $U_{n-1}$
of dimension $(n-1)^2-1$ is $SU_{n-1}$. Hence $H$ is conjugate to
the subgroup
of matrices of the form (\ref{spec}) if $n\ge 3$ and to $H_{1,0}^2$ for $n=2$. Now let  $\hbox{dim}\,H_1=(n-1)^2$, i.e., $H_1=U_{n-1}$. Consider the Lie algebra ${\frak h}$ of $H$. Up to conjugation, it consists of matrices of the form
\begin{equation}
\left(\begin{array}{cc}
{\frak A}& 0\\
0& l({\frak A})
\end{array}\right),\label{mat2}
\end{equation}
where ${\frak A}\in{\frak u}_{n-1}$ and $l({\frak A})\not\equiv 0$
is a linear function  of the matrix elements of ${\frak A}$ ranging  in $i\RR$. Clearly, $l({\frak A})$ must vanish on the derived algebra of ${\frak u}_{n-1}$, which is ${\frak{su}}_{n-1}$. Hence matrices (\ref{mat2}) form  a Lie algebra if and only if $l({\frak A})=c\cdot\hbox{trace}\,{\frak A}$, where $c\in\RR\setminus\{0\}$. Such an algebra can be the Lie algebra of a closed subgroup of $U_{n-1}\times U_1$  only if $c\in\QQ\setminus\{0\}$. Hence $H$ is conjugate to $H_{k_1,k_2}^n$ for some $k_1,k_2\in\ZZ$, where one can always assume that $k_1>0$ and $(k_1,k_2)=1$.

Now let $m=1$. We  shall proceed as in the proof of Lemma 2.1 in \cite{IKra}.
Let ${\frak h}^{\CC}:={\frak h}+i{\frak h}\subset{\frak{gl}}_n$ be the complexification of ${\frak h}$, where ${\frak{gl}}_n:={\frak{gl}}_n(\CC)$. The algebra ${\frak h}^{\CC}$ acts irreducibly on $\CC^n$ and by a theorem of E. Cartan (see, e.g., \cite{GG}), ${\frak h}^{\CC}$ is either semisimple or  the direct sum of the center ${\frak c}$ of ${\frak{gl}}_n$ and a semisimple ideal ${\frak t}$. Clearly, the action of the ideal ${\frak t}$ on $\CC^n$ is irreducible.

Assume first that ${\frak h}^{\CC}$ is semisimple, and let ${\frak
h}^{\CC}={\frak h}_1\oplus\dots\oplus{\frak h}_k$ be its decomposition
into the direct sum of simple ideals. Then the natural irreducible $n$-dimensional representation of
${\frak h}^{\CC}$ (given by the embedding of ${\frak h}^{\CC}$ in ${\frak{gl}}_n$) is the tensor product of some irreducible faithful representations of the ${\frak h}_j$ (see,
e.g., \cite{GG}). Let $n_j$ be the dimension of the corresponding representation  of ${\frak h}_j$, $j=1,\dots,k$. Then $n_j\ge 2$, $\hbox{dim}_{\CC}\,{\frak h}_j\le n_j^2-1$, and $n=n_1\cdot...\cdot n_k$. The following observation is simple.

\begin{quote}
{\bf Claim:} {\sl If $n=n_1\cdot...\cdot n_k$, $k\ge 2$, $n_j\ge 2$
for $j=1,\dots,k$, then $\sum_{j=1}^k n_j^2\le n^2-2n$.}
\end{quote}

Since $\hbox{dim}_{\CC}\,{\frak h}^{\CC}=(n-1)^2$, it follows from the
above claim that $k=1$, i.e. ${\frak h}^{\CC}$ is
simple. The minimal dimensions of irreducible faithful
representations of complex simple Lie algebras ${\frak s}$ are well-known (see, e.g., \cite{OV}). In the table below $V$ denotes representations of minimal dimension.

\begin{center}
\begin{tabular}{|l|c|c|}
\hline
\multicolumn{1}{|c|}{${\frak s}$}&
\multicolumn{1}{c|}{ $\hbox{dim}\,V$}&
\multicolumn{1}{c|}{$\hbox{dim}\,{\frak s}$}
\\ \hline
${\frak{sl}}_k$\,\,$k\ge 2$ & $k$ & $k^2-1$
\\ \hline
${\frak o}_k$\,\, $k\ge 7$&  $k$ & $k(k-1)/2$
\\ \hline
${\frak{sp}}_{2k}$\,\,$k\ge 2$ & $2k$ & $2k^2+k$
\\ \hline
${\frak e}_6$ & 27 & 78
\\ \hline
${\frak e}_7$ & 56 & 133
\\ \hline
${\frak e}_8$ & 248 & 248
\\ \hline
${\frak f}_4$ & 26 & 52
\\ \hline
${\frak g}_2$ & 7 & 14
\\ \hline
\end{tabular}
\end{center}

Since $\hbox{dim}_{\CC}\,{\frak h}^{\CC}=(n-1)^2$, it follows that
none  of the above possibilities realize. Therefore, ${\frak h}^{\CC}={\frak c}\oplus{\frak t}$, where $\hbox{dim}\,{\frak t}=n^2-2n$. Then, if $n=2$, we obtain that $H$ coincides with the center of $U_2$ which is impossible since its action on $\CC^2$ is then not irreducible. Assuming that $n\ge 3$ and repeating the above argument for ${\frak t}$ in place of ${\frak h}^{\CC}$, we see that ${\frak t}$ can only be isomorphic to ${\frak{sl}}_{n-1}$. But ${\frak{sl}}_{n-1}$ does not have an irreducible $n$-dimensional representation unless $n=3$. 

Thus, $n=3$ and ${\frak h}^{\CC}\simeq\CC\oplus{\frak{sl}}_2\simeq\CC\oplus{\frak{so}}_3$. Further, we observe that every irreducible $3$-dimensional representation of ${\frak{so}}_3$ is equivalent to its defining representation. This implies that $H$ is conjugate in $GL_3(\CC)$ to $e^{i\RR}SO_3(\RR)$. Since $H\subset U_3$ it is straightforward to show that the conjugating element can be chosen to belong to $U_3$.

The proof of the proposition is complete.\qed
\smallskip\\

Let $M$ be a connected complex manifold of dimension $n\ge 2$, and suppose that a connected Lie group $G\subset\hbox{Aut}(M)$ with $d_G=n^2+1$ acts properly on $M$. Fix $p\in M$, consider the linear isotropy subgroup $LG_p$, and choose coordinates in $T_p(M)$ so that $LG_p\subset U_n$. We say that the pair $(M,G)$ (or the action of $G$ on $M$) is of type I, II or III, if the connected identity component $LG_p^0$ of the group $LG_p$ is conjugate in $U_n$ to a subgroup listed in I, II or III of Proposition \ref{un}, respectively. Since $M$ is $G$-homogeneous, this definition is independent of the choice of $p$.

We will now separately consider actions of each type. Recall from the introduction that $G(\CC^n)$ is the group of holomorphic isometries of $\CC^n$ with respect to the Euclidean metric (see (\ref{groupcn})) and $G(\CC\PP^n)$ is the group of holomorphic isometries of $\CC\PP^n$ with respect to the Fubini-Study metric (see (\ref{groupcpn})). These groups and some of their subgroups will frequently occur throughout the paper. In particular, we will need the subgroup $G_1(\CC^n)\simeq SU_n\ltimes\CC^n$ that consists of all elements of $G(\CC^n)$ with $U\in SU_n$ and the subgroup $G_2(\CC^3)\simeq e^{i\RR}SO_3(\RR)\ltimes\CC^3$ that consists of all elements of $G(\CC^3)$ with $U\in e^{i\RR}SO_3(\RR)$. We usually write $G\left(\CC\right)$ instead of $G\left(\CC^1\right)$ and $G_1\left(\CC\right)$ instead of $G_1\left(\CC^1\right)$.      

\section{Actions of Type I}\label{type1}
\setcounter{equation}{0}

In this section we prove the following theorem.

\begin{theorem}\label{classtype1}{\sl Let $M$ be a connected complex manifold of dimension 3 and $G\subset\hbox{Aut}(M)$ a connected Lie group with $d_G=10$ that acts properly on $M$. If the pair $(M,G)$ is of type I, then it is equivalent to one of the following:
\vspace{0.3cm}

\noindent (i) $(\mathscr{S},\hbox{Aut}(\mathscr{S}))$, where $\mathscr{S}$ is the Siegel space
$$
\mathscr{S}:=\left\{(z_1,z_2,z_3)\in\CC^3: Z\overline{Z}\ll\hbox{id}\right\},
$$
with
$$
Z:=\left(
\begin{array}{ll}
z_1 & z_2\\
z_2 & z_3
\end{array}
\right)
$$
(here $\hbox{Aut}(\mathscr{S})$ is isomorphic to $Sp_4(\RR)/\ZZ_2$);
\vspace{0.3cm}

\noindent (ii) $({\cal Q}_3,SO_5(\RR))$, where ${\cal Q}_3$ is the complex quadric in $\CC\PP^4$, and $SO_5(\RR)$ is realized as a maximal compact subgroup of $\hbox{Aut}({\cal Q}_3)^0\simeq PSO_5(\CC)$;
\vspace{0.3cm}

\noindent (iii) $(\CC^3,G_2(\CC^3))$.}
\end{theorem}

\noindent{\bf Proof:} Fix a $G$-invariant Hermitian metric on $M$. Since $LG_q$ for every $q\in M$ contains the element $-\hbox{id}$, the manifold $M$ equipped with this metric becomes a Hermitian symmetric space. The group $LG_p^0$ acts irreducibly on $T_p(M)$, and therefore $M$ either is an irreducible Hermitian symmetric space, or is equivalent (holomorphically and isometrically) to $\CC^3$ with the flat metric.

If $M$ is an irreducible Hermitian symmetric space, it follows from the general theory of Riemannian symmetric spaces that $G$ coincides with the connected identity component of the group of holomorphic isometries of $M$ (see Theorem 1.1 in Chapter V of \cite{H}). Now E. Cartan's classification of irreducible Hermitian symmetric spaces implies that $(M,G)$ is equivalent to either $(\mathscr{S},\hbox{Aut}(\mathscr{S}))$ or $({\cal Q}_3,SO_5(\RR))$ (see Chapter IX of \cite{H}). 

Let $M$ be equivalent to $\CC^3$ and let $F$ be an equivalence map. The map $F$ transforms $G$ into a closed subgroup of $G(\CC^3)$ (recall that $G(\CC^3)$ is the full group of holomorphic isometries of $\CC^3$ with respect to the flat metric). Let $p_0\in M$ be such that $F(p_0)=0$. Then $F$ transforms $G_{p_0}^0$ into a closed subgroup $H$ of $U_3\subset G(\CC^3)$ isomorphic to $e^{i\RR}SO_3(\RR)$ and acting irreducibly on $T_0(\CC^3)$. By Proposition \ref{un}, the subgroup $H$ is conjugate in $U_3$ to the standard embedding of $e^{i\RR}SO_3(\RR)$ in $U_3$, and hence there exists an equivalence map $\hat F$ between $M$ and $\CC^3$ that transforms $G_{p_0}^0$ into $e^{i\RR}SO_3(\RR)$.

Let ${\frak g}$ be the Lie algebra (isomorphic to the Lie algebra of $G$) of fundamental vector fields of the action of the group $\hat G:=\hat F\circ G\circ \hat F^{-1}$ on $\CC^3$, that is, ${\frak g}$ consists of all holomorphic vector fields $X$ on $\CC^3$ for which there exists an element $a$ of the Lie algebra of $\hat G$ such that for all $z\in\CC^3$ we have
$$
X(z)=\frac{d}{dt}\Bigl[\exp(ta)(z)\Bigr]\Bigr|_{t=0}.
$$
Since $\hat G\subset G(\CC^3)$, the algebra ${\frak g}$ is generated by $\langle Z_0\rangle\oplus {\frak{so}}_3(\RR)$ and some affine holomorphic vector fields $V_j$, $j=1,\dots,6$, that do not vanish at the origin. Here
$$
Z_0:=i\sum_{k=1}^3 z_k\,\partial/\partial z_k,
$$ 
and ${\frak{so}}_3(\RR)$ is realized as the algebra of fundamental vector fields of the standard action of $SO_3(\RR)$ on $\CC^3$. Considering $[Z_0,[V_j,Z_0]]$ instead of $V_j$, we can assume that $V_j$ are constant vector fields for all $j$ (cf. the proof of Satz 4.9 in \cite{Ka}). It then follows that $\hat G=G_2(\CC^3)$. 

The proof is complete.\qed

\section{Actions of Type II}\label{type2}
\setcounter{equation}{0}

In this section we obtain the following result.

\begin{theorem}\label{classtype2}{\sl Let $M$ be a connected complex manifold of dimension $n\ge 3$ and $G\subset\hbox{Aut}(M)$ a connected Lie group with $d_G=n^2+1$ that acts properly on $M$. If the pair $(M,G)$ is of type II, then it is equivalent to\linebreak $(\CC^{n-1}\times M',G_1(\CC^{n-1})\times G')$, where $M'$ is one of $\BB^1$, $\CC$, $\CC\PP^1$, and $G'$ is one of the groups $\hbox{Aut}(\BB^1)$, $G(\CC)$, $G(\CC\PP^1)$, respectively.}
\end{theorem}

We start with the following lemma that clarifies the structure of full linear isotropy subgroups for actions of type II.

\begin{lemma}\label{fullisotropyii}\sl Let $H\subset U_n$, with $n\ge 3$, be a closed subgroup, and let $H^0=SU_{n-1}\times U_1$. Then for some $m\in\NN$ the group $H$ consists of all matrices of the form
$$
\left(\begin{array}{cc}
\alpha A & 0\\
0 & a
\end{array}\right),
$$
where $A\in SU_{n-1}$, $a\in U_1$, $\alpha^m=1$.
\end{lemma}

\noindent{\bf Proof:} The proof is similar to that of Lemma 4.4 in \cite{IKru}. Let $C_1,\dots,C_K$ be the connected   components of $H$ with $C_1=H^0=SU_{n-1}\times U_1$. Clearly, there exist $g_1=\hbox{id},g_2,\dots,g_K$ in $U_n$
such that
$C_j=g_jH^0$, $j=1,\dots,K$. Moreover,  for each  pair of indices $i,j$ there exists an index $l$ such that $g_iH^0\cdot
g_jH^0=g_lH^0$, and therefore
\begin{equation}
g_l^{-1}g_iH^0g_j=H^0.\label{compii}
\end{equation}
Applying each side of (\ref{compii}) to the vector $v:=(0,\dots,0,1)$, which is an eigenvector of every element of $H^0$, we obtain that for every $h\in H^0$ there exists $\beta(h)\in\CC$ such that
$$
g_l^{-1}g_ihg_jv=\beta(h)v,
$$
or, equivalently,
$$
hg_jv=\beta(h)g_i^{-1}g_lv,
$$
which implies that $g_jv=(0,\dots,0, a_j)$, where $|a_j|=1$, $j=1,\dots,K$. Hence $g_j$ has the form
$$
g_j=\left(\begin{array}{cc}
A_j & 0\\
0 & a_j
\end{array}\right),
$$
where $A_j\in U_{n-1}$. Multiplying $g_j$ by an appropriate element of $H^0$, we can assume that $a_j=1$ and $A_j=\alpha_j\cdot\hbox{id}$, with $|\alpha_j|=1$, $j=1,\dots,K$.

Clearly, all elements in $H$ of the form
\begin{equation}
\left(\begin{array}{cc}
t\cdot\hbox{id} & 0\\
0 & 1
\end{array}\right),\label{formm3}
\end{equation}
where $|t|=1$ form a finite subgroup and therefore the corresponding numbers $t$ form a group of roots of unity of some order $m$.

The proof of the lemma is complete.\qed
\smallskip\\

\noindent{\bf Proof of Theorem \ref{classtype2}:} Fix $p\in M$, set $H:=G_p$, and identify $M$ as a smooth manifold with $G/H$ (in particular, we identify $T_p(M)$ with $T_H(G/H)$). By means of this identification we introduce a $G$-invariant complex structure on $G/H$. Let $\Pi_{G,H}: G\ra G/H$ be the factorization map. For every element $g\in G$ we denote by ${\mathscr L}_{g}$ the action of $g$ on $G/H$. Let ${\frak g}$ be the Lie algebra of $G$. Since the subgroup $\hbox{Ad}(H)\subset GL(\RR,{\frak g})$ is compact, there exists an $\hbox{Ad}(H)$-invariant scalar product on ${\frak g}$. Let ${\frak h}\subset{\frak g}$ be the Lie algebra of $H$ and ${\frak h}^{\perp}$ the orthogonal complement to ${\frak h}$ in ${\frak g}$. Since ${\frak h}$ is $\hbox{Ad}(H)$-invariant, so is ${\frak h}^{\perp}$. The map $\Phi:=d\Pi_{G,H}(\hbox{id})|_{{\frak h}^{\perp}}$ is a linear isomorphism between ${\frak h}^{\perp}$ and $T_H(G/H)$, and for every $h\in H$ transforms the operator $\hbox{Ad}(h)$ on ${\frak h}^{\perp}$ into the operator $d{\mathscr L}_h(H)$ on $T_H(G/H)$ (see e.g. \cite{On}). Analogously, for $c\in{\frak h}$, the map $\hbox{ad}(c)$ on ${\frak h}^{\perp}$ is transformed into ${\mathscr L}(c)$ on $T_H(G/H)$, where ${\mathscr L}$ is the differential of the map $h\mapsto d{\mathscr L}_h(H)$ at $\hbox{id}\in H$.

By Lemma \ref{fullisotropyii}, at every $q\in M$ there are exactly two non-trivial proper $LG_q$-invariant complex subspaces ${\cal L}_1(q)$ and ${\cal L}_2(q)$ in $T_q(M)$. Here ${\cal L}_1(q)$ denotes the invariant subspace of dimension $n-1$, and ${\cal L}_2(q)$ the invariant complex line. Choosing such subspaces at every point $q\in M$ we obtain two real-analytic $G$-invariant distributions ${\cal L}_1$ and ${\cal L}_2$ of $(n-1)$- and 1-dimensional complex subspaces on $M$, respectively. Lifting ${\cal L}_1$ and ${\cal L}_2$ to $G$ by means of $\Pi_{G,H}$, we obtain distributions ${\cal S}_1$ and ${\cal S}_2$ of real $(n^2-1)$- and $(n^2-2n+3)$-dimensional subspaces on $G$, respectively. Since these distributions are invariant under left translations on $G$, we will think of them as linear subspaces of ${\frak g}$. We have ${\frak g}={\cal S}_1+{\cal S}_2$ and ${\cal S}_1\cap{\cal S}_2={\frak h}$. Let ${\frak h}_j^{\perp}:={\frak h}^{\perp}\cap{\cal S}_j$, $j=1,2$. Clearly, ${\frak h}^{\perp}={\frak h}_1^{\perp}+{\frak h}_2^{\perp}$, $\hbox{dim}\,{\frak h}_1^{\perp}=2(n-1)$, $\hbox{dim}\,{\frak h}_2^{\perp}=2$, and ${\frak h}_j^{\perp}$ is $\hbox{Ad}(H)$-invariant for each $j$.

Fix complex coordinates $(\xi_1,\dots,\xi_n)$ in $T_H(G/H)$ in which $LH^0$ is given by $SU_{n-1}\times U_1$. Accordingly, we have ${\frak h}={\frak h}_1\oplus{\frak h}_2$, where ${\frak h}_1:={\frak{su}}_{n-1}$, ${\frak h}_2:={\frak u}_1$. Clearly, $\Phi$ maps ${\frak h}_1^{\perp}$ and ${\frak h}_2^{\perp}$ onto $\{\xi_n=0\}$ and $\{\xi_1=\dots=\xi_{n-1}=0\}$, respectively, and the following holds
\begin{equation}
\begin{array}{ll}
[{\frak h}_j^{\perp},{\frak h}_j]\subset {\frak h}_j^{\perp},& j=1,2,\\

[{\frak h}_j^{\perp},{\frak h}_k]=0, & j\ne k.
\end{array}\label{hperphii}
\end{equation}

Set ${\cal S}_j':={\frak h}_j^{\perp}+{\frak h}_j$ for $j=1,2$. We will now show that ${\cal S}_j'$ is an ideal in ${\frak g}$ for each $j$. For any two elements $v_1,v_2\in{\frak g}$ we write
\begin{equation}
[v_1,v_2]=u_1+u_2+w_1+w_2,\label{abcii}
\end{equation}
where $u_k\in{\frak h}_k^{\perp}$ and $w_k\in{\frak h}_k$, $k=1,2$. Suppose first that $v_1\in{\frak h}_1^{\perp}$, $v_2\in{\frak h}_2^{\perp}$. Applying to (\ref{abcii}) the element $\varphi_0$ of $\hbox{Ad}(H^0)$ that acts trivially on ${\frak h}_1^{\perp}$ and coincides with $-\hbox{id}$ on ${\frak h}_2^{\perp}$, we obtain that $u_1=0$, $w_1=0$, $w_2=0$. Next, applying to (\ref{abcii}) an element of $\hbox{Ad}(H^0)$ that acts trivially on ${\frak h}_2^{\perp}$ and transforms $v_1$ into $-v_1$ we obtain that $u_2=0$. Thus
\begin{equation}
[{\frak h}_1^{\perp},{\frak h}_2^{\perp}]=0.\label{h12}
\end{equation}
Let $v_1,v_2\in{\frak h}_1^{\perp}$. In this case the application of the transformation $\varphi_0$ to (\ref{abcii}) yields $u_2=0$. We now apply the Jacobi identity to $v_1,v_2,v$, where $v$ is an arbitrary element of $\in{\frak h}_2^{\perp}$. Then (\ref{hperphii}), (\ref{h12}) imply that $[w_2,v]=0$, and hence $w_2=0$. Thus
\begin{equation}
[{\frak h}_1^{\perp},{\frak h}_1^{\perp}]\subset{\cal S}_1'.\label{h11}
\end{equation}
Let finally $v_1,v_2\in{\frak h}_2^{\perp}$. Applying to (\ref{abcii}) elements of $\hbox{Ad}(H^0)$ that act trivially on ${\frak h}_2^{\perp}$ we see that $u_1$ and $w_1$ are invariant under all such transformations. Therefore $u_1=0$, $w_1=0$, and we have obtained     
\begin{equation}
[{\frak h}_2^{\perp},{\frak h}_2^{\perp}]\subset{\cal S}_2'.\label{h22}
\end{equation}

Identities (\ref{hperphii}), (\ref{h12}), (\ref{h11}), (\ref{h22}) yield that ${\cal S}_j'$ is an ideal in ${\frak g}$ for each $j$. Thus, for each $j$ the distribution ${\cal L}_j$ is integrable, and its integral manifolds form a foliation ${\frak F}_j$ of $M$ by connected complex submanifolds of dimension $n-1$ for $j=1$ and by connected complex curves for $j=2$. For $q\in M$ we denote by ${\frak F}_j(q)$ the leaf of the $j$th foliation passing through $q$.

Let $G_j$ be the (possibly non-closed) normal connected subgroup of $G$ with Lie algebra ${\cal S}_j'$ for $j=1,2$. For every $q\in M$ the leaf ${\frak F}_j(q)$ coincides with the orbit $G_jq$. Since the tangent space to $G_jq$ at $q'\in G_jq$ is spanned by the values of the holomorphic fundamental vector fields of the $G_j$-action at $q'$, it follows that ${\frak F}_j$ is a holomorphic foliation. Clearly, every two orbits of $G_j$ are holomorphically equivalent. The ineffectivity kernel $K_j$ of the action of $G_j$ on $G_jp$ is discrete for each $j$. Since $G_1/K_1$ acts properly on $G_1p$, it follows from \cite{I5} that the orbit $G_1p$ is holomorphically equivalent to $\CC^{n-1}$ by means of a map that transforms $G_1/K_1$ into the group $G_1(\CC^{n-1})$. Furthermore, $G_2p$ is holomorphically equivalent to one of $\BB^1$, $\CC$, $\CC\PP^1$ by means of a map that transforms $G_2/K_2$ into one of $\hbox{Aut}(\BB^1)$, $G(\CC)$, $G(\CC\PP^1)$, respectively.

We will now show that each $G_j$ is closed in $G$. We assume that $j=1$; for $j=2$ the proof is similar. Let ${\frak U}$ be a neighborhood of $0$ in ${\frak g}$ where the exponential map into $G$ is a diffeomorphism, and let ${\frak V}:=\exp({\frak U})$. To prove that $G_1$ is closed in $G$ it is sufficient to show that for some neighborhood ${\frak W}$ of $e\in G$, ${\frak W}\subset{\frak V}$, we have $G_1\cap {\frak W}= \exp({\cal S}_1'\cap {\frak U})\cap {\frak W}$. Assuming the opposite we obtain a sequence $\{g_j\}$ of elements of $G_1$ converging to $e$ in $G$ such that for every $j$ we have $g_j=\exp(a_j)$ with $a_j\in {\frak U}\setminus{\cal S}_1'$. Let $p_j:=g_jp$. For a neighborhood ${\cal V}$ of $p$ we denote by $N_p$ and $N_{p_j}$ the connected components of $G_1p\cap{\cal V}$ containing $p$ and $p_j$, respectively. We will now show that there exists a neighborhood ${\cal V}$ of $p$ such that $N_{p_j}\ne N_p$ for large $j$.

Let ${\frak U}''\subset {\frak U}'\subset {\frak U}$ be neighborhoods of 0 in ${\frak g}$ such that: (a) $\exp({\frak U}'')\cdot\exp({\frak U}'')\subset\exp({\frak U}')$; (b) $\exp({\frak U}'')\cdot\exp({\frak U}')\subset\exp({\frak U})$; (c) ${\frak U}'=-{\frak U}'$;  (d) $H\cap\exp({\frak U}')\subset\exp({\cal S}_1'\cap {\frak U}')$. We now choose ${\cal V}$ so that $N_p\subset\exp({\cal S}_1'\cap {\frak U}'')p$. Suppose that $p_j\in N_p$. Then we have $p_j=sp$ for some $s\in\exp({\cal S}_1'\cap {\frak U}'')$ and hence $t:=g_j^{-1}s$ is an element of $H$. For large $j$ we have $g_j^{-1}\in\exp({\frak U}'')$. Condition (a) now implies that $t\in\exp({\frak U}')$ and hence by (c), (d) we have $t^{-1}\in\exp({\cal S}_1'\cap {\frak U}')$. Therefore, by (b) we obtain $g_j\in\exp({\cal S}_1'\cap {\frak U})$ which contradicts our choice of $g_j$. Thus, for large $j$ we have $N_{p_j}\ne N_p$, and thus the orbit $G_1p$ accumulates to itself (we will use this term in the future in analogous situations). Below we will show that this is in fact impossible thus obtaining a contradiction.

Consider the set $S:=G_1p\cap G_2p$. The set $S$ contains a non-constant sequence converging to $p$. Clearly, $H^0$ preserves $S$. Since the $H^0$-orbit of a point in $S$ cannot have positive dimension, the subgroup $H^0$ fixes every point in $S$.  At the same time, any compact subgroup of dimension $n^2-2n$ in $G_1(\CC^{n-1})$ fixes exactly one point in $\CC^{n-1}$. This contradiction shows that $G_j$ is closed in $G$ for each $j$. Therefore, the action of $G_j$ on $M$ is proper and hence every leaf of ${\frak F}_j$ is closed in $M$, for $j=1,2$.

We will now show that the subgroup $K_j$ is in fact trivial for each $j=1,2$. Let first $j=1$. Since $G_1/K_1$ is isomorphic to the simply-connected group $G_1(\CC^{n-1})\simeq SU_{n-1}\ltimes\CC^{n-1}$ and since $G_1$ covers $G_1/K_1$ with fiber $K_1$, it follows that $K_1$ is trivial. Let $j=2$. If $G_2/K_2$ is isomorphic to $G(\CC)$, the triviality of $K_2$ follows as above. Further, the action of $G_{2\,p}^0$ on $G_2p$ is effective, and thus we have $K_2\setminus\{e\}\subset G_{2\,p}\setminus G_{2\,p}^0$. Suppose that $G_2/K_2$ is isomorphic to $\hbox{Aut}(\BB^1)$. Every maximal compact subgroup of $\hbox{Aut}(\BB^1)$ is 1-dimensional, hence so is every maximal compact subgroup of $G_2$. Since $G_{2\,p}^0$ is 1-dimensional, it is  maximal compact in $G_2$. Therefore $G_{2\,p}$ is connected, which implies that $K_2$ is trivial. Suppose next that $G_2/K_2$ is isomorphic to $G(\CC\PP^1)\simeq PSU_2$. If $K_2$ is non-trivial, then $G_2\simeq SU_2$ and $K_2\simeq\ZZ_2$. Then $G_{2\,p}^0$ is conjugate in $G_2$ (upon the identification of $G_2$ with $SU_2$) to the subgroup of matrices of the form
$$
\left(
\begin{array}{cc}
1/b & 0\\
0 & b
\end{array}
\right),
$$
where $|b|=1$ (see e.g. Lemma 2.1 of \cite{IKru}). Since this subgroup contains the center of $SU_2$, the subgroup $G_{2\,p}^0$ contains the center of $G_2$. In particular, $K_2\subset G_{2\,p}^0$ which contradicts the non-triviality of $K_2$.
Thus, $G_1$ is isomorphic to $G_1(\CC^{n-1})$ and $G_2$ is isomorphic to one of $\hbox{Aut}(\BB^1)$, $G(\CC)$, $G(\CC\PP^1)$.

Next, since ${\frak g}={\cal S}_1'\oplus{\cal S}_2'$ and $G_1$, $G_2$ are closed, the group $G$ is a locally direct product of $G_1$ and $G_2$. We claim that ${\mathscr T}:=G_1\cap G_2$ is trivial. Indeed, ${\mathscr T}$ is a discrete normal subgroup of each of $G_1$, $G_2$. However, every discrete normal subgroup of each of $G_1(\CC^{n-1})$, $\hbox{Aut}(\BB^1)$, $G (\CC)$, $G(\CC\PP^1)$ is trivial, since the center of each of these groups is trivial. Hence $G=G_1\times G_2$.

We will now observe that for every $q_1,q_2\in M$ the orbits $G_1q_1$ and $G_2q_2$ intersect at exactly one point. Let $g\in G$ be an element such that $gq_2=q_1$. It can be uniquely represented in the form $g=g_1g_2$ with $g_j\in G_j$ for $j=1,2$, and therefore we have $g_2q_2=g_1^{-1}q_1$. Hence the intersection $G_1q_1\cap G_2q_2$ is non-empty. Next, the fact that for every $q\in M$ the intersection $G_1q\cap G_2q$ consists of $q$ alone follows by the argument used at the end of the proof of the closedness of $G_1$, $G_2$.

Let $F_1$ be a biholomorphic map from $G_1p$ onto $\CC^{n-1}$ that transforms $G_1$ into $G_1(\CC^{n-1})$, and $F_2$ a biholomorphic map from $G_2p$ onto $M'$, where $M'$ is one of $\BB^1$, $\CC$, $\CC\PP^1$, that transforms $G_2$ into $G'$, where $G'$ is one of $\hbox{Aut}(\BB^1)$, $G(\CC)$, $G(\CC\PP^1)$, respectively. We will now construct a biholomorphic map ${\cal F}$ from $M$ onto $\CC^{n-1}\times M'$. For $q\in M$ consider $G_2q$ and let $r$ be the unique point of intersection of $G_1p$ and $G_2q$. Let $g\in G_1$ be an element such that $r=gp$. Then we set ${\cal F}(q):=(F_1(r), F_2(g^{-1}q))$. Clearly, ${\cal F}$ is a well-defined diffeomorphism from $M$ onto $\CC^{n-1}\times M'$. Since the foliation ${\frak F}_j$ is holomorphic for each $j$, the map ${\cal F}$ is in fact holomorphic. By construction, ${\cal F}$ transforms $G$ into $G_1(\CC^{n-1})\times G'$.

The proof is complete.\qed

\section{Actions of Type III}\label{type3}
\setcounter{equation}{0}

We start this section with a large number of examples of actions of type III. Some of the examples can be naturally combined into classes and some of the actions form parametric families. In what follows $n\ge 2$.
\vspace{0.5cm}

\noindent {\bf (i).} Here both the manifolds and the groups are represented as direct products.
\vspace{0.3cm}

{\bf (ia).} $M=M'\times\CC$, where $M'$ is one of $\BB^{n-1}$, $\CC^{n-1}$, $\CC\PP^{n-1}$, and $G=G'\times G_1(\CC)$, where $G'$ is one of the groups $\hbox{Aut}(\BB^{n-1})$, $G(\CC^{n-1})$, $G(\CC\PP^{n-1})$, respectively.
\vspace{0.3cm}

{\bf (ib).} $M=M'\times\CC^*$, where $M'$ is as in (ia), and $G=G'\times \CC^*$, where $G'$ is as in (ia). 
\vspace{0.3cm}

{\bf (ic).} $M=M'\times\TT$, where $M'$ is as in (ia) and $\TT$ is an elliptic curve; $G=G'\times \hbox{Aut}(\TT)^0$, where $G'$ is as in (ia).
\vspace{0.3cm}

{\bf (id).} $M=M'\times{\cal P}_{>}$\,, where $M'$ is as in (ia) and ${\cal P}_{>}:=\left\{\xi\in\CC:\hbox{Re}\,\xi>0\right\}$; $G=G'\times G({{\cal P}_{>}})$, where $G'$ as in (ia) and $G({{\cal P}_{>}})$ is the group of all maps of the form
\begin{equation}
\xi\mapsto \lambda\xi+ia,\label{gofp}
\end{equation}
with $a\in\RR$, $\lambda>0$.
\vspace{0.5cm}

\noindent{\bf (ii).} Parts (iib) and (iic) of this example are obtained by passing to quotients in Part (iia).

{\bf (iia).} $M=\BB^{n-1}\times\CC$, and $G$ consists of all maps of the form
$$
\begin{array}{lll}
z'&\mapsto&\displaystyle\frac{Az'+b}{cz'+d},\\
\vspace{0mm}&&\\
z_n&\mapsto&\displaystyle z_n+\ln(cz'+d)+a,
\end{array}
$$
where
\begin{equation}
\left(\begin{array}{cc}
A& b\\
c& d
\end{array}
\right)
\in SU_{n-1,1},\label{autballn-1}
\end{equation}
$z':=(z_1,\dots,z_{n-1})$ and $a\in\CC$. We denote this group by $G\left(\BB^{n-1}\times\CC\right)$. In fact, for $T\in\CC$ one can consider the following family of groups acting on $\BB^{n-1}\times\CC$
\begin{equation}
\begin{array}{lll}
z'&\mapsto&\displaystyle\frac{Az'+b}{cz'+d},\\
\vspace{0mm}&&\\
z_n&\mapsto&\displaystyle z_n+T\ln(cz'+d)+a,
\end{array}\label{groupbcs}
\end{equation}
where $A,a,b,c,d$ are as above. Example (ia) for $M'=\BB^{n-1}$ is included in this family for $T=0$. If $T\ne 0$, then conjugating group (\ref{groupbcs}) in $\hbox{Aut}(\BB^{n-1}\times\CC)$ by the automorphism
\begin{equation}
\begin{array}{lll}
z'&\mapsto & z'\\
z_n&\mapsto & z_n/T,
\end{array}\label{divs}
\end{equation}
we can assume that $T=1$.  
\vspace{0.3cm}

{\bf (iib).} $M=\BB^{n-1}\times\CC^*$, and for a fixed $T\in\CC^*$ the group $G$ consists of all maps of the form
\begin{equation}
\begin{array}{lll}
z'&\mapsto&\displaystyle\frac{Az'+b}{cz'+d},\\
\vspace{0mm}&&\\
z_n&\mapsto&\displaystyle\chi (cz'+d)^Tz_n,
\end{array}\label{groupbcstar}
\end{equation}
where $A,b,c,d$ are as in (iia) and $\chi\in\CC^*$. Example (ib) for $M'=\BB^{n-1}$ can be included in this family for $T=0$. This family is obtained from (\ref{groupbcs}) by passing to a quotient in the last variable.
\vspace{0.3cm}

{\bf (iic).} $M=\BB^{n-1}\times\TT$, where $\TT$ is an elliptic curve, and for a fixed $T\in\CC^*$ the group $G$ consists of all maps of the form
\begin{equation}
\begin{array}{lll}
z'&\mapsto&\displaystyle\frac{Az'+b}{cz'+d},\\
\vspace{0mm}&&\\

[z_n]&\mapsto&\displaystyle
\left[\chi (cz'+d)^Tz_n\right],
\end{array}\label{groupbtorus}
\end{equation}
where $A,b,c,d,\chi$ are as in (iib), $\TT$ is obtained from $\CC^*$ by taking the quotient with respect to the equivalence relation $z_n\sim dz_n$, for some $d\in\CC^*$, $|d|\ne 1$, and $[z_n]\in\TT$ is the equivalence class of a point $z_n\in\CC^*$. Example (ic) for $M'=\BB^{n-1}$ can be included in this family for $T=0$. Clearly, after passing to the quotient (\ref{groupbcstar}) turns into (\ref{groupbtorus}). 
\vspace{0.5cm}

\noindent{\bf (iii).} Part (iiib) of this example is obtained by passing to a quotient in\linebreak Part (iiia).  
\vspace{0.3cm}

{\bf (iiia).} $M=\CC^n$, and $G$ consists of all maps of the form
$$
\begin{array}{lll}
z'&\mapsto & e^{\hbox{\tiny Re}\,b}Uz'+a,\\
z_n&\mapsto & z_n+b,
\end{array}
$$
where $U\in U_{n-1}$, $a\in\CC^{n-1}$, $b\in\CC$. In fact, for $T\in\CC$ one can consider the following family of groups acting on $\CC^n$
\begin{equation}
\begin{array}{lll}
z'&\mapsto & e^{\hbox{\tiny Re}\,(T b)}Uz'+a,\\
z_n&\mapsto & z_n+b,
\end{array}\label{g3cns}
\end{equation}
where $U$, $a$, $b$ are as above. Example (ia) for $M'=\CC^{n-1}$ is included in this family for $T=0$. If $T\ne 0$, then conjugating group (\ref{g3cns}) in $\hbox{Aut}(\CC^n)$ by the automorphism
$$
\begin{array}{lll}
z'&\mapsto & z'\\
z_n&\mapsto & T z_n,
\end{array}
$$
we can assume that $T=1$.  
\vspace{0.3cm}

{\bf (iiib).} $M=\CC^{n-1}\times\CC^*$, and for a fixed $T\in\RR^*$ the group $G$ consists of all maps of the form
$$
\begin{array}{lll}
z'&\mapsto & e^{T\,\hbox{\tiny Re}\,b}Uz'+a,\\
z_n&\mapsto & e^{b}z_n,
\end{array}
$$
where $U, a, b$ are as in (iiia). Example (ib) for $M'=\CC^{n-1}$ can be included in this family for $T=0$. This family is obtained from (\ref{g3cns}) for $T\in\RR^*$ by passing to a quotient in the last variable. 
\vspace{0.5cm}

\noindent{\bf (iv).} Parts (ivb) and (ivc) of this example are obtained by passing to quotients in Part (iva).
\vspace{0.3cm}

{\bf (iva).} $M=\CC^n$, and $G$ consists of all maps of the form
$$
\begin{array}{lll}
z' & \mapsto & Uz'+a,\\
z_n & \mapsto & z_n+\langle Uz',a\rangle+b,
\end{array}
$$
where $U\in U_{n-1}$, $a\in\CC^{n-1}$, $b\in\CC$, and $\langle\cdot\,,\cdot\rangle$ is the inner product in $\CC^{n-1}$. We denote this group by ${\frak G}(\CC^n)$. In fact, for $T\in\CC$ one can consider the following family of groups acting on $\CC^n$
\begin{equation}
\begin{array}{lll}
z' & \mapsto & Uz'+a,\\
z_n & \mapsto & z_n+T\langle Uz',a\rangle+b,
\end{array}\label{g4cns}
\end{equation}
where $U$, $a$, $b$ are as above. Example (ia) for $M'=\CC^{n-1}$ is included in this family for $T=0$. If $T\ne 0$, then conjugating group (\ref{g4cns}) in $\hbox{Aut}(\CC^n)$ by automorphism (\ref{divs}), we can assume that $T=1$.  
\vspace{0.3cm}

{\bf (ivb).} $M=\CC^{n-1}\times\CC^*$, and for a fixed $0\le\tau<2\pi$ the group $G$ consists of all maps of the form
\begin{equation}
\begin{array}{lll}
z' & \mapsto & Uz'+a,\\
z_n & \mapsto &\chi\exp\Bigl(e^{i\tau}\langle Uz',a\rangle\Bigr)z_n,
\end{array}\label{groupcstartau}
\end{equation}
where $U,a$ are as in (iva) and $\chi\in\CC^*$. In fact, for $T\in\CC$ one can consider the following family of groups acting on $\CC^{n-1}\times\CC^*$
\begin{equation}
\begin{array}{lll}
z' & \mapsto & Uz'+a,\\
z_n & \mapsto &\chi\exp\Bigl(T\langle Uz',a\rangle\Bigr)z_n,
\end{array}\label{g4cnsstar}
\end{equation}
where $U,a,\chi$ are as above. Example (ib) for $M'=\CC^{n-1}$ is included in this family for $T=0$. For $T\ne 0$ this family is obtained from (\ref{g4cns}) by passing to a quotient in the last variable. Furthermore, conjugating group (\ref{g4cnsstar}) for $T\ne 0$ in $\hbox{Aut}(\CC^{n-1}\times\CC^*)$ by the automorphism
$$
\begin{array}{lll}
z'&\mapsto & \sqrt{|T|}z'\\
z_n&\mapsto & z_n,
\end{array}
$$
we obtain the group defined in (\ref{groupcstartau}) for $\tau=\hbox{arg}\,T$.
\vspace{0.3cm}

{\bf (ivc).} $M=\CC^{n-1}\times\TT$, where $\TT$ is an elliptic curve, and for a fixed $0\le\tau<2\pi$ the group $G$ consists of all maps of the form
\begin{equation}
\begin{array}{lll}
z' & \mapsto & Uz'+a,\\

[z_n] & \mapsto &\left[\chi\exp\Bigl(e^{i\tau}\langle Uz',a\rangle\Bigr)z_n\right],
\end{array}\label{grouptorustau}
\end{equation}
where $U,a,\chi$ are as in (ivb), $\TT$ is obtained from $\CC^*$ by taking the quotient with respect to the equivalence relation $z_n\sim dz_n$, for some $d\in\CC^*$, $|d|\ne 1$, and $[z_n]\in\TT$ is the equivalence class of a point $z_n\in\CC^*$. In fact, for $T\in\CC$ one can consider the following family of groups acting on $\CC^{n-1}\times\TT$
\begin{equation}
\begin{array}{lll}
z' & \mapsto & Uz'+a,\\

[z_n] & \mapsto &\left[\chi\exp\Bigl(T\langle Uz',a\rangle\Bigr)z_n\right],
\end{array}\label{g4cnst}
\end{equation}
where $U,a,\chi$ are as above. Example (ic) for $M'=\CC^{n-1}$ is included in this family for $T=0$. For $T\ne 0$ this family is obtained from (\ref{g4cnsstar}) by passing to the quotient described above. Furthermore, conjugating group (\ref{g4cnst}) for $T\ne 0$ in $\hbox{Aut}(\CC^{n-1}\times\TT)$ by the automorphism
$$
\begin{array}{lll}
z'&\mapsto & \sqrt{|T|}z'\\
\xi &\mapsto & \xi,
\end{array}
$$
where $\xi\in\TT$, we obtain the group defined in (\ref{grouptorustau}) for $\tau=\hbox{arg}\,T$.
\vspace{0.5cm}

\noindent{\bf (v).} $M=\CC^{n-1}\times{\cal P}_{>}$\,, and for a fixed $T\in\RR^*$ the group $G$ consists of all maps of the form
\begin{equation}
\begin{array}{lll}
z'&\mapsto & \lambda^TUz'+a,\\
z_n&\mapsto & \lambda z_n+ib,
\end{array}\label{groupv}
\end{equation}
where $U\in U_{n-1}$, $a\in\CC^{n-1}$, $b\in\RR$, $\lambda>0$. Example (id) for $M'=\CC^{n-1}$ can be included in this family for $T=0$.
\vspace{0.5cm}

\noindent {\bf (vi).} $M=\CC^n$, and for fixed $k_1,k_2\in\ZZ$, $(k_1,k_2)=1$, $k_1>0$, $k_2\ne 0$, the group $G$ consists of all maps of the form (\ref{groupcn}) with $U\in H_{k_1,k_2}^n$ (see (\ref{mat})). We denote this group by $G_{k_1,k_2}(\CC^n)$. Example (ia) for $M'=\CC^{n-1}$ can be included in this family for $k_2=0$.  
\vspace{0.5cm}

\noindent{\bf (vii).} Part (viib) of this example is obtained by passing to a quotient in\linebreak Part (viia).  
\vspace{0.3cm}

{\bf (viia).} $M=\CC^{n{}*{}}/\ZZ_l$, where $\CC^{n{}*{}}:=\CC^n\setminus\{0\}$, $l\in\NN$,  and the group $G$ consists of all maps of the form
\begin{equation}
\{z\}\mapsto\{\lambda Uz\},\label{gviiia}
\end{equation}
where $U\in U_n$, $\lambda>0$, and $\{z\}\in\CC^{n{}*{}}/\ZZ_l$ is the equivalence class of a point $z\in\CC^{n{}*{}}$.
\vspace{0.3cm}

{\bf (viib).} $M=M_d/\ZZ_l$, where $M_d$ is the Hopf manifold $\CC^{n{}*{}}/\{z{}\sim{}dz\}$, for $d\in\CC^*$, $|d|\ne 1$, and $l\in\NN$; the group $G$ consists of all maps of the form 
\begin{equation}
\{[z]\}\mapsto\left\{[\lambda Uz\right]\},\label{gviiib}
\end{equation}
where $U,\lambda$ are as in (viia), $[z]\in M_d$ denotes the equivalence class of a point $z\in\CC^{n{}*{}}$, and $\{[z]\}\in M_d/\ZZ_l$ denotes the equivalence class of $[z]\in M_d$.   
\vspace{0.5cm}

\noindent{\bf (viii).} In this example the manifolds are the open orbits of the action of a group of affine transformations on $\CC^n$. Let $G_{\cal P}$ be the group of all maps of the form
$$
\begin{array}{lll}
z' & \mapsto & \lambda Uz'+a,\\
z_n & \mapsto & \lambda^2z_n+2\lambda\langle Uz',a\rangle+|a|^2+ib,
\end{array}
$$
where $U\in U_{n-1}$, $a\in\CC^{n-1}$, $b\in\RR$, $\lambda>0$.\vspace{0.3cm}

{\bf (viiia).} $M={\cal P}_{>}^n$\,, $G=G_{\cal P}$, where
\begin{equation}
{\cal P}_{>}^n:=\left\{(z',z_n)\in\CC^{n-1}\times\CC: \hbox{Re}\,z_n>|z'|^2\right\}.\label{p>}
\end{equation}
Observe that ${\cal P}_{>}^n$ is holomorphically equivalent to $\BB^n$. 
\vspace{0.3cm}

{\bf (viiib).} $M={\cal P}_{<}^n$\,, $G=G_{\cal P}$, where
\begin{equation}
{\cal P}_{<}^n:=\left\{(z',z_n)\in\CC^{n-1}\times\CC: \hbox{Re}\,z_n<|z'|^2\right\}.\label{p<}
\end{equation}
Observe that ${\cal P}_{<}^n$ is holomorphically equivalent to $\CC\PP^n\setminus(\overline{\BB^n}\cup L)$, where $L$ is a complex hyperplane tangent to $\partial\BB^n$ at some point.
\vspace{0.5cm}

\noindent {\bf (ix).} Here $n=2$, $M=\BB^1\times\CC$, and $G$ consists of all maps of the form
$$
\begin{array}{lll}
z_1 & \mapsto & \displaystyle\frac{az_1+b}{\overline{b}z_1+\overline{a}},\\
\vspace{0.1cm}\\
z_2 & \mapsto & \displaystyle\frac{z_2+cz_1+\overline{c}}{\overline{b}z_1+\overline{a}},\\
\end{array}
$$
where $a,b\in\CC$, $|a|^2-|b|^2=1$, $c\in\CC$. We denote this group by ${\frak G}(\BB^1\times\CC)$.
\vspace{0.5cm}

\noindent{\bf (x).} Here $n=3$, $M=\CC\PP^3$, and $G$ consists of all maps of the form (\ref{groupcpn}) for $n=3$ with $U\in Sp_2$ (where $Sp_2$ is the compact real form of $Sp_4(\CC)$). We denote this group of maps by $G_1(\CC\PP^3)$. It is isomorphic to $Sp_2/\ZZ_2$.
\vspace{0.5cm}

\noindent{\bf (xi).} Let $n=3$ and $(z:w)$ be homogeneous coordinates in $\CC\PP^3$ with $z=(z_1:z_2)$, $w=(w_1:w_2)$. Set $M=\CC\PP^3\setminus\{w=0\}$ and let $G$ be the group of all maps of the form
\begin{equation}
\begin{array}{lll}
z&\mapsto&Uz+Aw,\\
w&\mapsto& Vw,
\end{array}\label{gcp3cp1}
\end{equation}
where $U,V\in SU_2$, and 
$$
A=\left(
\begin{array}{lr}
a & i\overline{b}\\
b & -i\overline{a}
\end{array}
\right),
$$
for some $a,b\in\CC$. We denote this group by ${\frak G}\left(\CC\PP^3\setminus\CC\PP^1\right)$. 
\vspace{0.5cm}

\noindent{\bf (xii).} Here $n=3$, $M=\CC^3$, and $G$ consists of all maps of the form
$$
\begin{array}{lll}
z'&\mapsto&Uz'+a,\\
z_3&\mapsto&\hbox{det}\,U\,z_3+\left[
\left(
\begin{array}{ll}
0 & 1\\
-1 & 0
\end{array}
\right)
Uz'\right]\cdot a+b,
\end{array}
$$
where $z':=(z_1,z_2)$, $U\in U_2$, $a\in\CC^2$, $b\in\CC$, and $\cdot$ is the dot product in $\CC^2$. We denote this group of maps by $G_3(\CC^3)$.
\vspace{0.5cm}

The result of this section is the following theorem.

\begin{theorem}\label{classtype3}\sl Let $M$ be a connected complex manifold of dimension $n\ge2$ and $G\subset\hbox{Aut}(M)$ a connected Lie group with $d_G=n^2+1$ that acts properly on $M$. If the pair $(M,G)$ is of type III, then it is equivalent to one of the pairs listed in (i)--(xii) above.
\end{theorem}

The following lemma is analogous to Lemma \ref{fullisotropyii}, and we omit the proof.

\begin{lemma}\label{fullisotropy}\sl Let $H\subset U_n$, with $n\ge 2$, be a closed subgroup, and let $H^0=H_{k_1,k_2}^n$ for some $k_1,k_2$. Assume that for $n=2$ we have $k_1\ne \pm k_2$ (that is, for $n=2$ the cases $k_1=k_2=1$ and $k_1=1$, $k_2=-1$ are excluded from consideration). Then for some $m\in\NN$ the group $H$ consists of all matrices of the form
$$
\left(\begin{array}{cc}
A & 0\\
0 & \alpha\cdot a
\end{array}\right),
$$
where $A\in U_{n-1}$, $a\in(\det A)^{\frac{k_2}{k_1}}$, $\alpha^m=1$.
\end{lemma}

\begin{remark}\label{remlemma}\rm For $n=2$ and either $k_1=1$, $k_2=1$ or $k_1=1$, $k_2=-1$ it is easy to construct an example of a closed disconnected subgroup $H\subset U_2$ with $H^0=H^2_{k_1,k_2}$ that contains matrices with non-zero anti-diagonal elements.
\end{remark}

\noindent {\bf Proof of Theorem \ref{classtype3}:} We will consider two cases. Case 1 occupies most of the proof and splits into a number of subcases, as shown in the following diagram.

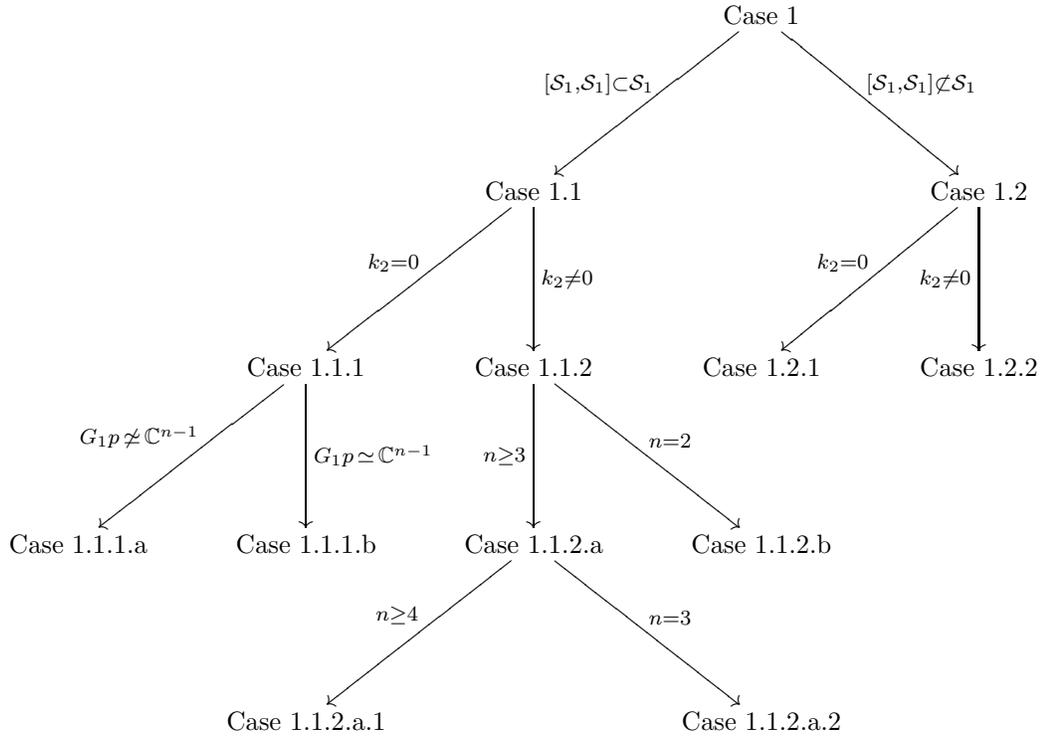
\xymatrix{
&&&\hbox{\footnotesize Case 1}\ar@{->}[ldd]_{[{\cal S}_1,{\cal S}_1]\subset{\cal S}_1} \ar@{->}[rdd]^{[{\cal S}_1,{\cal S}_1]\not\subset{\cal S}_1}\\
 \\
&&\hbox{\footnotesize Case 1.1}\ar@{->}[ldd]_{k_2=0} \ar@{->}[dd]^{k_2\ne 0}&&\hbox{\footnotesize Case 1.2}\ar@{->}[ldd]_{k_2=0} \ar@{->}[dd]_{k_2\ne 0}\\
\\
&\hbox{\footnotesize Case 1.1.1}\ar@{->}[ldd]_{G_1p\,\not\simeq\,\CC^{n-1}} \ar@{->}[dd]^{G_1p\,\simeq\,\CC^{n-1}}&\hbox{\footnotesize Case 1.1.2}\ar@{->}[dd]_{n\ge 3} \ar@{->}[rdd]^{n=2}&\hbox{\footnotesize Case 1.2.1}&\hbox{\footnotesize Case 1.2.2}\\
\\
\hbox{\footnotesize Case 1.1.1.a}&\hbox{\footnotesize Case 1.1.1.b}&\hbox{\footnotesize Case 1.1.2.a}\ar@{->}[ldd]_{n\ge 4} \ar@{->}[rdd]^{n=3}&\hbox{\footnotesize Case 1.1.2.b}\\
\\
&\hbox{\footnotesize Case 1.1.2.a.1}&&\hbox{\footnotesize Case 1.1.2.a.2}
}

\vspace{0.3cm}

{\bf Case 1.} {\it Suppose that either $n\ge 3$, or $n=2$ and $k_1\ne \pm k_2$.} Fix $p\in M$ and introduce $H$, $\Phi$, ${\frak h}$ and ${\frak h}^{\perp}$ as in the proof of Theorem \ref{classtype2}. By Lemma \ref{fullisotropy}, at every $q\in M$ there are exactly two non-trivial proper $LG_q$-invariant complex subspaces ${\cal L}_1(q)$ and ${\cal L}_2(q)$ in $T_q(M)$. For $n\ge 3$, ${\cal L}_1(q)$ denotes, as before, the invariant subspace of dimension $n-1$, and ${\cal L}_2(q)$ the invariant complex line. For $n=2$, ${\cal L}_1(q)$ and ${\cal L}_2(q)$ denote the invariant complex lines such that for every $\varphi\in LG_q^0$ we have $\varphi(\xi_j)=u^{k_j}\xi_j$, for $j=1,2$, where $\xi_j\in{\cal L}_j(q)$ and $|u|=1$. Choosing such subspaces at every point $q\in M$ we obtain two real-analytic $G$-invariant distributions ${\cal L}_1$ and ${\cal L}_2$ of $(n-1)$- and 1-dimensional complex subspaces on $M$, respectively. As in the proof of Theorem \ref{classtype2} these distribution lead to linear subspaces ${\cal S}_1$ and ${\cal S}_2$ of ${\frak g}$. Let, as before, ${\frak h}_j^{\perp}:={\frak h}^{\perp}\cap{\cal S}_j$, $j=1,2$. Clearly, ${\frak h}_j^{\perp}$ is $\hbox{Ad}(H)$-invariant for each $j$, and therefore we have
\begin{equation}
[{\frak h}_j^{\perp},{\frak h}]\subset {\frak h}_j^{\perp},\quad j=1,2.\label{hperph}
\end{equation}

We fix complex coordinates $(\xi_1,\dots,\xi_n)$ in $T_H(G/H)$ in which $LH^0$ is given by $H_{k_1,k_2}^n$ and thus realize ${\frak h}$ as the Lie algebra of matrices of the form
\begin{equation}
\left(\begin{array}{cc}
{\frak A}& 0\\
0& k_2/k_1\cdot\hbox{trace}\,{\frak A}
\end{array}\right),\label{frakA}
\end{equation}
where ${\frak A}\in{\frak u}_{n-1}$. We now choose bases in ${\frak h}_1^{\perp}$ and ${\frak h}_2^{\perp}$  as follows. Set
$$
a_k:=\Phi^{-1}(e_k),\quad b_k:=\Phi^{-1}(ie_k),\quad k=1,\dots,n,
$$
where $e_k$ denotes the $k$th coordinate vector for $k=1,\dots,n$. Then $\{a_r,b_r\}$, $r=1,\dots,n-1$, is a basis in ${\frak h}_1^{\perp}$ and $\{a_n,b_n\}$ is a basis in ${\frak h}_2^{\perp}$. In this basis the pull-back to ${\frak h}^{\perp}$ of the operator of complex structure on $T_H(G/H)$ is given by
\begin{equation}
\begin{array}{llr}
a_k&\mapsto& b_k,\\
b_k &\mapsto& -a_k,\\
\end{array}\label{j0}
\end{equation}
for $k=1,\dots,n$. 

Observe that $[a_n,b_n]$ is $\hbox{Ad}(H)$-invariant, and therefore lies in ${\cal S}_2$. Together with (\ref{hperph}) this implies that ${\cal S}_2$ is a Lie subalgebra of ${\frak g}$. Thus, the distribution ${\cal L}_2$ is integrable, and its integral manifolds form a foliation ${\frak F}_2$ of $M$ by connected complex curves. For $q\in M$ we denote by ${\frak F}_2(q)$ the leaf of the foliation passing through $q$. Since the foliation ${\frak F}_2$ is $G$-invariant, its leaves are pairwise holomorphically equivalent. 

Let $G_2$ be the (possibly non-closed) connected subgroup of $G$ with Lie algebra ${\cal S}_2$ (observe that $H^0\subset G_2$). The leaf ${\frak F}_2(p)$ passing through the distinguished point $p$ coincides with the orbit $G_2p$, and for $g\in G$ the leaf ${\frak F}_2(gp)$ is an orbit of the group $gG_2g^{-1}$. The ineffectivity kernel $K_2$ of the action of $G_2$ on $G_2p$ is either $(n^2-2n)$-dimensional (if $k_2\ne 0$) or $(n-1)^2$-dimensional (if $k_2=0$). Hence $G_2p$ is holomorphically equivalent to one of: (a) $\BB^1$, $\CC$, $\CC\PP^1$, if $k_2\ne 0$ by means of a map that transforms $G_2/K_2$ into one of $\hbox{Aut}(\BB^1)$, $G(\CC)$, $G(\CC\PP^1)$, respectively; (b) ${\cal P}_{>}$ , $\CC$, $\CC^*$, or an elliptic curve $\TT$, if $k_2=0$ by means of a map that transforms $G_2/K_2$ into one of $G({\cal P}_{>})$ (see (\ref{gofp})), $G_1(\CC)$, $\hbox{Aut}(\CC^*)^0$, or $\hbox{Aut}(\TT)^0$, respectively.

Further, for each $1\le r\le n-1$ we fix a non-zero element $t_r\in{\frak h}$ such that the one-parameter subgroup $H_r$ of $H^0$ arising from $t_r$ is represented in $LH^0$ by transformations that act on ${\cal L}_1(p)$ as follows
$$
\begin{array}{lll}
\xi_s&\mapsto&\xi_s,\quad s=1,\dots,n-1,\quad s\ne r,\\
\xi_r&\mapsto&e^{i\psi}\xi_r,
\end{array}
$$  
where $\psi\in\RR$. We choose $t_r$ so that the element
\begin{equation}
t_0:=\sum_{r=1}^{n-1} t_r\label{t0}
\end{equation}
generates the center ${\frak z}$ of ${\frak h}$. We have
\begin{equation}
\begin{array}{llc}
[t_r,t_s]&=&0,\\

[t_r,a_s]&=&\mu\delta_r^s b_r,\\

[t_r,b_s]&=&-\mu\delta_r^s a_r,
\end{array}\label{identities1}
\end{equation}
for some $\mu\in\RR^*$, where $\delta_r^s$ is the Kronecker symbol and $r,s=1,\dots,n-1$. 

We also note that ${\frak h}={\frak z}\oplus{\frak h}'$, where ${\frak h}'$ is the derived algebra of ${\frak h}$ (which is is isomorphic to ${\frak{su}}_{n-1}$), and the following holds
$$
[{\frak h}',{\frak h}_2^{\perp}]=0.
$$

{\bf Case 1.1.} {\it Suppose that ${\cal S}_1$ is a Lie subalgebra of ${\frak g}$.} In this case the distribution ${\cal L}_1$ is integrable and gives rise to a foliation ${\frak F}_1$ of $M$ by connected complex submanifolds of dimension $n-1$. For $q\in M$ we denote by ${\frak F}_1(q)$ the leaf of the foliation ${\frak F}_1$ passing through $q$. Since the foliation ${\frak F}_1$ is $G$-invariant, its leaves are pairwise holomorphically equivalent. Let $G_1$ be the connected subgroup of $G$ with Lie algebra ${\cal S}_1$ (clearly, we have $H^0\subset G_1$). The leaf ${\frak F}_1(p)$ coincides with the orbit $G_1p$, and for $g\in G$ the leaf ${\frak F}_1(gp)$ is an orbit of the group $gG_1g^{-1}$. The group $G_1$ acts properly on $G_1p$ and the ineffectivity kernel $K_1$ of this action is finite, hence by \cite{Ka} the orbit $G_1p$ is holomorphically equivalent to one of $\BB^{n-1}$, $\CC^{n-1}$, $\CC\PP^{n-1}$ by means of a map that transforms $G_1/K_1$ into one of the groups $\hbox{Aut}(\BB^{n-1})$, $G(\CC^{n-1})$, $G(\CC\PP^{n-1})$, respectively.

We will now show that each $G_j$ is a closed subgroup of $G$. Assuming the opposite, we obtain that the orbit $G_jp$ accumulates to itself (in the sense explained in the proof of Theorem \ref{classtype2} above). Consider the set $S:=G_1p\cap G_2p$. The set $S$ contains a non-constant sequence converging to $p$. Clearly, $H^0$ preserves $S$. Since the $H^0$-orbit of a point in $S$ cannot have positive dimension, the subgroup $H^0$ fixes every point in $S$.  At the same time, any compact subgroup of dimension $(n-1)^2$ in $\hbox{Aut}(\BB^{n-1})$ or $G(\CC^{n-1})$ fixes exactly one point in $\BB^{n-1}$ or $\CC^{n-1}$, respectively (every such subgroup is maximal compact), and any closed subgroup in $G(\CC\PP^{n-1})$ of dimension $(n-1)^2$ fixes at most two points in $\CC\PP^{n-1}$ (a two-point set can only occur for $n=2$). This contradiction shows that $G_j$ is closed in $G$ for each $j$. Therefore, the action of $G_j$ on $M$ is proper and hence every leaf of ${\frak F}_j$ is closed in $M$, for $j=1,2$.

Further, for any $v_1,v_2\in{\frak h}_1^{\perp}$ consider the commutator $[v_1,v_2]$ and write it in the general form $[v_1,v_2]=a+c$, where $a\in{\frak h}_1^{\perp}$, $c\in{\frak h}$. Applying to this identity the element of $\hbox{Ad}(Z)$, where $Z$ is the center of $H^0$, that acts as $-\hbox{id}$ on ${\frak h}_1^{\perp}$, we see that $a=0$, hence 
\begin{equation}
[{\frak h}_1^{\perp},{\frak h}_1^{\perp}]\subset{\frak h}.\label{v1v2}
\end{equation}
Let ${\frak P}_r$ be the subspace of ${\frak h}_1$ spanned by $a_r,b_r$, $r=1,\dots,n-1$. If $v_1\in {\frak P}_r$ and $v_2\in {\frak P}_s$ for $r\ne s$, then, applying to $[v_1,v_2]$ the elements of $\hbox{Ad}(H_r)$ and $\hbox{Ad}(H_s)$ that act as $-\hbox{id}$ on ${\frak P}_r$ and ${\frak P}_s$, respectively, we see that $c$ changes sign under each of these transformations, and therefore is equal to 0 (see (\ref{frakA})). Thus, we have
\begin{equation}
[{\frak P}_r,{\frak P}_s]=0,\quad r,s=1,\dots,n-1,\, r\ne s.\label{PrPs}
\end{equation}

Next, we observe that, since every commutator $[a_r,b_r]$ for $r=1,\dots,n-1$ is $\hbox{Ad}(H_s)$-invariant for $s=1,\dots,n-1$, the following holds
\begin{equation}
[a_r,b_r]=\sum_{s=1}^{n-1}\sigma_r^st_s,\label{formarbr}
\end{equation}
for some $\sigma_r^s\in\RR$. Applying to $[a_r,b_r]$ the transformations from $\hbox{Ad}(H^0)$ that interchange $a_r$ with $a_s$ and $b_r$ with $b_s$ for $s=1,\dots,n-1$, $s\ne r$, leaving the other elements of the basis $\{a_k,b_k\}_{k=1,\dots,n-1}$ in ${\frak h}_1^{\perp}$ fixed, we see that $\sigma_r^r=\sigma_s^s:=\sigma$ for $r,s=1,\dots,n-1$. Further, the Jacobi identity applied to all triples $\{a_r,b_r,a_s\}$ with $r,s=1,\dots,n-1$, $r\ne s$, together with (\ref{identities1}), (\ref{PrPs}) implies that $\sigma_r^s=0$ for $r\ne s$, and therefore we have
\begin{equation}
[a_r,b_r]=\sigma t_r,\quad r=1,\dots,n-1.\label{arbr}
\end{equation}

Relations (\ref{identities1}) and (\ref{arbr}) imply that the subspace ${\frak g}_r$ of ${\cal S}_1$ spanned by $a_r,b_r,t_r$ is in fact a Lie subalgebra for $r=1,\dots,n-1$. It is isomorphic to the Lie algebra of $\hbox{Aut}(\BB^1)$ if and only if $\mu\sigma<0$, to the Lie algebra of $G(\CC\PP^1)$ if and only if $\mu\sigma>0$ and to the Lie algebra of $G(\CC)$ if and only if $\sigma=0$. Let ${\frak G}_r$ be the connected subgroup of $G_1$ with Lie algebra ${\frak g}_r$. The tangent space $T_p({\frak G}_rp)$ at $p$ to the orbit ${\frak G}_rp$ is a complex line in $T_p(M)$, and $LH_r$ acts on $T_p({\frak G}_rp)$ effectively. Therefore, the orbit ${\frak G}_rp$ is a complex curve in $G_1p$ and there exists a discrete subgroup ${\frak K}_r$ of ${\frak G}_r$ such that ${\frak G}_rp$ is holomorphically equivalent to one of $\BB^1$, $\CC\PP^1$, $\CC$ by means of a map that transforms ${\frak G}_r/{\frak K}_r$ into one of the groups $\hbox{Aut}(\BB^1)$, $G(\CC\PP^1)$, $G(\CC)$, respectively. Hence the restriction of a $G_1$-invariant Hermitian metric on $G_1p$ to ${\frak G}_rp$ has a positive, negative, or zero curvature, respectively. This shows that $\mu\sigma<0$ if and only if $G_1p$ is equivalent to $\BB^{n-1}$; $\mu\sigma>0$ if and only if $G_1p$ is equivalent to $\CC\PP^{n-1}$; and $\sigma=0$ if and only if $G_1p$ is equivalent to $\CC^{n-1}$. 
    
{\bf Case 1.1.1.} {\it Suppose that $k_2=0$.} In this case $H^0$ acts on $G_1p$ effectively, and applying an argument from the proof of Theorem 2.1 in \cite{I5} (see also the proof of Theorem \ref{classtype2} above), we obtain that $K_1$ is trivial. Furthermore, we have
\begin{equation}
[{\frak h}_2^{\perp},{\frak h}]=0.\label{h2perph}
\end{equation}
Next, consider the commutator $[v_1,v_2]$, where $v_j\in{\frak h}_j^{\perp}$, $j=1,2$, and write it in the general form $[v_1,v_2]=a+b+c$, where $a\in{\frak h}_1^{\perp}$, $b\in{\frak h}_2^{\perp}$, $c\in{\frak h}$. Applying to this identity the element of $\hbox{Ad}(Z)$ that acts as $-\hbox{id}$ on ${\frak h}_1^{\perp}$ and taking into account that $\hbox{Ad}(Z)$ acts trivially on ${\frak h}_2^{\perp}$ and ${\frak h}$, we see that $b=c=0$, that is,
\begin{equation}
[{\frak h}_1^{\perp},{\frak h}_2^{\perp}]\subset{\frak h}_1^{\perp}.\label{h1h2}
\end{equation}
Relations (\ref{h2perph}) and (\ref{h1h2}) imply that ${\cal S}_1$ is an ideal in ${\frak g}$,  and hence $G_1$ is normal in $G$. This yields that every leaf of ${\frak F}_1$ is a $G_1$-orbit, and therefore ${\frak F}_1$ is a holomorphic foliation.

As a special case of (\ref{h1h2}), for every $v\in{\frak h}_2^{\perp}$ and $1\le r\le n-1$ we have
$$
\begin{array}{lll}
[a_r,v]&=&\displaystyle\sum_{s=1}^{n-1}(\lambda_sa_s+\omega_sb_s),\\ 
\vspace{-0.3cm}\\

[b_r,v]&=&\displaystyle\sum_{s=1}^{n-1}(\lambda'_sa_s+\omega'_sb_s),
\end{array}
$$
where $\lambda_s,\lambda'_s,\omega_s,\omega'_s\in\RR$. Applying to the above identities the family of operators $\hbox{Ad}(H_r)$ we see that $\lambda_s=\lambda'_s=\omega_s=\omega'_s=0$ for $s\ne r$ and that $\lambda'_r=-\omega_r$, $\omega'_r=\lambda_r$. Hence for $1\le r\le n-1$ we have
\begin{equation}
\begin{array}{ll}
[a_r,a_n]=\alpha_r a_r+\beta_r b_r,& [b_r,a_n]=-\beta_r a_r+\alpha_r b_r,\\ 

[a_r,b_n]=\gamma_r a_r+\delta_r b_r,& [b_r,b_n]=-\delta_r a_r+\gamma_r b_r,
\end{array}\label{identities2}
\end{equation}
for some $\alpha_r,\beta_r,\gamma_r,\delta_r\in\RR$. Applying to $[a_r,a_n]$ and $[a_r,b_n]$ a transformation from $\hbox{Ad}(H^0)$ that interchanges $a_r$ with $a_s$ and $b_r$ with $b_s$, we see that 
\begin{equation}
\begin{array}{ll}
\alpha_r=\alpha_s:=\alpha, & \beta_r=\beta_s:=\beta,\\
\gamma_s=\gamma_s:=\gamma, & \delta_r=\delta_s:=\delta,
\end{array}\label{thesame}
\end{equation}
for $r,s=1,\dots,n-1$.

Further, we have
\begin{equation}
[a_n,b_n]=\rho a_n+\nu b_n+t,\label{id3}
\end{equation}
where $\rho,\nu\in\RR$ and $t\in{\frak h}$. Since the group $\hbox{Ad}(H^0)$ acts trivially on ${\frak h}_2^{\perp}$, it fixes $t$, which implies that $t$ lies in the center of ${\frak h}$ and therefore
\begin{equation}
t=\eta t_0,\label{id3prime}
\end{equation} 
for some $\eta\in\RR$. For $r=1,\dots,n-1$ the Jacobi identity applied to the triple $\{a_r, a_n, b_n\}$, together with (\ref{t0}), (\ref{identities1}), (\ref{identities2}), (\ref{id3}), (\ref{id3prime}) implies
\begin{equation}
\begin{array}{l}
\alpha\rho+\gamma\nu=0,\\
\beta\rho+\delta\nu=\mu\eta.
\end{array}\label{id4}
\end{equation}
Define
\begin{equation}
\begin{array}{lll}
a_n'&:=&a_n+\displaystyle\frac{\beta}{\mu}t_0,\\
\vspace{-0.3cm}\\
b_n'&:=&b_n+\displaystyle\frac{\delta}{\mu}t_0.
\end{array}\label{anprimebnprime}
\end{equation}
It then follows from (\ref{h2perph}), (\ref{identities2}), (\ref{thesame}), (\ref{id3}), (\ref{id3prime}) and the second identity in (\ref{id4}) that
\begin{equation}
\begin{array}{ll}
[a_r,a_n']=\alpha a_r,& [b_r,a_n']=\alpha b_r,\\ 

[a_r,b_n']=\gamma a_r,& [b_r,b_n']=\gamma b_r,\\

[a_n',b_n']=\rho a_n'+\nu b_n'.&
\end{array}\label{identities44}
\end{equation}
In particular, the subspace ${\cal S}_2'$ of ${\frak g}$ spanned by $a_n'$, $b_n'$ is a Lie subalgebra in ${\frak g}$. Clearly, ${\frak g}$ is a semidirect sum of ${\cal S}_1$ and ${\cal S}_2'$. Let $G_2'$ be the connected subgroup of $G$ with Lie algebra ${\cal S}_2'$. Observe that ${\frak F}_2(p)=G_2'p$ and that for $g\in G$ the leaf ${\frak F}_2(gp)$ is an orbit of $gG_2'g^{-1}$. Since every leaf in ${\frak F}_2$ is closed in $M$, it follows that $G_2'$ is closed.

Further, the Jacobi identity applied to the triples $\{a_r,b_r,a_n\}$ and $\{a_r,b_r,b_n\}$ for $r=1,\dots,n-1$, together with (\ref{arbr}), (\ref{h2perph}), (\ref{identities2}), (\ref{thesame}) yields
\begin{equation}
\alpha\sigma=0,\quad \gamma\sigma=0.\label{idalphasigma}
\end{equation}

{\bf Case 1.1.1.a.} {\it Suppose that $G_1p$ is equivalent to one of $\BB^{n-1}$, $\CC\PP^{n-1}$.} In this case $\sigma\ne 0$, and (\ref{idalphasigma}) gives $\alpha=\gamma=0$. Then (\ref{h2perph}), (\ref{identities44}) imply that the Lie subalgebra ${\cal S}_2'$ is in fact an ideal in ${\frak g}$, and therefore $G_2'$ is normal in $G$. This yields that ${\frak F}_2$ is a holomorphic foliation. 

Since ${\frak g}={\cal S}_1\oplus{\cal S}_2'$ and $G_1$, $G_2'$ are closed, the group $G$ is a locally direct product of $G_1$ and $G_2'$. Let ${\mathscr T}:=G_1\cap G_2'$. Clearly, ${\mathscr T}$ is a discrete normal (hence central) subgroup of $G_1$. However, the center of each of $\hbox{Aut}(\BB^{n-1})$, $G(\CC\PP^{n-1})$ is trivial. Hence ${\mathscr T}$ is trivial and therefore $G=G_1\times G_2'$.

Let $K_2'$ be the ineffectivity kernel of the action of $G_2'$ on $G_2'p$. Then for every $g\in K_2'$ its differential $dg(p)$ at $p$ in the coordinates $(\xi_1,\dots,\xi_n)$ in $T_p(M)$ has the form
\begin{equation}
\left(
\begin{array}{ll}
A & 0\\
0 & 1
\end{array}
\right)\label{stabiliq}
\end{equation}
where $A\in U_{n-1}$ (note that $T_p(G_2'p)$ is given by $\{\xi_1=0,\dots,\xi_{n-1}=0\}$). Hence the element $g$ lies in $H^0\subset G_1$, and therefore is trivial. Thus, $K_2'$ is trivial.

We will now show that for every $q_1,q_2\in M$ the orbits $G_1q_1$ and $G_2'q_2$ intersect at exactly one point. Clearly, $G_1q_1\cap G_2'q_2$ is non-empty (see the proof of Theorem \ref{classtype2}). We will now show that for every $q\in M$ we have $G_1q\cap G_2'q=\{q\}$. Before proceeding, we remark that for $n\ge 3$ and for $n=2$ with $G_1p$ equivalent to $\BB^1$ this statement follows as in the proof of Theorem \ref{classtype2}. Below we give a proof that works for $n=2$ with $G_1p$ equivalent to $\CC\PP^1$ as well. Let for some $q\in M$ the intersection $G_1q\cap G_2'q$ contain a point $q'\ne q$. Let $g_1\in G_1$ be an element such that $g_1q=q'$. Clearly, $g_1$ preserves $G_2'q$. Since $g_1\in G_1$ and $G=G_1\times G_2'$, the element $g_1$ commutes with every element of $G_2'$. Consider the restriction $g_1':=g_1|_{G_2'q}$. Recall that $G_2'q$ is holomorphically equivalent to one of ${\cal P}_{>}$ , $\CC$, $\CC^*$, or an elliptic curve $\TT$ by means of a map $f$ that transforms $G_2'$ into one of $G({\cal P}_{>})$, $G_1(\CC)$, $\hbox{Aut}(\CC^*)^0$, or $\hbox{Aut}(\TT)^0$, respectively. Then $g_1'$ is transformed by $f$ into an automorphism of the corresponding curve that commutes with every automorphism in the corresponding group.

If $G_2'q$ is equivalent to ${\cal P}_{>}$ this implies that $g_1'$ is transformed into the identity, since the centralizer of $G({\cal P}_{>})$ in $\hbox{Aut}({\cal P}_{>})$ is trivial. Hence $g_1'$ acts on $G_2'q$ trivially, which contradicts our choice of $g_1$. If $G_2'q$ is equivalent to $\CC$, it follows that $g_1'\in G_1(\CC)$. Therefore, there exists $g_2\in G_2'$ such that $g_1'=g_2|_{G_2'q}$. Thus, $g:=g_1\circ g_2^{-1}$ acts trivially on $G_2'q$. Hence $dg(q)$ in suitable coordinates in $T_q(M)$ has the form (\ref{stabiliq}). Therefore, $g$ lies in a subgroup conjugate to $H^0$, hence $g\in G_1$. This implies that $g_2$ is trivial, which is impossible. If $G_2'q$ is equivalent to either $\CC^*$ or $\TT$, an analogous argument leads to a contradiction. Thus, we have shown that for every $q_1,q_2\in M$ the orbits $G_1q_1$ and $G_2'q_2$ intersect at exactly one point.

Let $F_1$ be a biholomorphic map from $G_1p$ onto one of $\BB^{n-1}$, $\CC\PP^{n-1}$ that transforms $G_1$ into one of $\hbox{Aut}(\BB^{n-1})$, $G(\CC\PP^{n-1})$, respectively, and $F_2$ a biholomorphic map from $G_2'p$ onto one of  ${\cal P}_{>}$ , $\CC$, $\CC^*$, or an elliptic curve $\TT$, that transforms $G_2'$ into one of $G({\cal P}_{>})$, $G_1(\CC)$, $\hbox{Aut}(\CC^*)^0$, or $\hbox{Aut}(\TT)^0$, respectively. For $q\in M$ let $p_1:=G_1p\cap G_2'q$ and $p_2:=G_1q\cap G_2'p$. Setting $F(q):=(F_1(p_1),F_2(p_2))$ we obtain a map from $M$ onto the corresponding direct product that transforms $G$ into the product of the corresponding groups. Since each foliation ${\frak F}_j$ is holomorphic, the map $F$ is holomorphic. Thus, we have obtained the direct products listed in (ia), (ib), (ic), (id), where $M'$ is either $\BB^{n-1}$ or $\CC\PP^{n-1}$.

{\bf Case 1.1.1.b.} {\it Suppose that $G_1p$ is equivalent to $\CC^{n-1}$.} For every $q\in M$ consider the set
$$
S_q:=\left\{x\in M: G_{1\,x}=G_{1\,q}\right\}.
$$
Since for two distinct points $x,y$ lying in the same $G_1$-orbit we have $G_{1\,x}\ne G_{1\,y}$, the set $S_q$ intersects every $G_1$-orbit in $M$ at exactly one point. Furthermore, for every $x\in S_q$ the set $S_q$ coincides near $x$ with the leaf ${\frak F}_2(x)$. Hence $S_q$ is the union of some leaves of ${\frak F}_2$. We will now prove that $S_q$ in fact consists of a single leaf, that is, $S_q={\frak F}_2(q)$. Indeed, it is sufficient to show that every leaf ${\frak F}_2(x)\subset S_q$ intersects every $G_1$-orbit. We need to prove that the set of orbits that ${\frak F}_2(x)$ intersects is both open and closed. The openness is obvious. To show the closedness, let $G_1y_j$ be a sequence of $G_1$-orbits that intersect ${\frak F}_2(x)$ and accumulate to an orbit $G_1y$. Suppose that $G_1y$ does not intersect ${\frak F}_2(x)$. Let $y':=S_q\cap G_1y$. Then ${\frak F}_2(y')\subset S_q$ and ${\frak F}_2(y')\ne{\frak F}_2(x)$. Clearly, ${\frak F}_2(y')$ intersects every orbit $G_1y_j$ if $j$ is sufficiently large, and therefore $S_q$ intersects $G_1y_j$ at more than one point. This contradiction proves that ${\frak F}_2(x)$ intersects every $G_1$-orbit, hence $S_q={\frak F}_2(q)$ for every $q\in M$.  

Let $F_1: G_1p\ra\CC^{n-1}$ be a biholomorphic map that takes $p$ into the origin and transforms $G_1$ into $G(\CC^{n-1})$, and $F_2$ a biholomorphic map from $G_2'p$ onto one of  ${\cal P}_{>}$ , $\CC$, $\CC^*$, or an elliptic curve $\TT$, that transforms $G_2'/K_2'$ into one of $G({\cal P}_{>})$, $G_1(\CC)$, $\hbox{Aut}(\CC^*)^0$, or $\hbox{Aut}(\TT)^0$, respectively, where $K_2'$ denotes, as before, the ineffectivity kernel of the action of $G_2'$ on $G_2'p$. For $q\in M$ let $p_1:=G_1p\cap G_2'q$ and $p_2:=G_1q\cap G_2'p$. Setting $F(q):=(F_1(p_1),F_2(p_2))$ we obtain a diffeomorphism from $M$ onto a direct product $\CC^{n-1}\times M''$, where $M''$ is one of ${\cal P}_{>}$, $\CC$, $\CC^*$, $\TT$.

We will now show that $F$ is holomorphic. Let $\Pi_1: M\ra G_1p$ be the projection to $G_1p$ along the leaves of ${\frak F}_2$ and $\Pi_2: M\ra G_2'p$ be the projection to $G_2'p$ along the leaves of ${\frak F}_1$. In order to show that $F$ is holomorphic, we must prove that each $\Pi_j$ is holomorphic. By \cite{Ch} it is sufficient to show that each $\Pi_j$ is separately holomorphic on every leaf of each of ${\frak F}_1$, ${\frak F}_2$. Thus, we need to show that $\Pi_j$ is holomorphic on every leaf of ${\frak F}_j$, for $j=1,2$. Since the foliation ${\frak F}_1$ is holomorphic, the map $\Pi_2$ is holomorphic. 

We will now show that $\Pi_1$ is holomorphic on every leaf of ${\frak F}_1$. Fix $q\in M$ and let $F_0:G_1q\ra\CC^{n-1}$ be a biholomorphic map that transforms $G_1$ into $G(\CC^{n-1})$. Therefore, there exists an automorphism $\varphi$ of the group $G(\CC^{n-1})$ such that $\Pi_1|_{G_{{}_1}q}=F_1^{-1}\circ f_{\varphi} \circ F_0$, where $f_{\varphi}$ is the diffeomorphism of $\CC^{n-1}$ that maps $z\in\CC^{n-1}$ to the point $w\in\CC^{n-1}$ such that ${\cal I}_{w}=\varphi({\cal I}_z)$, with ${\cal I}_z$ being the isotropy subgroup of $z$ with respect to the action of $G(\CC^{n-1})$ on $\CC^{n-1}$.

Every automorphism of $G(\CC^{n-1})$ has either the form
\begin{equation}
Uz+a \mapsto g_0Ug_0^{-1}(z-b)+b +La,\label{phiform1}
\end{equation}
or the form 
\begin{equation}
Uz+a \mapsto g_0\overline{U}g_0^{-1}(z-b)+b +\overline{La},\label{phiform2}
\end{equation}
for some $g_0\in U_{n-1}$, $b\in\CC^{n-1}$ and $L\in GL_{n-1}(\CC)$ (here $U\in U_{n-1}$, $a\in\CC^{n-1}$). If $\varphi$ has the form (\ref{phiform1}), then $f_{\varphi}\in\hbox{Aut}(\CC^{n-1})$, and thus $\Pi_1|_{G_{{}_1}q}$ is holomorphic. If $\varphi$ has the form (\ref{phiform2}), then $f_{\varphi}$ is antiholomorphic and so is $\Pi_1|_{G_{{}_1}q}$. It is straightforward to show that the set of points $x\in M$ for which $\Pi_1|_{G_{{}_1}x}$ is anti-holomorphic is both open and closed, and thus coincides with $M$. Since $\Pi_1|_{G_{{}_1}p}$ is clearly holomorphic, we obtain  a contradiction. Hence, $\varphi$ in fact has the form (\ref{phiform1}) and $\Pi_1|_{G_{{}_1}q}$ is holomorphic. Thus, we have constructed a biholomorphic map from $M$ onto a direct product $\CC^{n-1}\times M''$, where $M''$ is one of ${\cal P}_{>}$, $\CC$, $\CC^*$, $\TT$.

We will now describe the subgroup of $\hbox{Aut}(\CC^{n-1}\times M'')$ into which the map $F$ transforms the group $G$. Let $G^F:=F\circ G\circ F^{-1}$ and $G_1^F:=F\circ G_1\circ F^{-1}$. For $\zeta_0\in\CC^{n-1}$ let $L_{\zeta_0}:=\{(\zeta,\xi)\in\CC^{n-1}\times M'':\zeta=\zeta_0\}$ and for $\xi_0\in M''$ let $L_{\xi_0}'':=\{(\zeta,\xi)\in\CC^{n-1}\times M'': \xi=\xi_0\}$. By construction, $L_0=F_2(G_2'p)$. Let $\hat\xi:=F_2(p)$; then we have $L_{\hat\xi}''=F_1(G_1p)$. For every $\xi_0\in M''$ the subset $L_{\xi_0}''$ is naturally identified with $\CC^{n-1}$, and, upon this identification, we have $G_1^F|_{L_{\xi_0}''}=f_{\xi_0}\circ G(\CC^{n-1})\circ f^{-1}_{\xi_0}$ for some $f_{\xi_0}\in\hbox{Aut}(\CC^{n-1})$, where $f_{\hat\xi}=\hbox{id}$.

Let $\Phi_{\xi_0}$ be the map
$$
\Phi_{\xi_0}:\, G_1^F|_{L_{\hat\xi}''}\ra G_1^F|_{L_{\xi_0}''},\,
g|_{L_{\hat\xi}''}\mapsto g|_{L_{\xi_0}''},\,\hbox{for $g\in G_1^F$}.
$$
Upon the identification of each of $L_{\hat\xi}''$, $L_{\xi_0}''$ with $\CC^{n-1}$, the map $\Phi_{\xi_0}$ takes the form
$$
T\mapsto f_{\xi_0}\circ\varphi(T) \circ f^{-1}_{\xi_0},
$$
where $T\in G(\CC^{n-1})$ and $\varphi$ is an automorphism of $G(\CC^{n-1})$ (see (\ref{phiform1}), (\ref{phiform2})). By the construction of the map $F$, for every point $(\zeta,\hat\xi)\in L_{\hat\xi}''$, every element $g\in G_1^F$ that fixes this point also fixes the point $(\zeta,\xi_0)\in L_{\xi_0}''$. It then follows that $\Phi_{\xi_0}$ is in fact the identity map, that is, $G_1^F$ consists of all maps of the form
\begin{equation}
\begin{array}{lll}
\zeta & \mapsto & U\zeta+a,\\
\xi & \mapsto & \xi,
\end{array}\label{g1transform}
\end{equation}
where $U\in U_{n-1}$, $a\in\CC^{n-1}$.

Let $\hat{\frak g}$ be the Lie algebra of holomorphic vector fields on $\CC^{n-1}\times M''$ arising from the action of $G^F$. Denote by $\hat a_r,\hat b_r,\hat a_n',\hat b_n',\hat t_0$ the holomorphic vector fields on $\CC^{n-1}\times M''$ arising from $a_r,b_r,a_n',b_n',t_0$, respectively, for $r=1,\dots,n-1$. Since $F$ transforms $H^0$ into the group of maps of the form (\ref{g1transform}) with $a=0$, we have
$$
\hat t_0=i\mu'\sum_{s=1}^{n-1}\zeta_s\,\partial/\partial \zeta_s,
$$
for some $\mu'\in\RR^*$. Considering the commutators $[\hat t_0,\hat a_r]$, $[\hat t_0,\hat b_r]$ and using (\ref{t0}), (\ref{identities1}) we see that
$$
\begin{array}{lll}
\hat a_r&=&\displaystyle\sum_{s=1}^{n-1}A_r^s\,\partial/\partial\zeta_s,\vspace{0.1cm}\\
\hat b_r&=\pm i&\displaystyle\sum_{s=1}^{n-1}A_r^s\,\partial/\partial\zeta_s,
\end{array}
$$
where $(A_r^s)_{r,s=1,\dots,n-1}\in GL_{n-1}(\CC)$. Next, using the identities $[\hat t_0,\hat a_n']=0$, $[\hat t_0,\hat b_n']=0$ (see (\ref{h2perph})), from (\ref{identities44}) we obtain
\begin{equation}
\begin{array}{lll}
\hat a_n'&=&\displaystyle\sum_{s=1}^{n-1}f_s\,\partial/\partial \zeta_s+a_n'',\\
\vspace{-0.3cm}\\
\hat b_n'&=&\displaystyle\sum_{s=1}^{n-1}g_s\,\partial/\partial \zeta_s+b_n'',
\end{array}\label{an''bn''}
\end{equation}
where $f_s$, $g_s$ are linear functions of $\zeta_1,\dots,\zeta_{n-1}$, and $a_n''$, $b_n''$ are holomorphic vector fields on $M''$. Moreover, it follows from the first identity in (\ref{id4}) that
\begin{equation}
\rho f_s+\nu g_s\equiv 0,\quad s=1,\dots,n-1.\label{rhofnug}
\end{equation}

Recall that $F$ maps $G_2'p$ biholomorphically onto $L_0$ and transforms $G_2'/K_2'$ into one of $G({\cal P}_{>})$, $G_1(\CC)$, $\hbox{Aut}(\CC^*)^0$, $\hbox{Aut}(\TT)^0$ (where we identify $L_0$ with one of ${\cal P}_{>}$, $\CC$, $\CC^*$, $\TT$, respectively). Therefore, $a_2''|_{L_0}$ and $b_2''|_{L_0}$ span the Lie algebra of holomorphic vector fields on $L_0$ arising from the action of one of these groups.
 
If $G_2'p$ is equivalent to ${\cal P}_{>}$, then we have
$$
a_n''=(x\xi+iy)\,\partial/\partial\xi,\quad b_n''=(u\xi+iv)\,\partial/\partial\xi,
$$
where $x,y,u,v\in\RR$, $xv-yu\ne 0$. Considering $[\hat a_n',\hat b_n']$ and using (\ref{identities44}), we obtain
\begin{equation}
\rho x+\nu u=0,\quad xv-yu=\rho y+\nu v.\label{relparamm}
\end{equation}
Since $\rho^2+\nu^2>0$ (the group $G({\cal P}_{>})$ is not abelian), it follows that $x=\nu$, $u=-\rho$. Now (\ref{an''bn''}), (\ref{rhofnug}), (\ref{relparamm}) yield
$$
\rho\hat a_n'+\nu \hat b_n'=i(\rho y+\nu v)\,\partial/\partial \xi,
$$
and therefore the algebra $\hat{\frak g}$ contains the vector field $i\,\partial/\partial \xi$. Hence $\hat{\frak g}$ is spanned by the vector fields arising from the action of $G_1^F$, the vector field $i\,\partial/\partial \xi$, and a vector field of the form
$$
\displaystyle\sum_{s=1}^{n-1}h_s\,\partial/\partial \zeta_s+\xi\,\partial/\partial \xi,
$$
where $h_s$ are linear functions of $\zeta_1,\dots,\zeta_{n-1}$. The one-parameter subgroup of $G^F$ arising from such a vector field has the form
$$
\begin{array}{lll}
\zeta&\mapsto&e^{R\tau}\zeta,\\
\xi&\mapsto&e^{\tau}\xi,
\end{array}
$$
where $R\in{\frak{gl}}_{n-1}$ and $\tau\in\RR$. Since $G_1^F$ is normal in $G^F$, it follows that the matrix $e^{R\tau}$ lies in the normalizer of $U_{n-1}$ in $GL_{n-1}(\CC)$ for all $\tau\in\RR$. Therefore, there exists $T\in\RR$ such that $e^{R\tau}\in e^{T\tau}\cdot U_{n-1}$ for all $\tau$. We have now shown that $G^F$ is the group $G(\CC^{n-1})\times G({\cal P}_{>})$ if $T=0$, and the group of maps of the form (\ref{groupv}) if $T\ne 0$, where $z'$ and $z_n$ are replaced with $\zeta$ and $\xi$, respectively. Thus, we have obtained Example (id) for $M'=\CC^{n-1}$, as well as Example (v).

If $G_2'p$ is equivalent to $\CC$, then $a_n''$ and $b_n''$ are constant vector fields. Therefore, $\hat{\frak g}$ is spanned by the vector fields arising from the action of $G_1^F$, and some vector fields
\begin{equation}
\begin{array}{l}
\displaystyle\sum_{s=1}^{n-1}\hat f_s\,\partial/\partial \zeta_s+\partial/\partial \xi,\\
\vspace{-0.3cm}\\
\displaystyle\sum_{s=1}^{n-1}\hat g_s\,\partial/\partial \zeta_s+i\,\partial/\partial \xi,
\end{array}\label{hatfhatg}
\end{equation}
where $\hat f_s$, $\hat g_s$ are linear functions of $\zeta_1,\dots,\zeta_{n-1}$. The one-parameter subgroups of $G^F$ arising from these vector fields have the forms
$$
\begin{array}{lll}
\zeta&\mapsto&e^{R_1\tau}\zeta,\\
\xi&\mapsto&\xi+\tau,
\end{array}\quad
\begin{array}{lll}
\zeta&\mapsto&e^{R_2\tau}\zeta,\\
\xi&\mapsto&\xi+i\tau,
\end{array}
$$
where $R_j\in{\frak{gl}}_{n-1}$ for $j=1,2$, and $\tau\in\RR$. Using the normality of $G_1^F$ in $G^F$ as above, we see that $G^F$ is the group of maps of the form (\ref{g3cns}) for some $T\in\CC$, where $z'$ and $z_n$ are replaced with $\zeta$ and $\xi$, respectively. Thus, we have obtained Example (ia) for $M'=\CC^{n-1}$ and Example (iiia).

If $G_2'p$ is equivalent to $\CC^*$, then we have
$$
a_n''=u\xi\,\partial/\partial\xi,\quad b_n''=v\xi\,\partial/\partial\xi,
$$       
where $u,v\in\CC^*$. Therefore, $\hat{\frak g}$ is spanned by the vector fields arising from the action of $G_1^F$, and some vector fields
$$
\begin{array}{l}
\displaystyle\sum_{s=1}^{n-1}\hat f_s\,\partial/\partial \zeta_s+\xi\,\partial/\partial \xi,\\
\vspace{-0.3cm}\\
\displaystyle\sum_{s=1}^{n-1}\hat g_s\,\partial/\partial \zeta_s+i\xi\,\partial/\partial \xi,
\end{array}
$$
where $\hat f_s$, $\hat g_s$ are as in (\ref{hatfhatg}). The one-parameter subgroups of $G^F$ arising from these vector fields have the forms
$$
\begin{array}{lll}
\zeta&\mapsto&e^{R_1\tau}\zeta,\\
\xi&\mapsto&e^{\tau}\xi,
\end{array}\quad
\begin{array}{lll}
\zeta&\mapsto&e^{R_2\tau}\zeta,\\
\xi&\mapsto&e^{i\tau}\xi,
\end{array}
$$
where $R_j\in{\frak{gl}}_{n-1}$ for $j=1,2$, and $\tau\in\RR$. Arguing as above, we see that $G^F$ is the group of maps of the form
$$
\begin{array}{lll}
\zeta&\mapsto & e^{\hbox{\tiny Re}\,(T b)}U\zeta+a,\\
\xi&\mapsto & e^{b}\xi,
\end{array}
$$
for some $T\in\CC$, where $U\in U_{n-1}$, $a\in\CC^{n-1}$, $b\in\CC$. The action of this group on $\CC^{n-1}\times\CC^*$ is proper if and only if $T\in\RR$. Thus, we have obtained Example (ib) for $M'=\CC^{n-1}$ and Example (iiib).

If $G_2'p$ is equivalent to $\TT$, then the lifts of $a_n''$ and $b_n''$ to $\CC$ are constant vector fields. Therefore, $\hat{\frak g}$ is spanned by the vector fields arising from the action of $G_1^F$, and some vector fields
$$
\begin{array}{l}
\displaystyle\sum_{s=1}^{n-1}\hat f_s\,\partial/\partial \zeta_s+a_n''',\\
\vspace{-0.3cm}\\
\displaystyle\sum_{s=1}^{n-1}\hat g_s\,\partial/\partial \zeta_s+b_n''',
\end{array}
$$
where $\hat f_s$, $\hat g_s$ are as in (\ref{hatfhatg}), and where $a_n'''$ and $b_n'''$ are vector fields on $\TT$ whose lifts to $\CC$ are $\partial/\partial \xi$ and $i\partial/\partial \xi$, respectively. The one-parameter subgroups of $G^F$ arising from these vector fields have the forms
$$
\begin{array}{lll}
\zeta&\mapsto&e^{R_1\tau}\zeta,\\
\xi&\mapsto&\Psi_{\tau}(\xi),
\end{array}\quad
\begin{array}{lll}
\zeta&\mapsto&e^{R_2\tau}\zeta,\\
\xi&\mapsto&\Psi_{i\tau}(\xi),
\end{array}
$$
where $R_j\in{\frak{gl}}_{n-1}$ for $j=1,2$, $\tau\in\RR$, and for $b\in\CC$ the map $\Psi_b$ is an element of $\hbox{Aut}(\TT)^0$ induced by the translation by $b$ on $\CC$. Arguing as earlier, we see that $G^F$ is the group of maps of the form
$$
\begin{array}{lll}
\zeta&\mapsto & e^{\hbox{\tiny Re}\,(T b)}U\zeta+a,\\
\xi&\mapsto & \Psi_b(\xi),
\end{array}
$$
for some $T\in\CC$, where $U\in U_{n-1}$, $a\in\CC^{n-1}$, $b\in\CC$. The action of such a group on $\CC\times\TT$ is proper only for $T=0$. Thus, we have obtained Example (ic) for $M'=\CC^{n-1}$.

{\bf Case 1.1.2.} {\it Suppose that $k_2\ne 0$.} Since $[a_n,b_n]$ is $\hbox{Ad}(H)$-invariant, we have
\begin{equation}
[a_n,b_n]=\eta t_0\label{anbnk2ne0}
\end{equation}
for some $\eta\in\RR$ (cf. (\ref{id3}), (\ref{id3prime})). Further, for every $r=1,\dots,n-1$ we have
\begin{equation}
[t_r,a_n]=\frac{k_2}{k_1}\mu b_n,\quad [t_r,b_n]=-\frac{k_2}{k_1}\mu a_n,\label{tran}
\end{equation} 
Next, consider	 the collections of commutators
$$
{\cal C}_r:=\left\{[a_r,a_n],\,[a_r,b_n],\,[b_r,a_n]\,[b_r,b_n]\right\},\quad r=1,\dots,n-1.
$$
We expand each commutator in ${\cal C}_r$ with respect to the decomposition ${\frak g}={\frak h}_1^{\perp}+{\frak h}_2^{\perp}+{\frak h}$ and see what terms can be present in the expansion.

{\bf Case 1.1.2.a.} {\it Suppose that $n\ge 3$.} In this case we apply transformations from $\hbox{Ad}(H_s)$ for $1\le s \le n-1$, $s\ne r$, to every element of ${\cal C}_r$. This immediately yields that the element's expansion does not have terms containing $a_l$, $b_l$ for $l\ne s,n$, and
\begin{equation}
\begin{array}{llrl}
[a_r,a_n]&\equiv&\alpha_r a_n+\beta_r b_n\,&\hbox{(mod ${\frak h}_1^{\perp}+{\frak h}$)},\\

[a_r,b_n]&\equiv&-\beta_r a_n+\alpha_r b_n\,&\hbox{(mod ${\frak h}_1^{\perp}+{\frak h}$)},\\

[b_r,a_n]&\equiv&\gamma_r a_n+\delta_r b_n\,&\hbox{(mod ${\frak h}_1^{\perp}+{\frak h}$)},\\

[b_r,b_n]&\equiv&-\delta_r a_n+\gamma_r b_n\,&\hbox{(mod ${\frak h}_1^{\perp}+{\frak h}$)},
\end{array}\label{relmod}
\end{equation}
for some $\alpha_r,\beta_r,\gamma_r,\delta_r\in\RR$.

{\bf Case 1.1.2.a.1.} {\it Suppose that $n\ge 4$.} In this case by varying $s$ in the above argument we obtain that every element of ${\cal C}_r$ has a zero projection to ${\frak h}_1^{\perp}$ for $r=1,\dots,n-1$. Next, we apply transformations from $\hbox{Ad}(Z)$ to each commutator in ${\cal C}_r$. 

Firstly, assume in addition that $k_1\ne\pm k_2(n-1)$ and $k_1\ne\pm 2k_2(n-1)$. For such $k_1,k_2$ a straightforward calculation implies that every element of ${\cal C}_r$ is in fact zero. Applying the Jacobi identity to the triple $\{a_r,b_r,a_n\}$ for $r=1,\dots,n-1$ and using (\ref{arbr}), (\ref{tran}), we see that $\sigma=0$. Therefore $G_1p$ is holomorphically equivalent to $\CC^{n-1}$ by means of a map that transforms $G_1/K_1$ into $G(\CC^{n-1})$. Next, applying the Jacobi identity to the triple $\{a_r,a_n,b_n\}$ and using (\ref{identities1}), (\ref{anbnk2ne0}), we obtain that $\eta=0$ which implies that $G_2p$ is holomorphically equivalent to $\CC$ by means of a map that transforms $G_2/K_2$ into $G(\CC)$. Therefore, ${\frak g}$ is isomorphic to the Lie algebra of the group $G_{k_1,k_2}(\CC^n)$ defined in Example (vi).

Every connected Lie group with Lie algebra isomorphic to that of the group $G_{k_1,k_2}(\CC^n)$ and containing an $(n-1)^2$-dimensional compact subgroup is isomorphic either to the group $G_{k_1,k_2}(\CC^n)$ or to its finitely-sheeted cover. Let ${\frak G}$ be one of these groups. Then every $(n-1)^2$-dimensional compact subgroup ${\frak K}\subset{\frak G}$ is maximal compact. Furthermore, ${\frak G}$ acts effectively on ${\frak G}/{\frak K}$ only if ${\frak G}=G_{k_1,k_2}(\CC^n)$. Therefore, the group $G$ is isomorphic to $G_{k_1,k_2}(\CC^n)$, and this isomorphism induces a diffeomorphism from $G/H$ onto $G_{k_1,k_2}(\CC^n)/K$, where $K$ is a maximal compact subgroup of $G_{k_1,k_2}(\CC^n)$. Next, $G_{k_1,k_2}(\CC^n)/K$ is $G_{k_1,k_2}(\CC^n)$-equivariantly diffeomorphic to $\CC^n$. The linearization of the action of $K$ on $\CC^n$ at its fixed point $p_K$ has exactly two non-trivial proper real invariant subspaces $L_1,L_2\subset T_{p_K}(\CC^n)$. Their dimensions are $2(n-1)$ and 2. Observe that there are exactly two non-trivial proper $LH^0$-invariant real subspaces in $T_p(M)$, and that these subspaces are in fact complex and have complex dimensions $(n-1)$ and 1. Therefore, $L_1$ and $L_2$ are complex subspaces of $T_{p_K}(\CC^n)$ with respect to the complex structure induced on $\CC^n$ by that of $M$. We will now determine all $G_{k_1,k_2}(\CC^n)$-invariant complex structures on $\CC^n$ for which $L_1$ and $L_2$ are complex.

Let $P$ be a Lie group with Lie algebra ${\frak p}$ and $Q$ a closed subgroup in $P$ with Lie algebra ${\frak q}$. If the manifold $P/Q$ admits a $P$-invariant complex structure, then there is a linear endomorphism $J$ of ${\frak p}$ satisfying the following conditions:
\begin{equation}
\begin{array}{l}
(1)\,J{\frak q}\subset{\frak q},\\
\vspace{0.1cm}\\
(2)\,J^2x\equiv -x\,\hbox{(mod ${\frak q}$) for all $x\in{\frak p}$},\\
\vspace{0.1cm}\\
(3)\,J[x,y]\equiv [x,Jy]\,\hbox{(mod ${\frak q}$) for all $x\in{\frak q}$, $y\in{\frak p}$},\\
\vspace{0.1cm}\\
(4)\,[Jx,Jy] \equiv [x,y]+J[Jx,y]+J[x,Jy]\,\hbox{(mod ${\frak q}$) for all $x,y\in{\frak p}$},
\end{array}\label{invarcomplstructure}
\end{equation}
and the map $J$ induces a linear transformation $\tilde J$ of ${\frak p}/{\frak q}$, which, upon the canonical identification of ${\frak p}/{\frak q}$ with the tangent space $T_Q(P/Q)$, coincides with the operator of complex structure on $T_Q(P/Q)$ (see \cite{Kosz}). For $P=G_{k_1,k_2}(\CC^n)$, $Q=K$ and $\tilde J$ leaving invariant the subspaces $L_1$ and $L_2$, conditions (1)--(3) of (\ref{invarcomplstructure}) lead to the following four structures: (a) the standard complex structure on $\CC^n$, (b) the complex structure conjugate to the standard one, (c) the complex structures obtained by conjugating the standard complex structure in either $z'$ or $z_n$ (note that condition (4) does not contribute any additional information in this situation). The manifold in case (b) is holomorphically equivalent to $\CC^n$ by means of a map that preserves the group $G_{k_1,k_2}(\CC^n)$. On the other hand, each of the two manifolds in case (c) is holomorphically equivalent to $\CC^n$ by means of a map that transforms the group $G_{k_1,k_2}(\CC^n)$ into $G_{k_1,-k_2}(\CC^n)$. Thus, we have obtained Example (vi) for $n\ge 4$ and $k_1\ne\pm k_2(n-1)$, $k_1\ne\pm 2k_2(n-1)$.

Assume next that $k_1=\pm 2k_2(n-1)$, that is, $k_1=2(n-1)$, $k_2=\pm1$. In this case the application of transformations from $\hbox{Ad}(Z)$ to the commutators in ${\cal C}_r$ does not immediately imply that they are all equal to zero; instead, after a short calculation it leads to the following relations 
$$
\begin{array}{llr}
[a_r,a_n]&=&\alpha'_r a_n+\beta'_r b_n,\\

[a_r,b_n]&=&\beta'_r a_n-\alpha'_r b_n,\\

[b_r,a_n]&=&\mp\beta'_r a_n\pm\alpha'_r b_n,\\

[b_r,b_n]&=&\pm\alpha'_r a_n\pm\beta'_r b_n,
\end{array}
$$
for some $\alpha'_r,\beta'_r\in\RR$, $r=1,\dots,n-1$. Together with (\ref{relmod}) these relations do in fact imply that all commutators in ${\cal C}_r$ are in fact equal to zero, for $r=1,\dots,n-1$. Arguing as earlier, we are now led to Example (vi) for $n\ge 4$ and $k_1=\pm 2k_2(n-1)$.

Assume finally that $k_1=\pm k_2(n-1)$, that is, $k_1=n-1$, $k_2=\pm1$. In this case the application of transformations from $\hbox{Ad}(Z)$ to the commutators in ${\cal C}_r$ only yields
\begin{equation}
\begin{array}{lll}
{\cal C}_r\subset{\frak h},\\
\vspace{-0.3cm}\\

[a_r,a_n]&=&\pm[b_r,b_n],\\
\vspace{-0.3cm}\\

[a_r,b_n]&=&\mp[b_r,a_n].
\end{array}\label{moreident}
\end{equation}
Applying, as earlier, the Jacobi identity to the triple $\{a_r,b_r,a_n\}$ and using (\ref{arbr}), (\ref{tran}) we obtain 
\begin{equation}
\frac{k_2}{k_1}\mu\sigma b_n+[[b_r,a_n],a_r]+[[a_n,a_r],b_r]=0,\quad r=1,\dots,n-1.\label{jacin1}
\end{equation}
Due to (\ref{hperph}), (\ref{moreident}), the last two terms in the left-hand side of (\ref{jacin1}) lie in ${\frak h}_1^{\perp}$, which implies that $\sigma=0$.
 
For $r=1,\dots,n-1$ let $\hat H_r$ be the subgroup of $H^0$ represented in $LH^0$ by transformations that preserve each of the variables $\xi_r$, $\xi_n$. Since $\hbox{Ad}(\hat H_r)$ preserves each element of ${\cal C}_r$, we have
$$
\begin{array}{lll}
[a_r,a_n]&=&\displaystyle\sum_{s=1}^{n-1}\nu_r^st_s,\\
\vspace{-0.3cm}\\

[b_r,a_n]&=&\displaystyle\sum_{s=1}^{n-1}\rho_r^st_s,
\end{array}
$$
for some $\nu_r^s,\rho_r^s\in\RR$ (see (\ref{frakA})). Now (\ref{identities1}), (\ref{jacin1}) imply that $\nu_r^r=0$, $\rho_r^r=0$ for $r=1,\dots,n-1$ (recall that $\sigma=0$). Next, from (\ref{identities1}), (\ref{PrPs}), (\ref{moreident}) and the Jacobi identity applied to the triples $\{a_r,a_s,a_n\}$, $\{b_r,b_s,a_n\}$ where $s=1,\dots,n-1$, $s\ne r$, we see that $\nu_r^s=\rho_r^s=0$. From relations (\ref{moreident}) we now see that for $k_1=\pm k_2(n-1)$ every commutator in ${\cal C}_r$ is zero for $r=1,\dots,n-1$. Arguing as earlier, we obtain Example (vi) for $n\ge 4$ and $k_1=\pm k_2(n-1)$.

{\bf Case 1.1.2.a.2.} {\it Suppose that $n=3$.} We first assume that $k_1\ne\pm k_2$. In this case, the application of transformations from $\hbox{Ad}(H_s)$ for $s\ne r$ to the elements of ${\cal C}_r$ yields  
$$
\begin{array}{llrl}
[a_r,a_3]&\equiv&\alpha_r a_3+\beta_r b_3\,&\hbox{(mod ${\frak h}$)},\\

[a_r,b_3]&\equiv&-\beta_r a_3+\alpha_r b_3\,&\hbox{(mod ${\frak h}$)},\\

[b_r,a_3]&\equiv&\gamma_r a_3+\delta_r b_3\,&\hbox{(mod ${\frak h}$)},\\

[b_r,b_3]&\equiv&-\delta_r a_3+\gamma_r b_3\,&\hbox{(mod ${\frak h}$)},
\end{array}
$$
for some $\alpha_r,\beta_r,\gamma_r,\delta_r\in\RR$, $r=1,2$ (cf. (\ref{relmod})). If, furthermore, we have $k_1\ne\pm 2k_2$, then the arguments used in Case 1.1.2.a.1 work for this case as well and lead to Example (vi) for $n=3$ and $k_1\ne\pm k_2$, $k_1\ne\pm 2k_2$.

Let $k_1=\pm 2k_2$, that is, $k_1=2$, $k_2=\pm 1$. As in Case 1.1.2.a.1, the application of transformations from $\hbox{Ad}(Z)$ to the commutators in ${\cal C}_r$ yields
\begin{equation}
\begin{array}{lll}
\alpha_r=0,&&\beta_r=0,\\
\gamma_r=0,&&\delta_r=0,\\
\vspace{-0.3cm}\\

[a_r,a_3]&=&\pm[b_r,b_3],\\
\vspace{-0.3cm}\\

[a_r,b_3]&=&\mp[b_r,a_3].
\end{array}\label{moreident3}
\end{equation}
(cf. (\ref{moreident})). We now apply transformations from $\hbox{Ad}(H_r)$ to the commutators in ${\cal C}_r$, for $r=1,2$. Writing explicitly the action of $\hbox{Ad}(H_r)$ on ${\frak h}$ and using (\ref{moreident3}), it is straightforward to see that all commutators in ${\cal C}_r$ are in fact zero, for $r=1,2$. We then obtain Example (vi) for $n=3$ and $k_1=\pm 2k_2$.

Assume now that $k_1=\pm k_2$, that is, $k_1=1$, $k_2=\pm 1$. In this case, the application of transformations from $\hbox{Ad}(H_s)$ for $s\ne r$ to the elements of ${\cal C}_r$ yields  
$$
\begin{array}{llrl}
[a_1,a_3]&\equiv&\tilde\alpha_1 a_2+\tilde\beta_1 b_2+\alpha_1 a_3+\beta_1 b_3\,&\hbox{(mod ${\frak h}$)},\\

[a_1,b_3]&\equiv&\mp\tilde\beta_1a_2\pm\tilde\alpha_1b_2
-\beta_1 a_3+\alpha_1 b_3\,&\hbox{(mod ${\frak h}$)},\\

[b_1,a_3]&\equiv&\tilde\gamma_1a_2+\tilde\delta_1b_2+
\gamma_1 a_3+\delta_1 b_3\,&\hbox{(mod ${\frak h}$)},\\

[b_1,b_3]&\equiv&\mp\tilde\delta_1a_2\pm\tilde\gamma_1b_2
-\delta_1 a_3+\gamma_1 b_3\,&\hbox{(mod ${\frak h}$)},\\

[a_2,a_3]&\equiv&\tilde\alpha_2 a_1+\tilde\beta_2 b_1+\alpha_2 a_3+\beta_2 b_3\,&\hbox{(mod ${\frak h}$)},\\

[a_2,b_3]&\equiv&\mp\tilde\beta_2a_1\pm\tilde\alpha_2b_1
-\beta_2 a_3+\alpha_2 b_3\,&\hbox{(mod ${\frak h}$)},\\

[b_2,a_3]&\equiv&\tilde\gamma_2a_1+\tilde\delta_2b_1+
\gamma_2 a_3+\delta_2 b_3\,&\hbox{(mod ${\frak h}$)},\\

[b_2,b_3]&\equiv&\mp\tilde\delta_2a_1\pm\tilde\gamma_2b_1
-\delta_2 a_3+\gamma_2 b_3\,&\hbox{(mod ${\frak h}$)},\\
\end{array}
$$
for some $\alpha_r,\tilde\alpha_r,\beta_r,\tilde\beta_r,\gamma_r,\tilde\gamma_r,\delta_r,\tilde\delta_r\in\RR$, $r=1,2$. Applying next transformations from $\hbox{Ad}(Z)$ to the commutators in ${\cal C}_r$ gives that all these commutators have a zero projection to ${\frak h}$ and that 
$$
\begin{array}{l}
\alpha_r=\beta_r=\gamma_r=\delta_r=0,\\
\tilde\gamma_r=\tilde\beta_r,\,\tilde\delta_r=-\tilde\alpha_r.
\end{array}
$$ 
for $r=1,2$.

If $k_1=1$, $k_2=-1$, then, applying condition (4) of (\ref{invarcomplstructure}) to the $G$-invariant complex structure on $G/H$ and the vectors $x=a_r$, $y=a_3$, and using the fact that the pull-back to ${\frak h}^{\perp}$ of the operator of complex structure on $T_H(G/H)$ is given by (\ref{j0}) in the basis $\{a_k,b_k\}_{k=1,2,3}$, we obtain that $\tilde\alpha_r=0$, $\tilde\beta_r=0$, for $r=1,2$. Thus, in this case all commutators in ${\cal C}_r$ are zero for $r=1,2$, and we obtain Example (vi) for $n=3$ and $k_1=1$, $k_2=\pm 1$.

Assume now that $k_1=1$, $k_2=1$. Let $\hat H$ be the subgroup of $H^0$ represented in $LH^0$ by transformations preserving the variable $\xi_3$ and acting on ${\cal L}_1(p)$ by matrices from the group $SU_2$. Applying transformations from $\hbox{Ad}(\hat H)$ to the commutators in ${\cal C}_r$, we see that $\tilde\alpha_2=-\tilde\alpha_1$ and $\tilde\beta_2=-\tilde\beta_1$. We set: $\alpha:=\tilde\alpha_1$, $\beta:=\tilde\beta_1$. If $\alpha=0$, $\beta=0$, then, arguing as earlier, we obtain Example (vi) for $n=3$ and $k_1=1$, $k_2=\pm 1$. Let $\alpha^2+\beta^2>0$. We have
\begin{equation}
\begin{array}{llrl}
[a_1,a_3]&=&\alpha a_2+\beta b_2,\\

[a_1,b_3]&=&-\beta a_2+\alpha b_2,\\

[b_1,a_3]&=&\beta a_2-\alpha b_2,\\

[b_1,b_3]&=&\alpha a_2+\beta b_2,\\

[a_2,a_3]&=&-\alpha a_1-\beta b_1,\\

[a_2,b_3]&=&\beta a_1-\alpha b_1,\\

[b_2,a_3]&=&-\beta a_1+\alpha b_1,\\

[b_2,b_3]&=&-\alpha a_1-\beta b_1.
\end{array}\label{reln3}
\end{equation} 
We now apply the Jacobi identity to the triples $\{a_1,a_2,a_3\}$, $\{a_1,b_2,a_3\}$. It follows from (\ref{PrPs}), (\ref{arbr}), (\ref{reln3}) that $\sigma=0$. Therefore $G_1p$ is holomorphically equivalent to $\CC^2$. Next, applying the Jacobi identity to the triple $\{a_1,a_3,b_3\}$ and using (\ref{identities1}), (\ref{anbnk2ne0}), (\ref{reln3}), we obtain $\eta\mu=2(\alpha^2+\beta^2)>0$. Therefore, $G_2p$ is equivalent to $\CC\PP^1$. It now follows that ${\frak g}$ is isomorphic to the Lie algebra of the group ${\frak G}\left(\CC\PP^3\setminus\CC\PP^1\right)$ defined in (\ref{gcp3cp1}).

Every connected Lie group that has such a Lie algebra and contains a subgroup isomorphic to $U_2$ is isomorphic to the group ${\frak G}\left(\CC\PP^3\setminus\CC\PP^1\right)$. Therefore, $G$ is isomorphic to ${\frak G}\left(\CC\PP^3\setminus\CC\PP^1\right)$. Under this isomorphism $H$ is mapped into a subgroup $K\subset{\frak G}\left(\CC\PP^3\setminus\CC\PP^1\right)$ for which $K^0$ is isomorphic to $U_2$. It follows from Lemma \ref{fullisotropy} that every connected component of $H$ has the form $gH^0$, where $g$ lies in the centralizer of $H^0$ in $G$. It is straightforward to see, however, that the centralizer of any subgroup of ${\frak G}\left(\CC\PP^3\setminus\CC\PP^1\right)$ isomorphic to $U_2$ coincides with the center of this subgroup, and therefore $K$ is connected. The isomorphism between $G$ and ${\frak G}\left(\CC\PP^3\setminus\CC\PP^1\right)$ induces a diffeomorphism between $G/H$ and ${\frak G}\left(\CC\PP^3\setminus\CC\PP^1\right)/K$.

Observe that there are two kinds of subgroups of ${\frak G}\left(\CC\PP^3\setminus\CC\PP^1\right)$ isomorphic to $U_2$. Subgroups of the first kind are conjugate to the subgroup $K_1$ for which $A=0$ and $V$ is diagonal; subgroups of the second kind are conjugate to the subgroup $K_2$ for which $A=0$ and $U$ is diagonal (see (\ref{gcp3cp1})). There exists an automorphism of ${\frak G}\left(\CC\PP^3\setminus\CC\PP^1\right)$ that maps $K_2$ into $K_1$. Indeed, every map of the form (\ref{gcp3cp1}) is determined by the matrix
$$
C_{U,V,A}:=\left(
\begin{array}{ll}
U & A\\
0 & V
\end{array}
\right).
$$
The automorphism of ${\frak G}\left(\CC\PP^3\setminus\CC\PP^1\right)$ that maps $K_2$ into $K_1$ corresponds to the following transformation of such matrices
$$
C_{U,V,A}\mapsto \left(C_{U,V,A}^{T^a}\right)^{-1},
$$
where $T^a$ denotes the operation of taking the transpose of a matrix with respect to its anti-diagonal. 

Thus, we can assume that an isomorphism between $G$ and ${\frak G}\left(\CC\PP^3\setminus\CC\PP^1\right)$ is chosen so that $H$ is mapped into a subgroup of the first kind. It is then clear that the manifold ${\frak G}\left(\CC\PP^3\setminus\CC\PP^1\right)/K$ is ${\frak G}\left(\CC\PP^3\setminus\CC\PP^1\right)$-equivariantly diffeomorphic to $\CC\PP^3\setminus\{w=0\}$ (see Example (xi)). The linearization of the action of $K$ on $\CC\PP^3\setminus\{w=0\}$ at its fixed point $p_K$ has exactly two non-trivial proper real invariant subspaces $L_1,L_2\subset T_{p_K}\left(\CC\PP^3\setminus\{w=0\}\right)$. Their dimensions are 2 and 4, and they are in fact complex subspaces with respect to the complex structure induced on $\CC\PP^3\setminus\{w=0\}$ by that of $M$. To determine all ${\frak G}\left(\CC\PP^3\setminus\CC\PP^1\right)$-invariant complex structures on $\CC\PP^3\setminus\{w=0\}$ for which $L_1$ and $L_2$ are complex, we again use conditions (\ref{invarcomplstructure}). In this case they lead either to the standard complex structure on $\CC\PP^3\setminus\{w=0\}$ or to its conjugate. The manifold corresponding to the latter structure clearly is holomorphically equivalent to $\CC\PP^3\setminus\{w=0\}$ by means of a map that preserves the group ${\frak G}\left(\CC\PP^3\setminus\CC\PP^1\right)$. Thus, we have obtained Example (xi).

{\bf Case 1.1.2.b.} {\it Suppose that $n=2$.} In this case the group $H$ is 1-dimensional and we apply transformations from $\hbox{Ad}(H^0)$ to the commutators in ${\cal C}_1$. If $k_1\ne\pm 2k_2$, $k_2\ne\pm 2k_1$, we obtain that every element in ${\cal C}_1$ is zero, which leads, as before, to Example (vi) for $n=2$ and $k_2\ne\pm k_1$, $k_1\ne\pm 2k_2$, $k_2\ne\pm 2k_1$.

Assume that $k_1=\pm 2k_2$, that is, $k_1=2$, $k_2=\pm 1$. In this case we obtain 
\begin{equation}
\begin{array}{llr}
[a_1,a_2]&=&\alpha a_2+\beta b_2,\\

[b_1,a_2]&=&\mp\beta a_2\pm\alpha b_2,\\ 

[a_1,b_2]&=&\beta a_2-\alpha b_2,\\

[b_1,b_2]&=&\pm\alpha a_2\pm\beta b_2.
\end{array}\label{identin2}
\end{equation}
If $k_1=2$, $k_2=-1$, then, applying condition (4) of (\ref{invarcomplstructure}) to the $G$-invariant complex structure on $G/H$ and the vectors $x=a_1$, $y=a_2$, and using the fact that the pull-back to ${\frak h}^{\perp}$ of the operator of complex structure on $T_H(G/H)$ is given by (\ref{j0}) in the basis $\{a_k,b_k\}_{k=1,2}$, we obtain that $\alpha=0$, $\beta=0$. Thus, in this case all commutators in ${\cal C}_1$ are zero, and we obtain Example (vi) for $n=2$ and $k_1=2$, $k_2=\pm 1$.

Assume now that $k_1=2$, $k_2=1$. If $\alpha=0$, $\beta=0$, then we obtain Example (vi) for $n=2$ and $k_1=2$, $k_2=\pm 1$. Let $\alpha^2+\beta^2>0$. Applying the Jacobi identity to the triple $\{a_1,b_1,a_2\}$ and using (\ref{arbr}), (\ref{tran}), (\ref{identin2}), we see that $\mu\sigma=-4(\alpha^2+\beta^2)<0$, and hence $G_1p$ is holomorphically equivalent to $\BB^1$. Furthermore, applying the Jacobi identity to the triple $\{a_1,a_2,b_2\}$ and using (\ref{identities1}), (\ref{anbnk2ne0}), (\ref{identin2}), we obtain that $\eta=0$. Therefore, $G_2p$ is equivalent to $\CC$. It now follows that ${\frak g}$ is isomorphic to the Lie algebra of the group ${\frak G}\left(\BB^1\times\CC\right)$ defined in Example (ix).

Every connected Lie group with Lie algebra isomorphic to that of the group ${\frak G}\left(\BB^1\times\CC\right)$ and containing a 1-dimensional compact subgroup is isomorphic either to the group ${\frak G}(\BB^1\times\CC)$ or to its finitely-sheeted cover. Arguing as in Case 1.1.2.a.1, we see that $G$ is isomorphic to ${\frak G}(\BB^1\times\CC)$. This isomorphism induces a diffeomorphism from $G/H$ onto ${\frak G}(\BB^1\times\CC)/K$, where $K$ is a maximal compact subgroup of ${\frak G}(\BB^1\times\CC)$. Next, ${\frak G}(\BB^1\times\CC)/K$ is ${\frak G}(\BB^1\times\CC)$-equivariantly diffeomorphic to $\BB^1\times\CC$. As in Cases 1.1.2.a.1 and 1.1.2.a.2, we are now led to determining all ${\frak G}(\BB^1\times\CC)$-invariant complex structures on $\BB^1\times\CC$ for which the non-trivial proper real invariant subspaces of the linearization of the action of $K$ at its fixed point are in fact complex lines. Conditions (\ref{invarcomplstructure}) yield in this case the standard complex structure on $\BB^1\times\CC$ and its conjugate. It is straightforward to see that the manifold arising from the latter complex structure is equivalent to $\BB^1\times\CC$ by means of a map that preserves the group ${\frak G}(\BB^1\times\CC)$. Thus, we have obtained Example (ix). Observe further that the case $k_2=\pm 2k_1$ differs from the case that we have just considered by the permutation of the indices and therefore leads to Example (ix) as well.

{\bf Case 1.2.} {\it Suppose that ${\cal S}_1$ is not a Lie subalgebra in ${\frak g}$.} For any $v_1,v_2\in{\frak h}_1^{\perp}$ consider the commutator $[v_1,v_2]$ and write it in the general form $[v_1,v_2]=a+b+c$, where $a\in{\frak h}_1^{\perp}$, $b\in{\frak h}_2^{\perp}$, $c\in{\frak h}$. Applying to this identity the element of $\hbox{Ad}(Z)$ that acts as $-\hbox{id}$ on ${\frak h}_1^{\perp}$, we see that $a=0$, hence the projection of $[{\frak h}_1^{\perp},{\frak h}_1^{\perp}]$ to ${\frak h}_1^{\perp}$ is trivial (cf. (\ref{v1v2})). 

{\bf Case 1.2.1.} {\it Suppose that $k_2=0$.} If in the above argument $v_1\in {\frak P}_r$ and $v_2\in {\frak P}_s$ for $r,s=1,\dots,n-1$, $r\ne s$, then we apply to the commutator $[v_1,v_2]$ the transformation from $\hbox{Ad}(H_r)$ that acts as $-\hbox{id}$ on ${\frak P}_r$. Since $\hbox{Ad}(H^0)$ acts trivially on ${\frak h}_2^{\perp}$, we immediately see that $[v_1,v_2]\in{\frak h}$. Arguing now as in Case 1.1, we obtain (\ref{PrPs}).

Since ${\cal S}_1$ is not a Lie subalgebra of ${\frak h}$ and (\ref{hperph}) holds, it follows that for some $1\le r\le n-1$ the commutator $[a_r,b_r]$ has a non-zero projection to ${\frak h}_2^{\perp}$. Write $[a_r,b_r]$ in the general form $[a_r,b_r]=b+c$, where $b,c$ are as above with $b\ne 0$, and let $[a_s,b_s]=b'+c'$ be an analogous decomposition for $1\le s\le n-1$, $s\ne r$. Applying to $[a_r,b_r]$ a transformation from $\hbox{Ad}(H^0)$ that interchanges $a_r$ with $a_s$ and $b_r$ with $b_s$, we see that $b'=b$. Observing, in addition, that $[a_r,b_r]$ is $\hbox{Ad}(H_s)$-invariant for $s=1,\dots,n-1$, we obtain
\begin{equation}
[a_r,b_r]=\omega a_n+\kappa b_n+\sum_{s=1}^{n-1}\sigma_r^st_s,\quad r=1,\dots,n-1,\label{arbrnew}
\end{equation}
for some $\omega,\kappa,\sigma_r^s\in\RR$ with $\omega^2+\kappa^2>0$ (cf. (\ref{arbr})). Applying to identity (\ref{arbrnew}) transformations from $\hbox{Ad}(H^0)$ that interchange $a_r$ with $a_s$ and $b_r$ with $b_s$ for $s=1,\dots,n-1$, $s\ne r$, leaving the other elements of the basis $\{a_k,b_k\}_{k=1,\dots,n}$ fixed, we see that $\sigma_r^r=\sigma_s^s:=\sigma$ for $r,s=1,\dots,n-1$.

Next, we note that relations (\ref{h2perph})--(\ref{id4}) remain true in this case as well. We now define $a_n'$ and $b_n'$ by formulas (\ref{anprimebnprime}) and observe, as in Case 1.1.1, that relations (\ref{identities44}) hold. Hence the subspace ${\cal S}_2'$ of ${\frak g}$ spanned by $a_n'$, $b_n'$ is a Lie subalgebra in ${\frak g}$. Clearly, we have ${\frak F}_2(p)=G_2'p$, and for $g\in G$ the leaf ${\frak F}_2(gp)$ is an orbit of $gG_2'g^{-1}$. The orbit $G_2'p$ is holomorphically equivalent to one of ${\cal P}_{>}$, $\CC$, $\CC^*$, or an elliptic curve $\TT$ by means of a map that transforms $G_2'/K_2'$ into one of $G({\cal P}_{>})$, $G_1(\CC)$, $\hbox{Aut}(\CC^*)^0$, or $\hbox{Aut}(\TT)^0$, respectively, where $K_2'$ is the discrete ineffectivity kernel of the $G_2'$-action on $G_2'p$.

Using $a_n'$, $b_n'$, we can rewrite identity (\ref{arbrnew}) as
$$
[a_r,b_r]=\omega a_n'+\kappa b_n'+\sum_{s=1}^{n-1}{\sigma'}_r^st_s,\quad r=1,\dots,n-1,
$$
for some ${\sigma'}_r^s\in\RR$, with ${\sigma'}_r^r=\sigma-(\beta+\delta)/\mu:=\sigma'$ for all $r$. Now the Jacobi identity applied to all triples $\{a_r,b_r,a_s\}$ with $r,s=1,\dots,n-1$, $r\ne s$, together with (\ref{identities1}), (\ref{PrPs}), (\ref{identities44}) implies that ${\sigma'}_r^s=0$, and therefore we have
\begin{equation}
[a_r,b_r]=\omega a_n'+\kappa b_n'+\sigma' t_r,\quad r=1,\dots,n-1.\label{arbrnewnew}
\end{equation}

Next, the Jacobi identity applied to the triples $\{a_r,b_r,a_n'\}$ and $\{a_r,b_r,b_n'\}$ for $r=1,\dots,n-1$, together with (\ref{h2perph}), (\ref{identities44}), (\ref{arbrnewnew}) yields
\begin{equation}
\begin{array}{ll}
\kappa\rho+2\alpha\omega=0, & \kappa\nu+2\alpha\kappa=0,\\
\omega\rho-2\gamma\omega=0, & \omega\nu-2\gamma\kappa=0,\\
\alpha\sigma'=0, & \gamma\sigma'=0,
\end{array}\label{newcommJac}
\end{equation}

Let ${\frak v}$ be the Lie algebra generated by ${\cal S}_1$. It follows from the above discussion that ${\frak v}$ has dimension $n^2$. Moreover, from (\ref{h2perph}), (\ref{identities44}), (\ref{newcommJac}) we see that ${\frak v}$ is in fact an ideal in ${\frak g}$. Let $V\subset G$ be the normal connected (possibly non-closed) subgroup with Lie algebra ${\frak v}$. Clearly, $H^0$ lies in $V$. The orbit $Vp$ is a (possibly non-closed) real hypersurface in $M$, and the complex subspaces in ${\cal L}_1$ at the points of $Vp$  give complex tangent spaces to $Vp$. Since $LH^0$ acts transitively on directions in ${\cal L}(p)$ and $[a_r,b_r]\ne 0$ for $r=1,\dots,n-1$, the hypersurface $Vp$ is strongly pseudoconvex. 

For a strongly pseudoconvex hypersurface $S$ we denote by $\hbox{Aut}_{CR}(S)$ the Lie group of $CR$-automorphisms of $S$. We have $\hbox{Aut}_{CR}(Vp)\ge n^2$, and therefore $\hbox{dim}\,V_p\ge (n-1)^2$. Hence $p$ is an umbilic point of $Vp$, which due to the homogeneity of $Vp$ yields that $Vp$ is spherical. Since $V$ is normal in $G$, it follows that every orbit of $V$ in $M$ is a spherical hypersurface $CR$-equivalent to $Vp$.

In \cite{I3} (see also \cite{I4}) we listed all $CR$-homogeneous spherical $CR$-manifolds $S$ of dimension $2n-1$ for which $\hbox{Aut}_{CR}(S)\ge n^2$, with $n\ge 2$. Proposition 3.3 of \cite{I4} gives that $Vp$ is $CR$-equivalent to one of:
$$
\begin{array}{ll}
\hbox{(a)}& \hbox{a lens space ${\frak L}_l:=S^{2n-1}/\ZZ_l$ for some $l\in\NN$};\\
\vspace{0cm}\\
\hbox{(b)} & \Sigma:=\left\{(z',z_n)\in\CC^{n-1}\times\CC:\hbox{Re}\,z_n=|z'|^2\right\};\\
\vspace{0cm}\\
\hbox{(c)} & \Upsilon:=\hbox{$\left\{(z',z_n)\in\CC^{n-1}\times\CC: |z_n|=\exp\left(|z'|^2\right)\right\}$};\\
\vspace{0cm}\\
\hbox{(d)} & \Omega:=\left\{(z',z_n)\in\CC^{n-1}\times\CC:|z'|^2+\exp\left(\hbox{Re}\,z_n\right)=1\right\};\\
\vspace{0cm}\\
\hbox{(e)} & E_{\varepsilon}:=\hbox{$\left\{(z',z_n)\in\CC^{n-1}\times\CC: |z'|^2+|z_n|^{\varepsilon}=1,\, z_n\ne 0\right\}$},\\
&\hspace{1.1cm}\hbox{for some $\varepsilon>0$}.
\end{array}
$$
Furthermore, since the group $V$ acts properly and effectively on $Vp$, Proposition 3.4 of \cite{I4} yields that a $CR$-equivalence map can be chosen to transform $V|_{Vp}$ into one of the following groups, respectively:
\vspace{0.5cm}

\noindent (a) ${\cal R}_{{\frak L}_l}:=U_n/\ZZ_l$;
\vspace{0.5cm}

\noindent (b) the group ${\cal R}_{\Sigma}$ which consists of  all maps of the form
$$
\begin{array}{lll}
z' & \mapsto & Uz'+a,\\
z_n & \mapsto & z_n+2\langle Uz',a\rangle+|a|^2+ib,
\end{array}
$$
where $U\in U_{n-1}$, $a\in\CC^{n-1}$, $b\in\RR$, and $\langle\cdot,\cdot\rangle$ denotes the inner product in $\CC^{n-1}$;
\vspace{0.5cm}

\noindent (c) the group ${\cal R}_{\Upsilon}$ which consists of  all maps of the form
$$
\begin{array}{lll}
z' & \mapsto & Uz'+a,\\
z_n & \mapsto &e^{i\psi}\exp\Bigl(2\langle Uz',a\rangle+|a|^2\Bigr)z_n,
\end{array}
$$
where $U,a$ are as above, and $\psi\in\RR$;
\vspace{0.5cm}

\noindent (d) the group ${\cal R}_{\Omega}$ which consists of  all maps of the form
$$
\begin{array}{lll}
z'&\mapsto&\displaystyle\frac{Az'+b}{cz'+d},\\
\vspace{0mm}&&\\
z_n&\mapsto&\displaystyle z_n-2\ln(cz'+d)+ia,
\end{array}
$$
where $A,b,c,d$ are as in (\ref{autballn-1}), and $a\in\RR$;
\vspace{0.5cm}

\noindent (e) the group ${\cal R}_{E_{\varepsilon}}$ which consists of  all maps of the form
$$
\begin{array}{lll}
z'&\mapsto&\displaystyle\frac{Az'+b}{cz'+d},\\
\vspace{0mm}&&\\
z_n&\mapsto&\displaystyle
\frac{e^{i\psi}z_n}{(cz'+d)^{2/\varepsilon}},
\end{array}
$$
where $A,b,c,d$ are as in (\ref{autballn-1}), and $\psi\in\RR$.

Suppose first that $V$ is a closed subgroup in $G$. Then $V$ acts properly on $M$ and therefore every orbit of $V$ is closed in $M$. Under this assumption, arguing as in Section 3.4 of \cite{I4} (see also Theorem 2.7 in \cite{IKru} for case (a) with $l=1$), one can determine $M$. We will treat each of cases (a)--(e) separately.

Suppose that $Vp$ is $CR$-equivalent to ${\frak L}_l$ for some $l\in\NN$. Then $M$ is holomorphically equivalent to one of 
$$
\begin{array}{l}
\hbox{(a1)}\,\, S_{R_1,R_2,l}:=\left\{z\in\CC^n:R_1<|z|<R_2\right\}/\ZZ_l,\,0\le R_1<R_2\le\infty,\\
\vspace{0.1cm}\\
\hbox{(a2)}\,\, M_d/\ZZ_l,\,d\in\CC^*, |d|\ne 1,
\end{array}
$$
where $M_d$ is a Hopf manifold (see Example (vii)). Observe that if in subcase (a1) we have either $R_1>0$ or $R_2<\infty$, then $\hbox{Aut}(S_{R_1,R_2,l})={\cal R}_{{\frak L}_l}$ and hence $S_{R_1,R_2,l}$ does not admit an action of an $(n^2+1)$-dimensional group for such $R_1$, $R_2$. Therefore, in subcase (a1) we in fact have $R_1=0$, $R_2=\infty$, (note that $S_{0,\infty,l}=\CC^{n{}*{}}/\ZZ_l$). Furthermore, in each of subcases (a1), (a2) an equivalence map can be chosen to transform $V$ into ${\cal R}_{{\frak L}_l}$.
 
Suppose next that $Vp$ is $CR$-equivalent to $\Sigma$. Then $M$ is holomorphically equivalent to one of 
$$
\begin{array}{ll}
\hbox{(b1)}\,\,\Sigma_{R_1,R_2}:=&\left\{(z',z_n)\in\CC^{n-1}\times\CC: R_1+|z'|^2<\hbox{Re}\,z_n<R_2+|z'|^2\right\},\\
&-\infty\le R_1<R_2\le\infty,\\
\vspace{0.1cm}\\
\hbox{(b2)}\,\,\CC^{n-1}\times\CC^*.&
\end{array}
$$
Observe that if in subcase (b1) we have $R_1>-\infty$ and $R_2<\infty$, then $\hbox{Aut}(\Sigma_{R_1,R_2})={\cal R}_{\Sigma}$ and hence $\Sigma_{R_1,R_2}$ does not admit an action of an $(n^2+1)$-dimensional group for such $R_1$, $R_2$. Therefore, in subcase (b1) we have either $R_1=-\infty$, $R_2<\infty$ or $R_1>-\infty$, $R_2=\infty$ or $R_1=-\infty$, $R_2=\infty$ (note that $\Sigma_{-\infty,R_2}$ is equivalent to ${\cal P}^n_{<}$ if $R_2<\infty$ (see (\ref{p<})), $\Sigma_{R_1,\infty}$ is equivalent to ${\cal P}^n_{>}$ if $R_1>-\infty$ (see (\ref{p>})), and $\Sigma_{-\infty,\infty}=\CC^n$). Furthermore, an equivalence map can be chosen so that it transforms $V$ in subcase (b1) into ${\cal R}_{\Sigma}$ and in subcase (b2) into the group of maps     
$$
\begin{array}{lll}
z' & \mapsto & Uz'+a,\\
z_n & \mapsto &\exp\Bigl(\Lambda(2\langle Uz',a\rangle+|a|^2+ib)\Bigr)z_n,
\end{array}
$$
where $U\in U_{n-1}$, $a\in\CC^{n-1}$, $b\in\RR$, and $\Lambda$ is a fixed complex number with $\hbox{Im}\,\Lambda\ne 0$.

Suppose that $Vp$ is $CR$-equivalent to $\Upsilon$. Then $M$ is holomorphically equivalent to one of
$$
\begin{array}{ll}
\hbox{(c1)}\,\,\Upsilon_{R_1,R_2}:=&\left\{(z',z_n)\in\CC^{n-1}\times\CC: R_1\exp\left({|z'|^2}\right)<|z_n|<R_2\exp\left({|z'|^2}\right)\right\},\\
&0\le R_1<R_2\le\infty,\\
\vspace{0.1cm}\\
\hbox{(c2)}\,\,\CC^{n-1}\times\TT,&
\end{array}
$$
where $\TT$ is an elliptic curve. Observe that if in subcase (c1) we have $R_1>0$ or $R_2<\infty$, then $\hbox{Aut}(\Upsilon_{R_1,R_2})={\cal R}_{\Upsilon}$ and hence $\Upsilon_{R_1,R_2}$ does not admit an action of an $(n^2+1)$-dimensional group for such $R_1$, $R_2$. Therefore, in subcase (c1) we have $R_1=0$, $R_2=\infty$, (note that $\Upsilon_{0,\infty}=\CC^{n-1}\times\CC^*$). Furthermore, an equivalence map can be chosen so that it transforms $V$ in subcase (c1) into ${\cal R}_{\Upsilon}$ and in subcase (c2) into the group of maps 
$$
\begin{array}{lll}
z' & \mapsto & Uz'+a,\\

[z_n] & \mapsto &\left[e^{i\psi}\exp\Bigl(2\langle Uz',a\rangle+|a|^2\Bigr)z_n\right],
\end{array}
$$
where $U\in U_{n-1}$, $a\in\CC^{n-1}$ and $\psi\in\RR$. Here $\TT$ is obtained from $\CC^*$ by taking the quotient with respect to the equivalence relation $z_n\sim dz_n$, for some $d\in\CC^*$, $|d|\ne 1$, and $[z_n]\in\TT$ is the equivalence class of a point $z_n\in\CC^*$.
     
Suppose next that $Vp$ is $CR$-equivalent to $\Omega$. Then $M$ is holomorphically equivalent to one of
$$
\begin{array}{lll}
\hbox{(d1)}\,\,&\Omega_{R_1,R_2}:=\Bigl\{(z',z_n)\in\CC^{n-1}\times\CC: |z'|<1,\\
&\hspace{2cm}R_1(1-|z'|^2)<\exp\left(\hbox{Re}\,z_n\right)<R_2(1-|z'|^2)\Bigr\},\\
&\hspace{2cm}0\le R_1<R_2\le\infty,&\\
\vspace{0.1cm}\\
\hbox{(d2)}\,\,&\BB^{n-1}\times\CC^*.&
\end{array}
$$
Observe that if in subcase (d1) we have $R_1>0$ or $R_2<\infty$, then $\hbox{Aut}(\Omega_{R_1,R_2})={\cal R}_{\Omega}$ and hence $\Omega_{R_1,R_2}$ does not admit an action of an $(n^2+1)$-dimensional group for such $R_1$, $R_2$. Therefore, in subcase (d1) we have $R_1=0$, $R_2=\infty$, (note that $\Omega_{0,\infty}=\BB^{n-1}\times\CC$). Furthermore, an equivalence map can be chosen so that it transforms $V$ subcase (d1) into ${\cal R}_{\Omega}$ and in subcase (d2) into the group of maps 
$$
\begin{array}{lll}
z'&\mapsto&\displaystyle\frac{Az'+b}{cz'+d},\\
\vspace{0mm}&&\\
z_n&\mapsto&\displaystyle\frac{e^{i\psi\Lambda}z_n}{(cz'+d)^{2\Lambda}},
\end{array}
$$
where $A,b,c,d$ are as in (\ref{autballn-1}), $\psi\in\RR$, and $\Lambda$ is a fixed complex number with $\hbox{Im}\,\Lambda\ne 0$.

Suppose that $Vp$ is $CR$-equivalent to $E_{\varepsilon}$ for some $\varepsilon>0$. Then $M$ is holomorphically equivalent to one of
$$
\begin{array}{lll}
\hbox{(e1)}\,\,&E_{R_1,R_2,\varepsilon}:=&\Bigl\{(z',z_n)\in\CC^{n-1}\times\CC: |z'|<1,\\
&&R_1(1-|z'|^2)^{1/\varepsilon}<|z_n|<R_2(1-|z'|^2)^{1/\varepsilon}\Bigr\},\\
&&0\le R_1<R_2\le\infty,\\
\vspace{0.1cm}\\
\hbox{(e2)}\,\,&\BB^{n-1}\times\TT,& 
\end{array}
$$   
where $\TT$ is an elliptic curve. Observe that if in subcase (e1) we have $R_1>0$ or $R_2<\infty$, then $\hbox{Aut}(E_{R_1,R_2,\varepsilon})={\cal R}_{E_{\varepsilon}}$ and hence $E_{R_1,R_2,\varepsilon}$ does not admit an action of an $(n^2+1)$-dimensional group for such $R_1$, $R_2$. Therefore, in subcase (e1) we have $R_1=0$, $R_2=\infty$, (note that $E_{0,\infty,\varepsilon}=\BB^{n-1}\times\CC^*$). Furthermore, an equivalence map can be chosen so that it transforms $V$ subcase (e1) into ${\cal R}_{E_{\varepsilon}}$ and in subcase (e2) into the group of maps
$$
\begin{array}{lll}
z'&\mapsto&\displaystyle\frac{Az'+b}{cz'+d},\\
\vspace{0mm}&&\\

[z_n]&\mapsto&\displaystyle
\left[\frac{e^{i\psi}z_n}{(cz'+d)^{2/\varepsilon}}\right],
\end{array}
$$
where $A,b,c,d$ are as in (\ref{autballn-1}), and $\psi\in\RR$. Here $\TT$ is obtained from $\CC^*$ by taking the quotient with respect to the equivalence relation $z_n\sim dz_n$, for some $d\in\CC^*$, $|d|\ne 1$, and $[z_n]\in\TT$ is the equivalence class of a point $z_n\in\CC^*$.

We will now determine the groups into which $G$ is transformed under the above equivalences. Observe that for every point $q\in M$ the fixed point set ${\cal F}_{V_q^0}$ of the identity component of the isotropy subgroup $V_q$ of $q$ with respect to the action of $V$ contains the leaf ${\frak F}_2(q)$. Examining the fixed point sets of the isotropy subgroups with respect to the action of the group into which $V$ is transformed in each of the cases listed above, we see that ${\cal F}_{V_q^0}$ is always a connected complex curve in $M$ and therefore ${\frak F}_2(q)={\cal F}_{V_q^0}$ for $q\in M$. This implies that the leaves of ${\frak F}_2$ are equivalent to $\CC^*$ in subcase (a1); to an elliptic curve $\TT$ in subcase (a2); to ${\cal P}_{>}$ if either $R_1=-\infty$, $R_2<\infty$ or $R_1>-\infty$, $R_2=\infty$, and to $\CC$ if $R_1=-\infty$, $R_2=\infty$ in subcase (b1); to $\CC^*$ in subcase (b2); to $\CC^*$ in subcase (c1); to an elliptic curve $\TT$ in subcase (c2); to $\CC$ in subcase (d1); to $\CC^*$ in subcase (d2), to $\CC^*$ in subcase (e1); and to an elliptic curve $\TT$ in subcase (e2). 

If the leaves of ${\frak F}_2$ are equivalent to one of $\CC$, $\CC^*$, $\TT$, then the algebra ${\cal S}_2'$ is commutative, and therefore $\rho=\nu=0$. Due to (\ref{newcommJac}) this implies that $\alpha=\gamma=0$. It now follows from (\ref{h2perph}), (\ref{identities44}) that the Lie subalgebra ${\cal S}_2'$ lies in the center of ${\frak g}$. Fix $v\in{\cal S}_2'$, $v\not\in{\frak v}$, and consider the one-parameter subgroup $\{{\cal G}(t)\}_{t\in\RR}$ of holomorphic transformation of $M$ arising from $v$. Since $\{{\cal G}(t)\}_{t\in\RR}$ lies in the center of $G$, in all the subcases, except for (b1) with either $R_1=-\infty$, $R_2<\infty$ or $R_1>-\infty$, $R_2=\infty$, $\{{\cal G}(t)\}_{t\in\RR}$ is transformed into a one-parameter subgroup of holomorphic automorphisms of the corresponding manifold whose every element commutes with every element of the group into which $V$ is transformed. Looking at the explicit forms of the groups that appear in all the subcases one can determine the form of the subgroup into which $\{{\cal G}(t)\}_{t\in\RR}$ is transformed. This yields that the group $G$ is transformed into the group defined in (\ref{gviiia}) in subcase (a1); into the group defined in (\ref{gviiib}) in subcase (a2); into the group defined in (\ref{g4cns}) for $T=2$ in subcase (b1) if $R_1=-\infty$, $R_2=\infty$; 
into the group defined in (\ref{g4cnsstar}) for $T=2\Lambda$ in subcase (b2), into the group defined in (\ref{g4cnsstar}) for $T=2$ in subcase (c1); into the group defined in (\ref{g4cnst}) for $T=2$ in subcase (c2); into the group defined in (\ref{groupbcs}) for $T=-2$ in subcase (d1); into the group defined in (\ref{groupbcstar}) for $T=-2\Lambda$ in subcase (d2); into the group defined in (\ref{groupbcstar}) for $T=-2/\varepsilon$ in subcase (e1); and into the group defined in (\ref{groupbtorus}) for $T=-2/\varepsilon$ in subcase (e2). Thus, these subcases lead to the following examples: (a1) to (viia), (a2) to (viib), (b1) for $R_1=-\infty$, $R_2=\infty$ to (iva), (b2) to (ivb) with $e^{i\tau}\not\in\RR$, (c1) to (ivb) with $e^{i\tau}\in\RR$, (c2) to (ivc) with $e^{i\tau}\in\RR$, (d1) to (iia), (d2) to (iib) for $T\in\CC^*\setminus\RR^*$, (e1) to (iib) for $T\in\RR^*$, and (e2) to (iic) for $T\in\RR^*$.

It remains to consider subcase (b1) with either $R_1=-\infty$, $R_2<\infty$ or $R_1>-\infty$, $R_2=\infty$. Applying a translation in $z_n$ we can assume that $R_2=0$ in the first case and $R_1=0$ in the second case. Every element of $\hbox{Aut}\left(\Sigma_{-\infty,0}\right)$ extends to an element of $\hbox{Aut}(\CC^n)$ that preserves the domain $\Sigma_{0,\infty}$ (which is equivalent to $\BB^n$). It is easy to see that every such map belongs to the group $G_{\cal P}$ defined in Example (viii). Hence if $R_1=-\infty$, $R_2=0$, the group $G$ is transformed into $G_{\cal P}$. Suppose now that $R_1=0$, $R_2=\infty$ (hence $M$ is holomorphically equivalent to $\BB^n$). Since the subgroup $V$ is normal in $G$, the foliation of $M$ by $V$-orbits is $G$-invariant. It is straightforward to show that every element of $\hbox{Aut}\left(\Sigma_{0,\infty}\right)$ that preserves the foliation of $\Sigma_{0,\infty}$ by ${\cal R}_{\Sigma}$-orbits lies in the group $G_{\cal P}$. Hence if $R_1=0$, $R_2=\infty$, the group $G$ is transformed into $G_{\cal P}$ as well. Thus, if exactly one of $R_1$, $R_2$ is infinite, subcase (b1) leads to Examples (viiia), (viiib).

Assume now that $V$ is not closed in $G$. This implies that the orbit $Vp$ accumulates to itself (in the sense explained in the proof of Theorem \ref{classtype2} above). Hence $Vp$ is non-compact and thus is $CR$-equivalent to one of $\Sigma$, $\Upsilon$, $\Omega$, $E_{\varepsilon}$, for some $\varepsilon>0$. Clearly, the intersection ${\mathscr T}:=Vp\cap G_2'p$ is $H^0$-invariant. Furthermore, since $H^0$ acts trivially on $G_2'p$, every point in ${\mathscr T}$ is a fixed point of $H^0$. A direct examination shows that for each of the groups ${\cal R}_{\Sigma}$, ${\cal R}_{\Upsilon}$, ${\cal R}_{\Omega}$, ${\cal R}_{E_{\varepsilon}}$ the fixed point set of the isotropy subgroup of every point in $\Sigma$, $\Upsilon$, $\Omega$, $E_{\varepsilon}$, respectively, is a connected smooth curve. It is straightforward to observe that every such a curve is an orbit of a one-parameter subgroup of one of ${\cal R}_{\Sigma}$, ${\cal R}_{\Upsilon}$, ${\cal R}_{\Omega}$, ${\cal R}_{E_{\varepsilon}}$, respectively, and that no proper subset of the curve is an orbit of a one-parameter subgroup of any of these groups.

Let $v$ be a non-zero vector tangent to each of $Vp$ and $G_2'p$ at $p$ and $\{{\cal G}(t)\}_{t\in\RR}$ a one-parameter subgroup of $G_2'$ such that $v=d/dt\,({\cal G}(t)p)|_{t=0}$. Since the foliation of $M$ by $V$-orbits is $G$-invariant and the foliation ${\frak F}_2$ is $G$-invariant as well, the vector field arising from the action of $\{{\cal G}(t)\}_{t\in\RR}$ on $M$ is tangent to each of $Vp$ and $G_2'p$ at all points. In particular, the orbit $\{{\cal G}(t)p\}_{t\in\RR}$ lies in ${\mathscr T}$. Since $\{{\cal G}(t)\}_{t\in\RR}$ preserves $Vp$, it is a one-parameter subgroup of $V$, and the discussion in the preceding paragraph yields that ${\mathscr T}$ coincides with the orbit $\{{\cal G}(t)p\}_{t\in\RR}$. Since $Vp$ accumulates to itself, so does the orbit $\{{\cal G}(t)p\}_{t\in\RR}$.

Recall that $G_2'p$ is holomorphically equivalent to one of ${\cal P}_{>}$ , $\CC$, $\CC^*$, or an elliptic curve $\TT$ by means of a map that transforms $G_2'/K_2'$ into one of $G({\cal P}_{>})$, $G_1(\CC)$, $\hbox{Aut}(\CC^*)^0$, or $\hbox{Aut}(\TT)^0$, respectively, where $K_2'$ is the ineffectivity kernel of the $G_2'$-action on $G_2'p$. It then follows that the only case in which $G_2'$ has a one-parameter
subgroup for which there exists an orbit in $G_2'p$ that accumulates to itself is when $G_2'p$ is equivalent to $\TT$. In this case ${\cal S}_2'$ is commutative, i.e. $\rho=\nu=0$, and hence $\alpha=\gamma=0$ (see (\ref{newcommJac})). Now (\ref{h2perph}), (\ref{identities44}) imply that the Lie subalgebra ${\cal S}_2'$ lies in the center of ${\frak g}$.

Since $V$ is not isomorphic to $R_{{\frak L}_l}$ for any $l\in\NN$, its Lie algebra ${\frak v}$ is isomorphic to either the Lie algebra of ${\cal R}_{\Sigma}$ or to the Lie algebra of ${\cal R}_{\Omega}$ (observe that ${\cal R}_{\Sigma}$ covers ${\cal R}_{\Upsilon}$ and ${\cal R}_{\Omega}$ covers ${\cal R}_{E_{\varepsilon}}$ for all $\varepsilon>0$). Therefore, ${\frak g}$ is isomorphic to either the Lie algebra of ${\frak G}(\CC^n)$ (see (\ref{g4cns}) for $T=1$) or to the Lie algebra of $G\left(\BB^{n-1}\times\CC\right)$ (see (\ref{groupbcs}) for $T=1$), respectively. 

Let $\tilde M$ be the universal cover of $M$ and $\Pi:\tilde M\ra M$ a covering map. We lift the complex structure on $M$ to $\tilde M$ and consider the group $\tilde G\subset\hbox{Aut}(\tilde M)$ that consists of all lifts of all elements of $G$ to $\tilde M$. If ${\mathscr G}$ is a $G$-invariant Riemannian metric on $M$, then the lift $\tilde{\mathscr G}$ of ${\mathscr G}$ to $\tilde M$ is a $\tilde G$-invariant Riemannian metric on $\tilde M$. To show that $\tilde G$ acts properly on $\tilde M$, we will prove that $\tilde G$ is a closed subgroup of the group $\hbox{Isom}(\tilde M,\tilde{\mathscr G})$ of all isometries of $\tilde M$ with respect to $\tilde{\mathscr G}$ endowed with the compact-open topology. Indeed, let $\{\tilde f_n\}$ be a sequence in $\tilde G$ that converges in $\hbox{Isom}(\tilde M,\tilde{\mathscr G})$ to an isometry $\tilde F$. Consider a sequence $\{f_n\}$ in $G$ such that $\tilde f_n$ is a lift of $f_n$ to $\tilde M$ for every $n$. It then follows that $\{f_n\}$ converges to a smooth map $F:M\ra M$ uniformly with all derivatives on compact subsets of $M$. Considering the sequences $\{\tilde f_n^{-1}\}$ and $\{f_n^{-1}\}$ of the inverse maps, we see that $F$ is in fact a diffeomorphism of $M$. Clearly, $\tilde F$ is a lift of $F$ to $\tilde M$. Since $G$ is closed in $\hbox{Diff}(M)$, we have $F\in G$, hence $\tilde F\in\tilde G$. Thus, $\tilde G$ acts effectively and properly on $\tilde M$ by holomorphic transformations and so does $G':=\tilde G^0$. Clearly, the group $G'$ covers $G$, and let ${\mathscr P}: G'\ra G$ be the natural covering map. Fix $\tilde p\in\Pi^{-1}(p)$ and let $H':=G'_{\tilde p}$. We identify $\tilde M$ as a smooth manifold with $G'/H'$. The subgroup $H'$ is connected, compact, and we have ${\mathscr P}(H')=H^0$.

Suppose first that ${\frak g}$ is isomorphic to the Lie algebra of ${\frak G}(\CC^n)$. One can verify that every connected group ${\frak G}$ with Lie algebra isomorphic to that of ${\frak G}(\CC^n)$ and having a compact subgroup ${\frak K}$ such that ${\frak G}/{\frak K}$ is simply-connected and ${\frak G}$ acts on ${\frak G}/{\frak K}$ effectively, is isomorphic to ${\frak G}(\CC^n)$ in which case ${\frak K}$ is a maximal compact subgroup of ${\frak G}$. For every maximal compact subgroup $K$ of ${\frak G}(\CC^n)$ the quotient ${\frak G}(\CC^n)/K$ is ${\frak G}(\CC^n)$-equivariantly diffeomorphic to $\CC^n$. The linearization of the action of $K$ on $\CC^n$ at its fixed point $p_K$ has exactly two non-trivial proper even-dimensional real invariant subspaces $L_1,L_2\subset T_{p_K}(\CC^n)$. Their dimensions are $2(n-1)$ and 2. As earlier, using conditions (\ref{invarcomplstructure}) we will find all ${\frak G}(\CC^n)$-invariant complex structures on $\CC^n$ for which $L_1$ and $L_2$ are complex subspaces. Every such structure turns out to be holomorphically equivalent, by means of a map preserving the group ${\frak G}(\CC^n)$, to a complex structure for which the corresponding operator at the origin has the form
\begin{equation}
\begin{array}{lll}
x_r&\mapsto& -y_r, \\
y_r &\mapsto& x_r, \\
x_n &\mapsto& w x_n+u y_n,\\
y_n &\mapsto& v x_n-w y_n,
\end{array}\label{endomorph}
\end{equation}
where $r=1,\dots,n-1$, $v,u,w\in\RR$, $w^2+uv=-1$, $x_k=\hbox{Re}\,z_k$, $y_k=\hbox{Im}\,z_k$ for $k=1,\dots,n$. One can now verify that the diffeomorphism
$$
\begin{array}{lll}
x_r&\mapsto&x_r,\\
\vspace{-0.3cm}\\
y_r&\mapsto&y_r,\\
\vspace{-0.3cm}\\
x_n&\mapsto&\displaystyle\frac{u+1}{2u}\sum_{r=1}^{n-1}(x_r^2+y_r^2)-\frac{1}{u}x_n,\\
\vspace{-0.3cm}\\
y_n&\mapsto&\displaystyle-\frac{w}{2u}\sum_{r=1}^{n-1}(x_r^2+y_r^2)+\frac{w}{u}x_n+y_n,
\end{array}
$$
where $r=1,\dots,n-1$, preserves the group ${\frak G}(\CC^n)$ and establishes holomorphic equivalence between $\CC^n$ equipped with the above complex structure and $\CC^n$ equipped with the standard complex structure.   

Thus, $\tilde M$ is holomorphically equivalent to $\CC^n$ by means of a map that transforms $G'$ into ${\frak G}(\CC^n)$. Let $\Gamma$ be the group of deck transformations of the cover $\Pi:\tilde M\ra M$. Since $\Gamma$ lies in the centralizer of $G'$ in $\hbox{Aut}(\tilde M)$, it follows that, upon the identification of $\tilde M$ with $\CC^n$ and $G'$ with ${\frak G}(\CC^n)$, every element of $\Gamma$ has the form 
\begin{equation}
\begin{array}{lll}
z_r & \mapsto & z_r,\\
z_n & \mapsto & z_n+c,
\end{array}\label{Gamma}
\end{equation}
where $r=1,\dots,n-1$ and $c\in\CC$. In particular, $\Gamma\subset {\frak G}(\CC^n)$, and therefore, upon the identification of $M$ with $\CC^n/\Gamma$, the group $G$ is identified with ${\frak G}(\CC^n)/\Gamma$. It is clear that $\CC^n/\Gamma$ admits a ${\frak G}(\CC^n)/\Gamma$-invariant foliation by elliptic curves only if $\Gamma$ has two generators. Thus, $M$ is holomorphically equivalent to $\CC^{n-1}\times\TT$, where $\TT$ is an elliptic curve, by means of a map that transforms $G$ into a group of maps of the form (\ref{g4cnst}) for some $T\in\CC^*$. Since $V$ is not closed in $G$, the number $T$ cannot be real for any factorization $\CC\ra\CC^*\ra\TT$. Thus, we have obtained Example (ivc) for $e^{i\tau}\not\in\RR$.

Suppose next that ${\frak g}$ is isomorphic to the Lie algebra of $G\left(\BB^{n-1}\times\CC\right)$. We will proceed as in the case considered above. It is straightforward to verify that every connected group ${\frak G}$ with Lie algebra isomorphic to that of $G\left(\BB^{n-1}\times\CC\right)$ and having a compact subgroup ${\frak K}$ such that ${\frak G}/{\frak K}$ is simply-connected and ${\frak G}$ acts on ${\frak G}/{\frak K}$ effectively, is isomorphic to $G\left(\BB^{n-1}\times\CC\right)$ in which case ${\frak K}$ is a maximal compact subgroup of ${\frak G}$. For every maximal compact subgroup $K$ of $G\left(\BB^{n-1}\times\CC\right)$ the quotient $G\left(\BB^{n-1}\times\CC\right)/K$ is $G\left(\BB^{n-1}\times\CC\right)$-equivariantly diffeomorphic to $\BB^{n-1}\times\CC$. The linearization of the action of $K$ on $\BB^{n-1}\times\CC$ at its fixed point $p_K$ has exactly two non-trivial proper even-dimensional real invariant subspaces $L_1,L_2\subset T_{p_K}(\BB^{n-1}\times\CC)$. Their dimensions are $2(n-1)$ and 2. Using conditions (\ref{invarcomplstructure}) we can determine all $G\left(\BB^{n-1}\times\CC\right)$-invariant complex structures on $\BB^{n-1}\times\CC$ for which $L_1$ and $L_2$ are complex subspaces. Every such structure turns out to be holomorphically equivalent, by means of a map preserving the group $G\left(\BB^{n-1}\times\CC\right)$, to a complex structure for which the corresponding operator at the origin has the form (\ref{endomorph}). One can now verify that the diffeomorphism
$$
\begin{array}{lll}
x_r&\mapsto&x_r,\\
\vspace{-0.3cm}\\
y_r&\mapsto&y_r,\\
\vspace{-0.3cm}\\
x_n&\mapsto&\displaystyle-\frac{u+1}{2u}\ln\left(1-\left(\sum_{r=1}^{n-1}(x_r^2+y_r^2)\right)\right)-\frac{1}{u}x_n,\\
\vspace{-0.3cm}\\
y_n&\mapsto&\displaystyle\frac{w}{2u}\ln\left(1-\left(\sum_{r=1}^{n-1}(x_r^2+y_r^2)\right)\right)+\frac{w}{u}x_n+y_n,
\end{array}
$$
where $r=1,\dots,n-1$, preserves the group $G\left(\BB^{n-1}\times\CC\right)$ and establishes holomorphic equivalence between $\BB^{n-1}\times\CC$ equipped with this complex structure and $\BB^{n-1}\times\CC$ equipped with the standard complex structure.

Thus, $\tilde M$ is holomorphically equivalent to $\BB^{n-1}\times\CC$ by means of a map that transforms $G'$ into $G\left(\BB^{n-1}\times\CC\right)$. Let, as before, $\Gamma$ be the group of deck transformations of the cover $\Pi:\tilde M\ra M$. Since $\Gamma$ lies in the centralizer of $G'$ in $\hbox{Aut}(\tilde M)$, it follows that, upon the identification of $\tilde M$ with $\BB^{n-1}\times\CC$ and $G'$ with $G\left(\BB^{n-1}\times\CC\right)$, every element of $\Gamma$ has the form (\ref{Gamma}). In particular, $\Gamma\subset G\left(\BB^{n-1}\times\CC\right)$, and therefore, upon the identification of $M$ with $\left(\BB^{n-1}\times\CC\right)/\Gamma$, the group $G$ is identified with $G\left(\BB^{n-1}\times\CC\right)/\Gamma$. Clearly, $\left(\BB^{n-1}\times\CC\right)/\Gamma$ admits a $G\left(\BB^{n-1}\times\CC\right)/\Gamma$-invariant foliation by elliptic curves only if $\Gamma$ has two generators. Thus, $M$ is holomorphically equivalent to $\BB^{n-1}\times\TT$, where $\TT$ is an elliptic curve, by means of a map that transforms $G$ into a group of maps of the form (\ref{groupbtorus}) for some $T\in\CC^*$. Since $V$ is not closed in $G$, the number $T$ cannot be real for any factorization $\CC\ra\CC^*\ra\TT$. Thus, we have obtained Example (iic) for $T\in\CC^*\setminus\RR^*$.

{\bf Case 1.2.2.} {\it Suppose that $k_2\ne 0$.} Since the commutators $[a_r,b_r]$,  are $\hbox{Ad}(H_s)$-invariant for $r,s=1,\dots,n-1$, it follows, as in Case 1.1, that (\ref{formarbr}) holds. Together with (\ref{hperph}) this implies that for $n=2$ the subspace ${\cal S}_1$ is a Lie subalgebra of ${\frak g}$, which contradicts our assumptions. Hence $n\ge 3$. We will now obtain a similar contradiction for $n\ge 4$ thus proving that in fact $n=3$.    

We showed at the beginning of Case 1.2 that $[{\frak h}_1^{\perp},{\frak h}_1^{\perp}]$ has a zero projection to ${\frak h}_1^{\perp}$. Let $v_1\in{\frak P}_r$, $v_2\in{\frak P}_s$, for $1\le r,s\le n-1$, $r\ne s$, and write the commutator $[v_1,v_2]$ in the general form $[v_1,v_2]=b+c$, where $b\in{\frak h}_2^{\perp}$, $c\in{\frak h}$. If $n\ge 4$, we apply to this commutator a transformation from $\hbox{Ad}(H^0)$ that acts as $-\hbox{id}$ on ${\frak h}_2^{\perp}$ and trivially on each of ${\frak P}_r$, ${\frak P}_s$. This immediately gives that $b=0$, that is
\begin{equation}
[{\frak P}_r,{\frak P}_s]\subset{\frak h},\quad r,s=1,\dots,n-1,\,r\ne s,\,n\ge 4.\label{PsPrnew}
\end{equation}
Relations (\ref{hperph}), (\ref{formarbr}), (\ref{PsPrnew}) imply that ${\cal S}_1$ is a Lie subalgebra of ${\frak g}$ for $n\ge 4$. Thus, we have shown that in Case 1.2.2 we have $n=3$.

Let ${\cal C}$ be the following collection of commutators
$$
{\cal C}:=\left\{[a_1,a_2],\,[a_1,b_2],\,[b_1,a_2],\,[b_1,b_2]\right\}.
$$
We now apply transformations from $\hbox{Ad}(H_1)$ to the elements of ${\cal C}$. If $k_1\ne\pm k_2$,  we immediately obtain that ${\cal C}\subset{\frak h}$. Taken together with (\ref{hperph}), (\ref{formarbr}), this implies that ${\cal S}_1$ is a subalgebra of ${\frak g}$, once again contradicting our assumptions. Hence we have $k_1=\pm k_2$, that is, $k_1=1$, $k_1=\pm 1$. The application of $\hbox{Ad}(H_r)$ for $r=1,2$ to the commutators in ${\cal C}$ then yields
\begin{equation}
\begin{array}{llrl}
[a_1,a_2]&\equiv&\gamma a_3+\delta b_3\,&\hbox{(mod ${\frak h}$)},\\

[a_1,b_2]&\equiv&\mp\delta a_3\pm\gamma b_3\,&\hbox{(mod ${\frak h}$)},\\

[b_1,a_2]&\equiv&\mp\delta a_3\pm\gamma b_3\,&\hbox{(mod ${\frak h}$)},\\

[b_1,b_2]&\equiv&-\gamma a_3-\delta b_3\,&\hbox{(mod ${\frak h}$)},
\end{array}\label{a1a2b1b2}
\end{equation}
for some $\gamma,\delta\in\RR$, $\gamma^2+\delta^2>0$.

If $k_1=1$, $k_2=-1$, then, applying condition (4) of (\ref{invarcomplstructure}) to the $G$-invariant complex structure on $G/H$ and the vectors $x=a_1$, $y=a_2$, and using the fact that the pull-back to ${\frak h}^{\perp}$ of the operator of complex structure on $T_H(G/H)$ is given by (\ref{j0}) in the basis $\{a_k,b_k\}_{k=1,2,3}$, we obtain that $\gamma=0$, $\delta=0$. This contradiction shows that in fact $k_1=1$, $k_2=1$.

Now the application of transformations from $\hbox{Ad}(Z)$ to the commutators in ${\cal C}$ implies
\begin{equation}
\begin{array}{llrl}
[a_1,a_2]&\equiv&[b_1,b_2]\,&\hbox{(mod ${\frak h}_2^{\perp}$)},\\
\vspace{-0.3cm}\\

[b_1,a_2]&\equiv&-[a_1,b_2]\,&\hbox{(mod ${\frak h}_2^{\perp}$)},\\
\end{array}\label{a1a2h}
\end{equation}
and we denote the ${\frak h}$-components of $[a_1,a_2]$ and $[b_1,a_2]$ by $c$ and $d$, respectively. Applying once again transformations from $\hbox{Ad}(H_1)$ to the elements of ${\cal C}$ and comparing the ${\frak h}$-components, we see that $c$ and $d$ are represented by matrices of the form (\ref{frakA}) for which the corresponding matrices ${\frak A}$ are
\begin{equation}
\left(
\begin{array}{ll}
0 & \nu\\
-\overline{\nu} & 0
\end{array}
\right),\quad
\left(
\begin{array}{ll}
0 & i\nu\\
i\overline{\nu} & 0
\end{array}
\right),\label{frakAcd}
\end{equation}
respectively, for some $\nu\in\CC$. 

Further, arguing as in Case 1.1.2.a.2, we observe that relations (\ref{reln3}) hold for some $\alpha,\beta\in\RR$ and that
$$
[a_3,b_3]=\frac{2(\alpha^2+\beta^2)}{\mu}t_0.
$$  
Next, applying to $[a_r,b_r]$, $r=1,2$, the transformation from $\hbox{Ad}(H^0)$ that interchanges $a_1$ with $a_2$ and $b_1$ with $b_2$, we see that (\ref{formarbr}) simplifies to
\begin{equation}
\begin{array}{lll}
[a_1,b_1]&=&\sigma t_1+\tau t_2,\\

[a_2,b_2]&=&\tau t_1+\sigma t_2,
\end{array}\label{arbrverynew}
\end{equation}
for some $\sigma,\tau\in\RR$.

Let $R$ be the element of ${\frak h}$ for which the matrix ${\frak A}$ in representation (\ref{frakA}) is
$
\left(
\begin{array}{ll}
0 & i\\
i & 0
\end{array}
\right).
$
It now follows from (\ref{identities1}), (\ref{reln3}), (\ref{a1a2b1b2}), (\ref{a1a2h}), (\ref{frakAcd}), (\ref{arbrverynew}) and the Jacobi identity applied to the triples $\{a_1,b_1,a_2\}$, $\{a_1,b_1,a_3\}$, $\{a_1,b_1,R\}$ that
\begin{equation}
\begin{array}{ll}
\tau=0,& \displaystyle\nu=-\frac{\mu\sigma}{2},\\
\alpha\delta+\beta\gamma=0,& \displaystyle\beta\delta-\alpha\gamma=\frac{\mu\sigma}{2}.
\end{array}\label{relpar33}
\end{equation}

Assume first that $\alpha=0$, $\beta=0$. In this case relations (\ref{relpar33}) yield that ${\frak g}$ is isomorphic to the Lie algebra of the group $G_3(\CC^3)$ defined in Example (xii). Every connected Lie group ${\frak G}$ with Lie algebra isomorphic to that of $G_3(\CC^3)$ and having a compact 4-dimensional subgroup ${\frak K}$ such that ${\frak G}$ acts on ${\frak G}/{\frak K}$ effectively, is isomorphic to $G_3(\CC^3)$ in which case ${\frak K}$ is a maximal compact subgroup of ${\frak G}$. Therefore, the group $G$ is isomorphic to $G_3(\CC^3)$, and this isomorphism induces a diffeomorphism from $G/H$ onto $G_3(\CC^3)/K$, where $K$ is a maximal compact subgroup of $G_3(\CC^3)$. Next, $G_3(\CC^3)/K$ is $G_3(\CC^3)$-equivariantly diffeomorphic to $\CC^3$. The linearization of the action of $K$ on $\CC^3$ at its fixed point $p_K$ has exactly two non-trivial proper real invariant subspaces $L_1,L_2\subset T_{p_K}(\CC^3)$. Their dimensions are 2 and 4. As earlier, all $G_3(\CC^3)$-invariant complex structures on $\CC^3$ for which $L_1$ and $L_2$ are complex subspaces can be found using conditions (\ref{invarcomplstructure}). In this case the conditions lead to either the standard complex structure on $\CC^3$ or to its conjugate. The manifold arising in the latter case is holomorphically equivalent to $\CC^3$ by means of a map that preserves the group $G_3(\CC^3)$. Thus, we have obtained Example (xii).

Assume now that $\alpha^2+\beta^2>0$. In this case relations (\ref{relpar33}) imply that ${\frak g}$ is isomorphic to the Lie algebra of $Sp_2$. Every connected Lie group with Lie algebra isomorphic to that of $Sp_2$ and containing a subgroup isomorphic to $U_2$, is isomorphic to $Sp_2/\ZZ_2$. Therefore, $G$ is isomorphic to $Sp_2/\ZZ_2$. Under this isomorphism $H$ is mapped into a subgroup $K\subset Sp_2/\ZZ_2$ for which $K^0$ is isomorphic to $U_2$. Now, arguing as in Case 1.1.2.a.2, we observe that $K$ is connected. Indeed, the centralizer of any subgroup of $Sp_2/\ZZ_2$ isomorphic to $U_2$ coincides with the center of the subgroup. The isomorphism between $G$ and $Sp_2/\ZZ_2$ induces a diffeomorphism between $G/H$ and $\left(Sp_2/\ZZ_2\right)/K$. Since $K$ is connected, the manifold $\left(Sp_2/\ZZ_2\right)/K$ is $Sp_2/\ZZ_2$-equivariantly diffeomorphic to $\CC\PP^3$. The linearization of the action of $K$ on $\CC\PP^3$ at its fixed point $p_K$ has exactly two non-trivial proper real invariant subspaces $L_1,L_2\subset T_{p_K}(\CC\PP^3)$. Their dimensions are 2 and 4. As earlier, to determine all $Sp_2/\ZZ_2$-invariant complex structures on $\CC\PP^3$ for which $L_1$ and $L_2$ are complex subspaces, we use conditions (\ref{invarcomplstructure}). They lead either to the standard complex structure on $\CC\PP^3$ or to its conjugate. The manifold corresponding to the latter structure clearly is holomorphically equivalent to $\CC\PP^3$ by means of a map that preserves the group $Sp_2/\ZZ_2$. Thus, we have obtained Example (x).

{\bf Case 2.} {\it Suppose that $k_1=\pm k_2$.} Assume first that $k_1=1$, $k_2=-1$. Clearly, $M^0:=G/H^0$ is a finite-to-one cover of $G/H$, and we lift to $M^0$ the $G$-invariant complex structure induced on $G/H$ by that on $M$. The group $G$ acts on $M^0$ properly and effectively by holomorphic transformations with connected isotropy subgroups, and the linearization of the action of the subgroup $H^0$ on $M^0$ at the point $H^0\in M^0$ is given by $H^2_{1,-1}$ in some coordinates in $T_{H^0}(M^0)$. This linearization has exactly two invariant complex lines $l_1$, $l_2$ in $T_{H^0}(M^0)$. Let
\begin{equation}
{\cal L}_j:=\left\{dg(H^0)l_j, g\in G\right\},\quad j=1,2.\label{newdistr}
\end{equation}
Clearly, ${\cal L}_j$ is a $G$-invariant distribution of complex lines on $M^0$ for each $j$. We will now reason for the manifold $M^0$ in the way we did for the manifold $M$ in Case 1, with ${\cal S}_j$ constructed from ${\cal L}_j$, $j=1,2$. First of all, the argument at the beginning of Case 1.2.2 shows that (\ref{arbr}) holds and that ${\cal S}_1$ is a Lie subalgebra of ${\frak g}$. We will therefore repeat for the manifold $M^0$ the arguments that we applied to $M$ in Case 1.1.2. Transforming the commutators in ${\cal C}_1$ by the elements of the group $\hbox{Ad}(H^0)$, we obtain
\begin{equation}
\begin{array}{llr}
{\cal C}_1\subset{\frak h},\\
\vspace{-0.3cm}\\

[a_1,a_2]&=&-[b_1,b_2],\\
\vspace{-0.3cm}\\

[a_1,b_2]&=&[b_1,a_2].
\end{array}\label{reln21}
\end{equation}
The Jacobi identity applied to the triple $\{a_1,b_1,a_2\}$ together with (\ref{identities1}), (\ref{arbr}), (\ref{reln21}) now implies that all commutators in ${\cal C}_1$ are zero. 

As in Case 1.1.2, we are now led to the problem of determining certain $G_{1,-1}(\CC^2)$-invariant complex structures on $\CC^2$. Observe, however, that in contrast with Case 1.1.2, for a maximal compact subgroup $K\subset G_{1,-1}(\CC^2)$ the linearization of $K$ at its fixed point $p_K$ has many non-trivial proper real invariant subspaces in $T_{p_K}(\CC^2)$. All such subspaces have dimension 2. Therefore, for every pair of invariant subspaces $L_1,L_2\subset T_{p_K}(\CC^2)$ we look for all $G_{1,-1}(\CC^2)$-invariant complex structures on $\CC^2$, with respect to which $L_1$ and $L_2$ are complex lines. As before, we use conditions (\ref{invarcomplstructure}). It turns out that for every complex structure of this kind there is a biholomorphic map onto $\CC^2$ that either preserves the group $G_{1,-1}(\CC^2)$ or transforms it into the group $G_{1,1}(\CC^2)$. Such a biholomorphic map is given by a composition $f_1\circ f_2$, where $f_1$ is one of the three maps
$$
\begin{array}{lll}
z_1&\mapsto&\overline{z_1},\\
z_2&\mapsto&\overline{z_2},
\end{array}\quad
\begin{array}{lll}
z_1&\mapsto&z_1,\\
z_2&\mapsto&\overline{z_2},
\end{array}\quad
\begin{array}{lll}
z_1&\mapsto&\overline{z_1},\\
z_2&\mapsto&z_2,
\end{array}
$$
and $f_2$ is a map of the form
$$
\begin{array}{lll}
z_1&\mapsto&\alpha z_1+\beta\overline{z_2},\\
z_2&\mapsto&\gamma z_1+\delta\overline{z_2},
\end{array}
$$
for some $\alpha,\beta,\gamma,\delta\in\CC$. As a result, we obtain that $M^0$ is holomorphically equivalent to $\CC^2$ by means of a map that transforms $G$ (regarded as a subgroup of $\hbox{Aut}(M^0)$) into one of the two groups $G_{1,\pm 1}(\CC^2)$. Since every compact 1-dimensional subgroup of $G_{1,\pm 1}(\CC^2)$ is connected, it follows that $H^0=H$, hence $M$ is holomorphically equivalent to $M^0$ by means of a map that transforms $G$ into $G_{1,\pm 1}(\CC^2)$. Thus, we have obtained Example (vi) for $n=2$ and $k_1=1$, $k_2=\pm 1$.  

Suppose now that $k_1=k_2=1$ and consider the manifold $M^0:=G/H^0$, as above. The group $G$ acts on $M^0$ properly and effectively by holomorphic transformations with connected isotropy subgroups, and the linearization of the action of the group $H^0$ on $M^0$ at the point $H^0\in M^0$ is given by $H_{1,1}^2$ in some coordinates in $T_{H^0}(M^0)$. Choose two arbitrary complex lines $l_1$, $l_2$ in $T_{H^0}(M^0)$ and define ${\cal L}_j$ for $j=1,2$ by formula (\ref{newdistr}). Clearly, ${\cal L}_j$ is a $G$-invariant distribution of complex lines on $M^0$ for each $j$. We will now reason for the manifold $M^0$ in the way we did for the manifold $M$ in Case 1. Again, we see that (\ref{arbr}) holds and that ${\cal S}_1$ is a Lie subalgebra of ${\frak g}$. We will therefore repeat for the manifold $M^0$ the arguments that we applied to $M$ in Case 1.1.2. Transforming the commutators in ${\cal C}_1$ by the elements of the group $\hbox{Ad}(H^0)$, we obtain
\begin{equation}
\begin{array}{llr}
{\cal C}_1\subset{\frak h},\\
\vspace{-0.3cm}\\

[a_1,a_2]&=&[b_1,b_2],\\
\vspace{-0.3cm}\\

[a_1,b_2]&=&-[b_1,a_2].
\end{array}\label{reln22}
\end{equation}
The Jacobi identity applied to the triple $\{a_1,b_1,a_2\}$ together with (\ref{identities1}), (\ref{arbr}), (\ref{reln22}) implies that all commutators in ${\cal C}_1$ are zero. 

As in Case 1.1.2, we are now required to determine certain $G_{1,1}(\CC^2)$-invariant complex structures on $\CC^2$. Again, in contrast with Case 1.1.2, for a maximal compact subgroup $K\subset G_{1,1}(\CC^2)$ the linearization of $K$ at its fixed point $p_K$ has many non-trivial proper real invariant subspaces in $T_{p_K}(\CC^2)$, all of dimension 2. Again, using conditions (\ref{invarcomplstructure}), for every pair of invariant subspaces $L_1,L_2\subset T_{p_K}(\CC^2)$ we look for all $G_{1,1}(\CC^2)$-invariant complex structures on $\CC^2$, with respect to which $L_1$ and $L_2$ are complex lines. It turns out that the resulting complex structures are as follows: (a) the standard complex structure on $\CC^2$, (b) the complex structure conjugate to the standard one, (c) the complex structures obtained by conjugating the standard complex structure on one complex line in $\CC^2$ while preserving the standard structure on another complex line. Clearly, the manifolds arising in (b) and (c) are holomorphically equivalent to $\CC^2$ by means of a map that either preserves the group $G_{1,1}(\CC^2)$ or transforms it into $G_{1,-1}(\CC^2)$. Hence we obtain that $M^0$ is holomorphically equivalent to $\CC^2$ by means of a map that transforms $G$ into one of the two groups $G_{1,\pm 1}(\CC^2)$.  As above, this leads to Example (vi) for $n=2$ and $k_1=1$, $k_2=\pm 1$.

The proof of the theorem is now complete.\qed

{\obeylines
Department of Mathematics
The Australian National University
Canberra, ACT 0200
AUSTRALIA
E-mail: alexander.isaev@maths.anu.edu.au
\hbox{ \ \ }
Department of Complex Analysis
Steklov Mathematical Institute
8 Gubkina St.
Moscow GSP-1 119991
RUSSIA
E-mail: kruzhil@mi.ras.ru
}

\end{document}